\documentclass[preprint]{elsarticle}
\biboptions{sort&compress}
\usepackage{amssymb}
\usepackage{subfigure}
\usepackage{caption}
\usepackage{amsmath}
 \usepackage{hyperref}
\usepackage{float}
\usepackage[dvipsnames]{xcolor}
 
\usepackage{amsthm}
\usepackage{bm}
\usepackage{soul}
\usepackage{todonotes}
\newtheorem{thm}{Theorem}[section]

\newtheorem{lem}{Lemma}[section]

\theoremstyle{definition}
\newtheorem{defn}{Definition}[section]
\newtheorem{example}{Example}[section]

\theoremstyle{remark}
\newtheorem{rem}{Remark}[section]

\begin{document}

\begin{frontmatter}

\title{An oscillation-free summation-by-parts finite difference method for hyperbolic conservation laws}

\author[1]{Hanyu Wu}
\author[2]{Kenneth Duru}
\author[3]{Siyang Wang}

\address[1]{School of Mathematics and Statistics,
Central South University, Changsha, China}

\address[2]{Department of Mathematical Sciences,
University of Texas at El Paso, El Paso, Texas, USA}

\address[3]{Department of Mathematics and Mathematical Statistics,
Ume{\aa} University, Ume{\aa}, Sweden}

\begin{abstract}
 An oscillation-free summation-by-parts (SBP) finite difference method is proposed to solve hyperbolic conservation laws. High-order numerical schemes may generate spurious oscillations near discontinuities, which can reduce the robustness and accuracy of numerical simulations. To address this issue, motivated by the projection-based damping technique in a discontinuous Galerkin method (see J. Lu, Y. Liu, and C. W. Shu [SIAM J. Numer. Anal. 59 (2021), pp. 1299-1324]
 ), we introduce a damping term within a multi-block SBP finite difference framework to suppress oscillations while preserving high-order accuracy.
 The proposed method preserves conservation and satisfies discrete entropy stability for the semi-discrete scheme.
  Moreover, the damping term adjusts the numerical dissipation according to the local smoothness of the numerical solution. In smooth regions, the damping term remains sufficiently small so that the high-order accuracy is preserved. Near discontinuities and steep gradients, the damping term increases automatically to control spurious oscillations.
Numerical experiments are presented for the linear advection equation, Burgers’ equation and shallow water equations in one and two space dimensions. The results demonstrate that the proposed method effectively suppresses spurious oscillations while maintaining stability, conservation and robustness for both scalar conservation laws and hyperbolic systems of conservation laws.
\end{abstract}

\begin{keyword}
 Hyperbolic conservation laws \sep Oscillation-free \sep High-order finite difference methods  \sep Summation-by-parts

\end{keyword}

\end{frontmatter}

\section{Introduction}
Hyperbolic conservation laws arise in many areas of science and engineering, including wave propagation, transport flows, geophysical flows and many other applications \cite{LeVeque1992}. Accurate and stable numerical methods are essential for simulating complex physical problems governed by hyperbolic conservation laws. Although many numerical methods have been developed for hyperbolic conservation laws, problems with discontinuities or steep gradients remain challenging. This difficulty is particularly relevant for high-order numerical methods, which provide high accuracy for smooth solutions but may generate spurious oscillations near discontinuities. This motivates the development of effective oscillation-control techniques to improve the performance of high-order numerical methods for non-smooth solutions.

To control spurious oscillations in high-order numerical schemes, various approaches have been developed, including monotonic upstream-centered schemes for conservation laws (MUSCL), total variation diminishing (TVD), essentially non-oscillatory (ENO), weighted essentially non-oscillatory (WENO), and total variation bounded (TVB) methods  \cite{Harten1983,VanLeer1979,Harten1987,Shu1987, Liu1994}.
Traditional MUSCL and TVD schemes can achieve second-order accuracy in smooth regions and capture discontinuities more sharply than first-order schemes while preserving monotonicity \cite{Harten1983, VanLeer1979}. The ENO procedure automatically selects the smoothest stencil from several candidates according to the local solution. In this way, arbitrarily high-order accuracy can be achieved in smooth regions while maintaining an essentially non-oscillatory behavior near discontinuities \cite{Harten1987}. WENO schemes were subsequently developed as an improvement over ENO schemes and further refined in \cite{Shu1987,Jiang1996}. WENO schemes combine all candidate stencils through nonlinear weights determined by local smoothness indicators. This allows more stencil information to be used and improves the accuracy in smooth regions while retaining the non-oscillatory property near discontinuities. WENO schemes have been widely used in finite difference methods, and WENO-type reconstruction has also been employed in discontinuous Galerkin methods. Limiting techniques are also used to control spurious oscillations near discontinuities, such as TVB and WENO limiters for discontinuous Galerkin methods \cite{CockburnShu1989,QiuShu2005}. Another approach is artificial viscosity, which adds numerical dissipation in discontinuous regions or steep gradients \cite{Guermond2011,Hiltebrand2014,Stiernstrom2021}. 
 The amount of numerical dissipation must be carefully chosen, since excessive dissipation may overly smooth the solution, whereas insufficient dissipation may fail to suppress spurious oscillations.
More recently, projection-based damping terms have been introduced into discontinuous Galerkin methods to control high-order components near discontinuities while becoming very small in smooth regions \cite{Lu2021,Liu2022,Liu2024}. An advantage of this approach is that the damping coefficients are automatically determined from local interface jumps.

In addition to controlling spurious oscillations, ensuring the stability of high-order spatial discretizations is another important consideration for hyperbolic conservation laws. Summation-by-parts (SBP) finite difference operators were established in the early works \cite{KreissScherer1974,KreissScherer1977,Strand1994} and provided an effective approach for constructing high-order finite difference discretizations of hyperbolic problems. However, the stable treatment of boundary and interface conditions remained challenging for high-order finite difference methods. To address this difficulty, the simultaneous approximation term (SAT) technique was introduced to impose boundary conditions weakly, enabling stable discretizations of more complicated problems \cite{Carpenter1994}. 
Since then, SBP-SAT methods have developed into a systematic framework for constructing high-order and stable finite difference schemes \cite{SvardNordstrom2014,DelReyFernandez2014}. 
The key idea of SBP operators is to mimic integration by parts at the discrete level, allowing stability estimates to be established for semi-discrete schemes in a similar way to the continuous case \cite{Strand1994,SvardNordstrom2014,DelReyFernandez2014}. 
In particular, the continuous stability analysis may rely on the chain rule and product rule for nonlinear conservation laws, which is not satisfied exactly by SBP operators at the discrete level. 
To construct entropy-stable schemes for nonlinear conservation laws, skew-symmetric and entropy-conserving split formulations have been widely investigated \cite{Fisher2013,Gassner2016,Sjogreen2019,Nordstrom2022,Ricardo2024}.
Based on the development of upwind SBP operators \cite{Mattsson2017}, recent studies have further investigated upwind and dual-pairing SBP methods for nonlinear conservation laws \cite{Hew2025,Lundgren2020,Glaubitz2025,Ranocha2025,Stewart2026}. 
These developments have established SBP-SAT methods as a high-order, stable and flexible framework for  hyperbolic conservation laws.

Despite these favorable properties, high-order SBP-SAT schemes may  generate spurious oscillations when steep gradient and discontinuities exist. To address this issue, energy-stable WENO (ESWENO) methods have been developed by combining the SBP property with WENO methodology \cite{Yamaleev2009a,Yamaleev2009b,Fisher2011}, and have been extended to upwind SBP operators for linear advection problems \cite{JiangWang2024}. 
\ Nevertheless, developing high-order SBP-SAT methods that effectively suppress spurious oscillations while maintaining entropy stability for nonlinear conservation laws remains challenging. More recently, Lu, Liu, and Shu introduced adaptive projection-based damping into discontinuous Galerkin (DG) methods to suppress spurious oscillations \cite{Lu2021}. This approach was subsequently extended to systems of hyperbolic conservation laws \cite{Liu2022} and further developed to ensure entropy stability \cite{Liu2024}.
In this work, inspired by the projection-based damping approach, we develop an oscillation-free SBP-SAT finite difference method for hyperbolic conservation laws within the multi-block framework, applying split-form discretizations for nonlinear equations. In particular, our damping terms project the entropy variables into a lower dimensional polynomial space, minimizing oscillations at shocks and discontinuities  while ensuring high-order accuracy in smooth regions. In addition, our method is entropy consistent, that is, it reproduces correct entropy loss at shocks.
We prove that the proposed scheme satisfies discrete conservation and entropy stability. Elaborate numerical  experiments corroborate the theoretical analysis, and demonstrate that the proposed method effectively suppresses spurious oscillations and enhances the robustness of the underlying SBP-SAT method.

The remainder of this paper is organized as follows. Section~\ref{Sec.2} introduces the mathematical problems considered in this work. Section~\ref{Sec.3} presents the SBP operators and the general SBP-SAT semi-discrete formulation for one-dimensional systems of hyperbolic conservation laws. Section~\ref{Sec.4} introduces a projection-based damping term into the SBP-SAT framework for both one- and two-dimensional systems and establishes the discrete conservation and entropy stability properties of the oscillation-free SBP-SAT schemes. Numerical results are presented in Section~\ref{Sec.5}.

 \section{Mathematical models} \label{Sec.2}
 Hyperbolic conservation laws are a class of partial differential equations describing the propagation and evolution of conserved quantities such as mass, momentum and energy. The equations take  the general form
\begin{equation}\label{hcl}
	\partial_t\mathbf{U}(\mathbf{x} ,t) +\nabla\cdot    \mathbf{f}(\mathbf{U}(\mathbf{x} ,t)) = 0, \quad \mathbf{x} \in \Omega, \quad t \in[0, T],
\end{equation} 
where $\mathbf{U}: \Omega \times [0, T]\to \mathbb{R}^{m} $  is an $m$-dimensional vector of conserved variables with $m\geq 1$, $\mathbf{f}:\mathbb{R}^m\rightarrow \mathbb{R}^{m\times d}$ is the flux function for the hyperbolic conservation laws and $\mathbf{x}\in\mathbb{R}^d$ represents the spatial variable for $d\geq 1$. 
Eq.~\eqref{hcl} satisfies the following important properties:
\begin{enumerate}
    \item Conservation. Integrating Eq.~\eqref{hcl} over the domain $\Omega$ yields
\begin{equation*}
\dfrac{\mathrm{d}}{\mathrm{dt}}\int_{\Omega}\mathbf{U}(\mathbf{x},t) \mathrm{d}\mathbf{x}
= - \oint_{\partial\Omega}(\mathbf{f}(\mathbf{U}(\mathbf{x},t))\cdot\mathbf{n})\mathrm{dx},
\end{equation*}
where $\mathbf{n}$ denotes the outward unit normal vector on the boundary $\partial\Omega$. 
The above identity shows that the change in the total quantity within $\Omega$ is exactly balanced by the flux through the boundary $\partial\Omega$. In particular, if the boundary flux vanishes, the total quantity is conserved in time.

\item Entropy estimate. 
To derive entropy estimate for Eq.~\eqref{hcl}, we introduce $e:\mathbb{R}^m\to\mathbb{R}$ is the entropy function and $\mathbf{q}:\mathbb{R}^m\to\mathbb{R}^d$ is the corresponding entropy flux satisfying $\partial_{\mathbf{U}} \mathbf{q}(\mathbf{U})=
\mathbf{g}^\top\partial_{\mathbf{U}}\mathbf{f}(\mathbf{U})$, with the entropy variables defined by $\mathbf{g}=\partial_{\mathbf{U}} e(\mathbf{U})$.

When $\mathbf{U}$ is smooth, multiplying Eq.~\eqref{hcl} by $\mathbf{g}^\top$ and  integrating over the domain $\Omega$ yields
\[
\dfrac{\mathrm{d}}{\mathrm{dt}}\int_\Omega e(\mathbf{U})\mathrm{d}\mathbf{x} =
-\oint_{\partial\Omega}
(\mathbf{q}\cdot\mathbf{n})\mathrm{d}\mathbf{x}.
\]
Therefore, for non-smooth solutions, the evolution of the total entropy is governed by the entropy flux through the boundary. If the entropy flux through the boundary vanishes, the total entropy is also conserved in time.

When $\mathbf{U}$ is non-smooth, we enforce entropy inequality for the weak solution \cite{Tadmor2003, LeVeque2002}:
\[
\dfrac{\mathrm{d}}{\mathrm{dt}}\int_\Omega e(\mathbf{U})\mathrm{d}\mathbf{x} \le
-\oint_{\partial\Omega}
(\mathbf{q}\cdot\mathbf{n})\mathrm{d}\mathbf{x}.
\]
Therefore, for non-smooth solutions, the total entropy is non-increasing over time and is bounded above by its initial value.
\end{enumerate}

 High-order numerical schemes for hyperbolic conservation laws may generate spurious oscillations near discontinuities. Oscillation-free methods suppress oscillations by introducing suitable damping terms. The added damping is designed to remain small in smooth regions and take effect near discontinuities, while preserving good properties such as conservation and entropy stability\cite{Lu2021,Liu2024}.
In this work, we focus on three representative examples with increasing complexity: the linear advection equation as a model for linear scalar conservation laws, Burgers' equation as a model for nonlinear scalar conservation laws and the shallow water equations as a model for nonlinear systems of conservation laws. These three equations provide a suitable framework for demonstrating the construction of the oscillation-free SBP-SAT schemes, analyzing their properties such as discrete conservation and entropy stability, and validating their effectiveness through numerical experiments.
\begin{example}
When the flux is given by $\mathbf f(U)=\mathbf a U$, the general conservation law~\eqref{hcl} becomes the linear advection equation
\begin {equation} \label{adv}
\partial_t U + \nabla\cdot (\mathbf{a}U) = 0,
\end {equation}
where $U: \mathbb{R}^d\times \mathbb{R}^{+} \to \mathbb{R}$ and $\mathbf{a}=(a_1,a_2,\ldots,a_d) \in \mathbb{R}^d $ is a constant velocity vector.
\end{example}
\begin{example}\label{bur}
With the nonlinear flux $\mathbf f(U)=\left(
\frac{1}{2}U^2,\frac{1}{2}U^2,\ldots,\frac{1}{2}U^2
\right)^\top$,  the general conservation law~\eqref{hcl} becomes the inviscid Burgers' equation, where $U: \mathbb{R}^d\times \mathbb{R}^{+} \to \mathbb{R}$. 
\end{example}

\begin{example}\label{sw}
With the conservative variables $\mathbf U=(h,m,n)^T,$
where $m=hU$ and $n=hV$, the nonlinear flux is given by
$\mathbf f(\mathbf U)=\left(\mathbf f_1(\mathbf U),
\mathbf f_2(\mathbf U)\right)$, where
$\mathbf f_1(\mathbf U)=\left[m, \dfrac{m^2}{h}+\dfrac{1}{2}gh^2
,\dfrac{mn}{h}\right]^\top$ and
$\mathbf f_2(\mathbf U)=\left[n, \dfrac{mn}{h},\dfrac{n^2}{h}+\dfrac{1}{2}gh^2\right]^\top$.
Then, the general conservation law~\eqref{hcl} becomes the two-dimensional shallow water equations. Here $U: \mathbb{R}^2\times \mathbb{R}^{+} \to \mathbb{R}^3$ and $\mathbf{f}: \mathbb{R}^2\times \mathbb{R}^{+} \to \mathbb{R}^{3\times 2}$.

\end{example}

\section{Spatial discretization}\label{Sec.3}
In this section, we first introduce the SBP operators used for the spatial discretization and then present the semi-discrete formulations in the one-dimensional case. For a rectangular domain $\Omega$, it is straightforward to extend the results to multi-dimensional problems by tensor products.

\subsection{SBP operators}
Consider the 1D domain $\Omega=[a, b]\subset \mathbb{R}$ with $a< b$ discretized by the grid {$\{x_i\}$} with $x_i = a+i\Delta x, i=0,1,\ldots,N$, $ \Delta x = \dfrac{b-a}{N}$. 
To define the SBP operators, we first introduce the boundary projection operators
\begin{equation*}
\begin{aligned}
\boldsymbol{\delta}_i =[0,\ldots,0, 1,0,
\ldots,0]^\top, \quad
B=\boldsymbol{\delta}_{N} \boldsymbol{\delta}_{N}^{\top} - \boldsymbol{\delta}_0 \boldsymbol{\delta}_0^{\top},
\end{aligned}
\end{equation*}
where $\mathbf{\delta}_i$ denotes the $i$-th unit vector in $\mathbb{R}^{N+1}$ and $B\in\mathbb R^{(N+1)\times(N+1)}$.  We also define the pointwise product of $\mathbf f,\mathbf g\in\mathbb{R}^{N+1}$ as $\mathbf{g}^2=\mathbf g\circ\mathbf{g}$ and $\mathbf{f}\circ\mathbf{g}
:=\left(f(x_0)g(x_0), f(x_1)g(x_1),\ldots,
f(x_N)g(x_N)\right)^\top $.

The precise definitions of the SBP operators are given below.
\begin{defn}\label{df1}
The difference operator $D = H^{-1}\left(Q+\dfrac{B}{2}\right)$ is called a $p$th-order  SBP operator if it approximates $\partial/\partial x$ with a $p$th-order accurate interior stencil, where $H$ is diagonal positive definite and $Q+Q^\top=0$.
\end{defn}

\begin{defn}\label{df2}
The difference operators $D_{+} = H^{-1} (Q_{+} + \dfrac{B}{2})$ and $D_{-} = H^{-1} (Q_{-} + \dfrac{B}{2})$ are called $p$th-order upwind SBP operators if they approximate $\partial/\partial x$ with $p$th-order accurate interior stencils, where $H$ is diagonal positive definite, $Q_{+}+Q_{-}^\top=0$ and $\dfrac{Q_{+}+Q_{+}^\top}{2}=S$ is symmetric negative semi-definite.
\end{defn}
Here, $H$ in matrix form can be written as $H:= \Delta x \ \text{diag}(w_0,w_1,  \ldots, w_N)$, where $w_i>0$ for $i = 0,1,\cdots,N$.

\begin{rem}
In this work, we use SBP operators satisfying Definition~\ref{df1} that are constructed from the upwind SBP operators of  \cite{Mattsson2017}. Each pair of $D_-$ and $D_+$ uses upwind finite difference stencils of order $p_{in}$ in the interior and one-sided stencils of order $p_b$ on a few grid points near the boundaries. According to Definition~\ref{df2}, their average, denoted by $D_1: = (D_{+}+D_{-})/2$ satisfies the  SBP property in Definition~\ref{df1}. More precisely, $D_1$ employs a central finite difference stencil of order $p_1$ in the interior and one-sided stencil of order $p_b$ near the boundaries, see Table \ref{accuracy}. The boundary stencils differ from those of the SBP operators originally constructed in \cite{KreissScherer1974,Strand1994}. For stable discretization of first-order hyperbolic systems, we typically observe the convergence rate in space $p=p_b+1$ \cite{Gustafsson2008}. The operators in Table 1 include $D_1$ with equal interior order but different boundary order, and vice versa, which allows us to isolate the effect of the boundary closure in the numerical experiments. In Section 5, the order labels 2–6 refer to $p_{in}$.

\begin{table}[h!]
\centering
\begin{tabular}{c|ccccc}
\hline
Interior order $p_{in}$ of $D_+$, $D_-$  & 2 & 3 & 4 & 5 & 6 \\
\hline
Interior order $p_1$ of $D_1$  & 2 & 4 & 4 & 6 & 6 \\
Boundary order $p_b$ & 1 & 1 & 2 & 2 & 3 \\
Expected convergence rate $p$ & 2 & 2 & 3 & 3 & 4 \\
\hline
\end{tabular}
\caption{Order of accuracy }
\label{accuracy}
\end{table}
\end{rem}

The SBP operator $D_1$ and the norm matrix $H$ possess the following three properties.

\textbf{(P.1)} The diagonal $H$ defines a discrete inner product and a quadrature \cite{Hicken2013},
\[
\langle \mathbf f,\mathbf g\rangle_H=\mathbf f^T H\mathbf g
\quad \forall \, \mathbf f,\mathbf g \in\mathbb R^{N+1},\quad
 \langle \mathbf1,\mathbf f\rangle_H= \Delta x
\sum_{i=0}^{N}w_i f(x_i)\approx\int_{\Omega}f(x)\mathrm dx.
\]

\textbf{(P.2)} The linear operator $D_1\in \mathbb{R}^{(N+1)\times(N+1)}$ maps a grid function $\mathbf f\in\mathbb{R}^{N+1}$ to an approximation of its first derivative at the grid points. For a sufficiently smooth function $f$, the approximation has interior order $p_1$ and boundary order $p_b$, 
\[
(D_1\mathbf f)_i
=\left.\frac{df}{dx}\right|_{x=x_i}
+\begin{cases}
\mathcal O(\Delta x^{p_1}), & \text{at the grid points in the interior},\\
\mathcal O(\Delta x^{p_b}), & \text{at the grid points near the boundaries}.
\end{cases}
\]
In particular, $D_1\mathbf1=\mathbf0$. For simplicity of notation, throughout this paper, $\mathbf 1$ and $\mathbf 0$ denote the all-ones and all-zeros vectors of appropriate dimensions, respectively. 

\textbf{(P.3)} The operator $D_1$ with $H$ satisfies  $\langle D_1\mathbf f,\mathbf g\rangle_H
+\langle \mathbf f,D_1\mathbf g\rangle_H=f(x_N)g(x_N)-f(x_0)g(x_0), \quad \forall \,\mathbf f,\mathbf g\in\mathbb R^{N+1}$.

\subsection{Multi-block SBP-SAT method}\label{Sec.3.2}

We first consider the one-dimensional system with  periodic boundary conditions
\begin{equation}\label{general_1d}
    \partial_t \mathbf U +\partial_x \mathbf f(\mathbf U)=0, \quad x\in\Omega,\quad t\in[0,T],
\end{equation}
where $\mathbf U\in\mathbb R^m$ denotes the vector of conserved variables and
$\mathbf f(\mathbf U)$ is the corresponding flux function.

For linear problems, such as the advection equation, the standard flux form can be discretized directly by the SBP operator. For nonlinear conservation
laws, the entropy analyses generally involve the chain rule. However, the chain rule does not hold at the discrete level. Therefore, to mimic the continuous analysis in the discrete setting, the nonlinear equations are written in the following split form 
\begin{equation}\label{general_split}
    \partial_t\mathbf U +\mathbf F(\mathbf U,\mathbf f(\mathbf U),\partial_x)=0,
\end{equation}
where $\mathbf F(\mathbf U,\mathbf f(\mathbf U),\partial_x)$ is a split form of $\partial_x\mathbf f(\mathbf U)$. At the continuous level, the split form is equivalent to the original conservation law~\eqref{general_1d} for smooth solutions and leads to the same conservation and entropy estimates without using the chain rule, see for example \cite{Stewart2026}.
For the three model equations introduced in Section~\ref{Sec.2}, the corresponding split forms satisfy the following conservation and entropy identities:
\begin{align}\label{cons_conservation}
\int_\Omega\mathbf F(\mathbf U,\mathbf f(\mathbf U),\partial_x)\mathrm{dx}=\int_\Omega\partial_x\mathbf f(\mathbf U)\mathrm{dx}
=\oint_{\partial\Omega}\bigl(\mathbf f(\mathbf U)\cdot\mathbf n\bigr)\mathrm{dx},
\end{align}
\begin{align}\label{cons_entropyest}
\int_\Omega\mathbf g^\top\mathbf F(\mathbf U,\mathbf f(\mathbf U),\partial_x)\mathrm{dx}
=\int_\Omega \partial_x q(\mathbf U)\mathrm{dx}
=\oint_{\partial\Omega}\bigl( q(\mathbf U)\cdot\mathbf n\bigr)\mathrm{dx},\quad  \mathrm x \in \Omega \subset \mathbb R.
\end{align}

In the following, we present the SBP-SAT discretization of \eqref{general_split} under a two-subdomain formulation. Its generalization to multi-block configuration is straightforward. We present the standard formulation without any oscillation-free damping in this section, and explain in detail the damping term in the next section. 

To this end, we partition $\Omega=[a, b]\subset \mathbb{R}$ with $a< b$ into two subdomains $\Omega = \Omega_1\cup\Omega_2$ with $\Omega_1=[a^{(1)},b^{(1)}]$, $\Omega_2=[a^{(2)},b^{(2)}]$, $a^{(1)} = a$, $b^{(2)} = b$. In general, the subdomains may have different lengths and grid spacings. For notational simplicity, we consider two subdomains of equal length, such that $|\Omega_1|=|\Omega_2|=(b-a)/2$, with a uniform grid spacing $\Delta x$ within each element. Suitable interface conditions, such as continuity of conserved variables and fluxes, are imposed at the interface point $x=b^{(1)}=a^{(2)}$, where the degrees of freedom are duplicated. Each subdomain $\Omega_k$ is discretized by the grid points $ x_i^{(k)} = a^{(k)}+i\Delta x, i=0,1,\ldots,N$, $ \Delta x = (b^{(k)}-a^{(k)})/N$ for $k=1,2$. The numerical solution is denoted by  $\mathbf{u}_{\mu}=[\mathbf{u}^{(1)}_{\mu},\mathbf{u}^{(2)}_{\mu}]^\top\in \mathbb{R}^{2(N+1)}$ for $1\le \mu \le m$. The numerical solution for the global domain $\Omega$ is defined by $\mathbf{u}=[\mathbf{u}_1,\mathbf{u}_2,\cdots, \mathbf{u}_m]^\top\in \mathbb{R}^{2m(N+1)}$. For convenience, $\mathbf e^{(k)}=[e_0^{(k)},e_1^{(k)},\ldots,e_N^{(k)}]^T\in\mathbb R^{N+1}$ denotes the entropy vector on the $k$-th subdomain. 
 The corresponding semi-discrete entropy is given by
$ E_h(t)=\sum_{k=1}^{2}\left\langle \mathbf 1, \mathbf e^{(k)} \right\rangle_{H}
\approx \sum_{k=1}^{2}\int_{\Omega_k} e(\mathbf U^{(k)})\, \mathrm{d}\mathbf{x}$, where $\mathbf U^{(k)}\in\mathbb R^{m(N+1)}$ denotes the exact solution vector containing all $m$ conserved variables on the $k$-th subdomain.
The semi-discrete SBP-SAT scheme for \eqref{general_split} can be written as 
\begin{equation}\label{generalform_nondp}
   \partial_t \mathbf{u}= -\mathbf{F}(\mathbf{u},\mathbf{f}(\mathbf{u}),\widetilde{D}) + \widetilde{H}^{-1}\mathbf{SAT}(\mathbf{u}) +\widetilde H^{-1} \mathbf{T}(\mathbf{u}),
\end{equation}
where $\widetilde{D}= I_m \otimes I_2  \otimes D_1$ and $\widetilde{H}=I_m \otimes I_2  \otimes H$. On the right-hand side, the first term is a discrete approximation of $-\mathbf F(\mathbf U,\mathbf f(\mathbf U),\partial_x)$ using the SBP operators; the second term imposes weakly the boundary and interface conditions using the SAT method (without any dissipation). 
To reproduce the conservation and entropy identities~\eqref{cons_conservation} and ~\eqref{cons_entropyest} at the discrete level,  the SAT terms are constructed to cancel the corresponding boundary and interface contributions from the discrete split form, such that the following two identities hold:
\begin{equation}\label{SATi1}
-\sum_{k=1}^{2}\left\langle\mathbf 1,\mathbf F_{\mu}
\left(\mathbf u^{(k)},\mathbf f(\mathbf u^{(k)}),D_1\right)
\right\rangle_H+
\sum_{k=1}^{2}\mathbf 1^\top\mathbf{SAT}_{\mu}^{(k)}(\mathbf u)=0,\quad\forall \mu=1,2,\cdots,m,
\end{equation}
\begin{equation}\label{SATi2}
-\sum_{k=1}^{2}\sum_{\mu=1}^{m}\left\langle\mathbf g_{\mu}^{(k)},
\mathbf F_{\mu}\left(\mathbf u^{(k)},\mathbf f(\mathbf u^{(k)}),D_1
\right)\right\rangle_H+\sum_{k=1}^{2}\sum_{\mu=1}^{m}
\left(\mathbf g_{\mu}^{(k)}\right)^\top\mathbf{SAT}_{\mu}^{(k)}(\mathbf u)=0.
\end{equation}
Here, $\mathbf F_{\mu}
\left(\mathbf u^{(k)},\mathbf f(\mathbf u^{(k)}), D_1\right)$ denotes the $\mu$-th variable of the split-form spatial discretization on the $k$-th subdomain, and  $\mathbf g_{\mu}^{(k)}$ is the the $\mu$-th entropy variable on subdomain $k$.  The specific SAT formulations for the model problems are provided  in~\ref{app1}.

For solutions containing discontinuities or sharp gradients, introducing numerical dissipation through upwind interface penalties is a common approach to suppressing spurious oscillations in the SBP-SAT method. The last term in~\eqref{generalform_nondp} represents an upwind penalty on the jumps of the entropy variables at the interfaces and periodic boundaries, taking the form
$\widetilde{H}^{-1}T(\mathbf u)= \dfrac{1}{2}\widetilde{H}^{-1}(\boldsymbol{\alpha} \otimes (\mathcal B_I+\mathcal B_n) \mathbf g)$,
where $\boldsymbol{\alpha}=\dfrac{1}{2}\operatorname{diag}
\left(\left[\alpha_1^{(1)}+\alpha_1^{(2)},\,\alpha_2^{(1)}+\alpha_2^{(2)},\,\ldots,\alpha_m^{(1)}+\alpha_m^{(2)}\right]\right)$ is a diagonal positive semi-definite matrix with $\alpha_{\mu}^{(k)}\geq 0$,  $\mathbf g\in \mathbb{R}^{2m(N+1)}$ is the discrete entropy-variable vector. 
The interface and boundary matrices  
are \[\mathcal B_I=
\begin{pmatrix}
-\boldsymbol \delta_{N}\boldsymbol \delta_{N}^\top
&\boldsymbol \delta_{N}\boldsymbol \delta_0^\top\\
 \boldsymbol \delta_0\boldsymbol \delta_{N}^\top
&-\boldsymbol \delta_0\boldsymbol \delta_0^\top
\end{pmatrix}, \quad
\mathcal B_n=
\begin{pmatrix}
-\boldsymbol \delta_0\boldsymbol \delta_0^\top
& \boldsymbol \delta_0\boldsymbol \delta_{N}^\top\\
\boldsymbol \delta_{N}\boldsymbol \delta_0^\top
&-\boldsymbol \delta_{N}\boldsymbol \delta_{N}^\top
\end{pmatrix}. \]
Both $\mathcal B_I$ and $\mathcal B_n$ are negative semi-definite. They also satisfy $\mathbf 1^\top\mathcal B_I=\mathbf 1^\top\mathcal B_n=\mathbf 0^\top$. 
It is worth noting that, for smooth solutions
satisfying the interface continuity and periodic
boundary conditions, the jumps of the entropy
variables vanish, yielding
\begin{equation*}
(\mathcal B_I+\mathcal B_n)\mathbf g_\mu
=\mathbf 0,\qquad \mu=1,\ldots,m.
\end{equation*}
Thus, the upwind interface dissipation vanishes for smooth exact solutions, and the semi-discrete scheme~\eqref{generalform_nondp} remains a consistent approximation of the original conservation law~\eqref{general_1d}.
Moreover, such interface dissipation may not be sufficient to eliminate spurious oscillations near discontinuities, even when the scheme satisfies a semi-discrete energy or entropy stability. This limitation is further illustrated by the numerical results in Section~\ref{Sec.5}.

Next, we present two theorems for the conservation and stability properties of \eqref{generalform_nondp}.
\begin{thm}\label{thm_cons_nodp}
The semi-discrete SBP-SAT discretization \eqref{generalform_nondp} satisfies discrete conservation for each conserved variable $\mu=1,\ldots,m$, that is, 
\[
\frac{\mathrm d}{\mathrm dt}
\sum_{k=1}^{2}\langle \mathbf 1,\mathbf u_{\mu}^{(k)}\rangle_H=0,
\]
if the SAT terms are constructed such that \eqref{SATi1} is satisfied. 
\end{thm}
\begin{proof}
For the $\mu$-th conserved variable, taking the discrete inner product of the corresponding semi-discrete equation~\eqref{generalform_nondp} with the constant grid function $\mathbf 1$ and summing over the two subdomains gives
\[\small
\begin{aligned}
\frac{\mathrm d}{\mathrm dt}\sum_{k=1}^{2}\langle\mathbf 1,\mathbf u_{\mu}^{(k)}\rangle_H=
\sum_{k=1}^{2}\left(-\left\langle\mathbf 1,\mathbf F_{\mu}
\left(\mathbf u^{(k)},\mathbf f(\mathbf u^{(k)}),
D_1\right)\right\rangle_H+
\mathbf 1^\top
\mathbf{SAT}_{\mu}(\mathbf u^{(k)})\right)
+\mathbf 1^\top
\mathbf T_{\mu}(\mathbf u).
\end{aligned}
\]
On the right-hand side, the first two terms cancel because of \eqref{SATi1}. For the upwind interface term, the $\mu$-th component has the form
$\mathbf T_{\mu}(\mathbf u)=\frac{\alpha_{\mu}}{2}\left(
\mathcal B_I+\mathcal B_n\right)\mathbf g_{\mu}$, where $\mathbf g_{\mu}=[\mathbf g_{\mu}^{(1)},\mathbf g_{\mu}^{(2)}]^\top\in\mathbb R^{2(N+1)}$ denotes the discrete vector of the $\mu$-th entropy variable. Since $\mathbf 1^\top\mathcal B_I=\mathbf 1^\top\mathcal B_n=\mathbf 0^\top$, it follows that
\[
\begin{aligned}
\mathbf 1^\top \mathbf T_{\mu}(\mathbf u)=\frac{\alpha_{\mu}}{2}\mathbf 1^\top \left(\mathcal B_I+\mathcal B_n\right)\mathbf g_{\mu}=0.
\end{aligned}
\]
Combining the above results gives
\[
\frac{\mathrm d}{\mathrm dt}\sum_{k=1}^{2}\langle
\mathbf 1,\mathbf u_{\mu}^{(k)}\rangle_H=0,\quad \mu=1,\ldots,m.
\]
Thus, the two-subdomain semi-discrete SBP-SAT formulation  is globally conservative for each conserved variable.
\end{proof}

Next, we consider the discrete entropy stability for  Eq.~\eqref{generalform_nondp}.
\begin{thm}\label{thm_L2_nodp}
The semi-discrete SBP-SAT discretization \eqref{generalform_nondp} satisfies a discrete entropy estimate,
\[
\frac{\mathrm d}{\mathrm dt}E_h(t)=
\sum_{k=1}^{2}\sum_{\mu=1}^{m}\left\langle\mathbf g_{\mu}^{(k)},
\partial_t\mathbf u_{\mu}^{(k)}\right\rangle_H\le 0,
\]
if the SAT terms are constructed such that \eqref{SATi2} is satisfied. 
\end{thm}

\begin{proof}
Taking the discrete inner product of the semi-discrete equations with the corresponding entropy variables and summing over all $m$ conserved variables and the two subdomains gives
\[
\begin{aligned}
\sum_{k=1}^{2}\sum_{\mu=1}^{m}\left\langle\mathbf g_{\mu}^{(k)},
\partial_t\mathbf u_{\mu}^{(k)}\right\rangle_H=&
-\sum_{k=1}^{2}\sum_{\mu=1}^{m}
\left\langle\mathbf g_{\mu}^{(k)},\mathbf F_{\mu}\left(
\mathbf u^{(k)},\mathbf f(\mathbf u^{(k)}),D_1\right)\right\rangle_H\\
&+\sum_{k=1}^{2}\sum_{\mu=1}^{m}
\left(\mathbf g_{\mu}^{(k)}\right)^\top\mathbf{SAT}_{\mu}^{(k)}(\mathbf u)+\sum_{\mu=1}^{m}\left(\mathbf g_{\mu}\right)^\top\mathbf T_{\mu}(\mathbf u).
\end{aligned}
\]
The first two terms on the right-hand side cancel because of \eqref{SATi2}. 
For the upwind interface term, the $\mu$-th component is $
\mathbf T_{\mu}(\mathbf u)=\frac{\alpha_{\mu}}{2}
\left(\mathcal B_I+\mathcal B_n\right)\mathbf g_{\mu}$, we have
\[
\sum_{\mu=1}^{m} \left(\mathbf g_{\mu}\right)^\top
\mathbf T_{\mu}(\mathbf u)=\frac12 \sum_{\mu=1}^{m}\alpha_{\mu}\mathbf g_{\mu}^\top\left(\mathcal B_I+\mathcal B_n\right)\mathbf g_{\mu}.
\]
Since $\alpha_{\mu}\ge0$ and $\mathcal B_I+\mathcal B_n$ is negative semi-definite, $\frac12\sum_{\mu=1}^{m}\alpha_{\mu}\mathbf g_{\mu}^\top
\left(\mathcal B_I+\mathcal B_n\right)\mathbf g_{\mu}\le0$.

Combining the above results yields
\[
\frac{\mathrm d}{\mathrm dt}E_h(t)=
\sum_{k=1}^{2}\sum_{\mu=1}^{m}
\langle\mathbf g_{\mu}^{(k)},
\partial_t\mathbf u_{\mu}^{(k)}
\rangle_H
\le0.
\]

Thus, the two-subdomain semi-discrete SBP-SAT formulation is entropy-stable.
\end{proof}

The SAT terms are in general not unique, and a variety of SAT formulations have been developed for  conservation laws \cite{Stewart2026,
Stiernstrom2021, Lundgren2020,Ranocha2017h,Ranocha2017,Carpenter1999, Berg2012}. 
    In the present study, the SAT terms in the Appendix should be regarded as one possible choice for the model problems. In particular, the SAT terms for Burgers' equation and the shallow water equations are chosen following \cite{Ranocha2017h} and \cite{Ranocha2017}, respectively. In addition, the SAT for upwind interface  is chosen from \cite{Stewart2026}.

\section{Oscillation-Free SBP-SAT finite difference methods}\label{Sec.4}
In the previous section, we presented the SBP-SAT semi-discrete formulations for a general one-dimensional system of hyperbolic conservation laws. These existing discretizations provide stable and accurate approximations for smooth solutions. However, the SBP-SAT scheme~\eqref{generalform_nondp} may still produce spurious numerical oscillations near discontinuities or sharp gradients. To suppress these oscillations, we introduce an additional damping term into the semi-discrete scheme~\eqref{generalform_nondp}, while preserving the desirable properties of the original discretization. The construction of the damping term is inspired by the projection-based damping method developed within the DG framework~\cite{Lu2021} where it was applied to the conservative form of nonlinear conservation laws. In this section, we extend the projection-based damping to the SBP-SAT method, applied to the split form and prove entropy stability. 

The remainder of this section is organized as follows.
In Section~\ref{4_1}, we first consider one-dimensional systems of conservation laws with $m$ conserved variables. The construction of the damping term is presented in detail, and several properties are stated, including discrete conservation and entropy stability.
In Section~\ref{4_2}, the two-dimensional construction is obtained from the one-dimensional formulation through tensor products for a rectangular domain $\Omega$. The corresponding oscillation-free SBP-SAT semi-discrete formulation is then presented, together with its conservation and stability properties.

\subsection{One-dimensional system}\label{4_1}
In this section, we consider a general one-dimensional system of hyperbolic conservation laws with $m$ conserved variables. We first construct the projection-based damping term for each conserved variable and then incorporate it into the SBP-SAT semi-discrete formulation~\eqref{generalform_nondp}. The discrete conservation and entropy stability properties of the oscillation-free SBP-SAT scheme are also established.
\subsubsection{Construction of the damping term in the oscillation-free SBP method}
For the two-subdomain problem, based on the SBP-SAT formulation~\eqref{generalform_nondp}, the semi-discrete oscillation-free SBP-SAT scheme can be written in the following general form
\begin{equation}\label{generalform_system}
\partial_t \mathbf{u}=-\mathbf{F}(\widetilde D,\mathbf{u})
+\widetilde{H}^{-1}\mathbf{SAT}(\mathbf{u})+\widetilde H^{-1}\mathbf{T}(\mathbf{u})+\mathbf{OF}.
\end{equation}
Here, $\mathbf u=\left[\mathbf u_1,\mathbf u_2,\ldots,
\mathbf u_m\right]^\top \in\mathbb R^{2m(N+1)}$ is the numerical solution, where $\mathbf u_{\mu}= \left[\mathbf u_{\mu}^{(1)},\mathbf u_{\mu}^{(2)}\right]^\top \in\mathbb R^{2(N+1)}$ for $ \mu=1,\ldots,m$ and $\mathbf u_{\mu}^{(k)}\in \mathbb R^{N+1}$ is the numerical solution of the $\mu$-th conserved variable on the $k$-th subdomain.
$\mathbf{OF}_{\mu}=\left[\mathbf {OF}_{\mu}^{(1)},\mathbf {OF}_{\mu}^{(2)}
\right]^\top$ is the damping term for the $\mu$-th conserved variable.

For each conserved variable per subdomain, the damping term takes the form
\[
\mathbf{OF}_{\mu}^{(k)}=-\sum_{l=0}^{s}
\dfrac{\sigma_{\mu}^{(k,l)}}{|\Omega_k|}\left(
\mathbf g_{\mu}^{(k)}-P^{\,l-1}\mathbf g_{\mu}^{(k)}
\right),\quad \mu=1,\ldots,m, \quad k=1,2.
\]
where $P^{\,l}$, $l\ge0$, denotes the discrete $L^2$ projection that is computed using the standard Legendre polynomials of degree no greater than $l$, $s \ll N$ is the highest degree of the Legendre polynomials used in the discrete $L^2$ projection, $|\Omega_k|$ is the width of the subdomain, and $\sigma_{\mu}^{(k,l)}\ge0$ is the scalar damping coefficient corresponding to the $\mu$-th conserved variable on the $k$-th subdomain.

It is worth noting that when $m=1$, the formulation~\eqref{generalform_system} reduces to the one-dimensional scalar case. For the linear advection and Burgers' equations considered in this paper, the energy and entropy functions are the same, and the corresponding entropy variable is $g=u$. Therefore, the damping term constructed from the entropy variable is equivalent to that constructed directly from the numerical solution.

For the $\mu$-th conserved variable on the subdomain $k$, the oscillation-free damping term is constructed component-wise. The construction of the corresponding discrete $L^2$ projection proceeds as follows. The grid points $\{x_i^{(k)}\}_{i=0}^{N}$ on the subdomain $\Omega_k=[a^{(k)},b^{(k)}]$, $k=1,2$, are mapped to the reference interval $[-1,1]$ through
\[
\xi_i=\dfrac{2(x_i^{(k)}-a^{(k)})}{b^{(k)}-a^{(k)}}-1,\quad
i=0,\dots,N.
\]
Let $\boldsymbol{\xi}=\left[\xi_0,\xi_1,\ldots,\xi_N\right]^\top$
denote the reference vector and $L_n$ be the standard Legendre polynomial of degree $n$ on $[-1,1]$. By calculating the values of $L_n$ at all reference points $\boldsymbol{\xi}$, we have
\[
L_n(\boldsymbol{\xi})=\left[L_n(\xi_0),L_n(\xi_1),\ldots,L_n(\xi_N)\right]^\top,
\]
which is referred to as the corresponding discrete Legendre basis vector.

The basis vectors satisfy
\[
L_{n+1}(\boldsymbol{\xi})=\dfrac{
(2n+1)\boldsymbol{\xi}\circ L_n(\boldsymbol{\xi})-nL_{n-1}(\boldsymbol{\xi})}{n+1}, \quad L_0(\boldsymbol{\xi})=\mathbf 1,
\quad L_1(\boldsymbol{\xi})=\boldsymbol{\xi}, \quad n\ge1.
\]
The discrete Legendre basis matrix is then defined by
\[
\Phi=\left[L_0(\boldsymbol{\xi}),L_1(\boldsymbol{\xi}),
\ldots,L_s(\boldsymbol{\xi})\right]\in\mathbb R^{(N+1)\times(s+1)}.
\]
For each \(l=0,\ldots,s\), let \(\Phi_l=\left[L_0(\boldsymbol{\xi}),L_1(\boldsymbol{\xi}),\ldots,L_l(\boldsymbol{\xi})\right]\in\mathbb R^{(N+1)\times(l+1)}\) denote the submatrix consisting of the first \(l+1\) columns of \(\Phi\). For the \(\mu\)-th conserved variable on the \(k\)-th subdomain, the discrete \(L^2\) projection \(P^l\mathbf g_\mu^{(k)}\) is defined by the following orthogonality condition:
\begin{align}\label{projection_system}
\left\langle \mathbf g_\mu^{(k)}-\Phi_l\mathbf c_{\mu,l}^{(k)},\,\Phi_l\mathbf v\right\rangle_H=0,
\quad \forall\,\mathbf v\in\mathbb R^{l+1}.
\end{align}
It follows that $\mathbf c_{\mu,l}^{(k)}=(\Phi_l^\top H\Phi_l)^{-1}\Phi_l^\top H\mathbf g_\mu^{(k)}$, where $\mathbf c_{\mu,l}^{(k)}\in \mathbb R^{l+1}$ denotes the coefficient vector of the discrete $L^2$ projection.
Therefore, the discrete \(L^2\) projection is given by
\[
P^l\mathbf g_\mu^{(k)}=\Phi_l\mathbf c_{\mu,l}^{(k)}.
\]
In particular, we define \(P^{-1}=P^0\).

In order to ensure that the damping coefficient $\sigma_{\mu}^{(k,l)}$ remains small in smooth regions and becomes large near discontinuities, we take
\[
\sigma_{\mu}^{(k,l)}
=C_l (\Delta x)^l\sqrt{[\mkern-2.8mu[
\mathcal D^l \mathbf z_{\mu}^{(k)}
]\mkern-2.8mu]_{L}^2 + [\mkern-2.8mu[
\mathcal D^l \mathbf z_{\mu}^{(k)}
]\mkern-2.8mu]_{R}^2 },\quad C_l=\dfrac{2(2l+1)}{(2s-1)l!},
\]
where 
$\mathbf z_{\mu}^{(k)}$ denotes an auxiliary variable constructed from the numerical solution to determine the adaptive damping coefficient for the $\mu$-th conserved variable on the $k$-th subdomain. Its physical dimensions are chosen to ensure dimensional consistency of the corresponding damping term. Its specific form depends on the governing equations. In particular, $\mathbf z_{\mu}^{(k)}=\mathbf u_{\mu}^{(k)}$ for the advection and Burgers’ equations, while the auxiliary variable $\mathbf z_{\mu}^{(k)}$ for the shallow water equations are given in~\eqref{z_1D}. 
$[\mkern-2.8mu[
\mathcal D^l \mathbf z_{\mu}^{(k)}
]\mkern-2.8mu]_{L}$ denotes the jump of the $l$-th derivative for $\mathbf z_{\mu}^{(k)}$ at the left endpoint and $[\mkern-2.8mu[
\mathcal D^l \mathbf z_{\mu}^{(k)}
]\mkern-2.8mu]_{R}$ represents the jump of the $l$-th derivative for $\mathbf z_{\mu}^{(k)}$ at the right endpoint.

At the right endpoint of the $k$-th subdomain, the right jump of the $l$-th derivative is defined by
\[
[\mkern-2.8mu[
\mathcal D^l \mathbf z_{\mu}^{(k)}
]\mkern-2.8mu]_{R}=\left(\mathbf d_{\mathcal I_L}^{(k+1,l)}
\right)^\top\mathbf z_{\mu}^{(k+1)}-
\left(\mathbf d_{\mathcal I_R}^{(k,l)}\right)^\top\mathbf z_{\mu}^{(k)}.
\]
For the left endpoint of the $k$-th subdomain,
\[
[\mkern-2.8mu[
\mathcal D^l \mathbf z_{\mu}^{(k)}
]\mkern-2.8mu]_{L}
=\left(\mathbf d_{\mathcal I_L}^{(k,l)}
\right)^\top\mathbf z_{\mu}^{(k)}-
\left(\mathbf d_{\mathcal I_R}^{(k-1,l)}
\right)^\top\mathbf z_{\mu}^{(k-1)}.
\]
Under periodic boundary conditions, we take $k+1=1$ when $k=2$ and $k-1=2$ when $k=1$.
The vectors $\mathbf d_{\mathcal I_L}^{(k,l)}$ and $\mathbf d_{\mathcal I_R}^{(k,l)}$ denote the coefficient vectors of $p$-th order one-sided finite difference approximations to the $l$-th derivative, where $p$ is the expected convergence order of the SBP operator $D_1$ listed in Table~\ref{accuracy}. These coefficients are determined using Taylor expansions following \cite{Gustafsson2008}.

To investigate whether the damping term affects the accuracy of the
SBP-SAT discretization in smooth regions, we estimate the damping
term using a smooth exact solution $\mathbf U$.
For the $\mu$-th conserved variable on each subdomain $\Omega_k$,
we introduce the projection residual
$\mathbf r_{\mu}^{(k,l)}
=\mathbf g_{\mu}^{(k)}-P^{\,l-1}\mathbf g_{\mu}^{(k)}$ and the discrete space
$ V_{l-1}=\left\{\sum_{i=0}^{l-1} v_i L_i(\boldsymbol{\xi}) :v_i\in\mathbb{R}, V_{-1}=V_0\right\}$.
By the orthogonality of the discrete $L^2$ projection,
\[\langle \mathbf g_{\mu}^{(k)} - P^{\,l-1}\mathbf g_{\mu}^{(k)},\mathbf v\rangle_H =0,\quad\forall\, \mathbf v\in V_{l-1}.\]
Throughout the paper, $\Delta x=|\Omega_k|/N$ denotes the
uniform grid spacing on $\Omega_k$, where
$|\Omega_k|=|\Omega|/K$ is the subdomain width.
Furthermore,
$g_{\mu,i}^{(k)}:=g_{\mu}(\mathbf U(x_i^{(k)}))$
denotes the value of the $\mu$-th entropy variable at
the grid point $x_i^{(k)}$ in $\Omega_k$.
The symbol $C$ denotes a generic positive constant,
which may take different values at different occurrences.

\begin{lem}\label{lemma_OF_acc}
Assume that the exact solution $\mathbf U(\cdot,t)$ is
sufficiently smooth such that
$g_{\mu}(\mathbf U),z_{\mu}(\mathbf U)
\in C^{s+p}(\Omega)$ for $\mu=1,\ldots,m$.
Then the damping term evaluated at the exact solution,
\[
\mathbf{OF}_{\mu}^{(k)}(\mathbf U)
=-\sum_{l=0}^{s}
\dfrac{\sigma_{\mu}^{(k,l)}(\mathbf U)}{|\Omega_k|}
\mathbf r_{\mu}^{(k,l)}(\mathbf U),
\]
satisfies
\[
\left\|\mathbf{OF}_{\mu}^{(k)}(\mathbf U)\right\|_H
\le C|\Omega_k|^{p+\frac32},
\qquad \mu=1,\ldots,m,\quad k=1,2,
\]
where $C$ is independent of $N$ and $|\Omega_k|$.
\end{lem}

The proof of Lemma~\ref{lemma_OF_acc}
is presented in Appendix~\ref{app_OF_acc}.

Consequently, for $m=1$, summing over the two subdomains yields $\left|\mathbf{OF}_{\mu}(\mathbf U)\right|_{\widetilde{H}}\le C|\Omega_k|^{p+1}$.

\subsubsection{Discrete conservation and entropy stability}
We now establish the discrete conservation and entropy stability for the oscillation-free SBP-SAT scheme~\eqref{generalform_system}. 

The first theorem proves the oscillation-free SBP-SAT method~\eqref{generalform_system} is conservative.
\begin{thm}\label{thm_cons}
Consider the general semi-discrete oscillation-free SBP-SAT
formulation~\eqref{generalform_system} with periodic boundary conditions. Then the semi-discrete scheme is conservative for each conserved variable, that is,
\begin{equation*}
\dfrac{\mathrm d}{\mathrm dt}
\sum_{k=1}^{2}
\left\langle\mathbf 1,\mathbf u_{\mu}^{(k)}
\right\rangle_H=0,
\qquad \mu=1,\ldots,m.
\end{equation*}
\end{thm}
\begin{proof}
For the conservation property, taking the discrete inner product of
the $\mu$-th component of Eq.~\eqref{generalform_system} with $\mathbf 1 $ gives
\begin{align*}
\dfrac{\mathrm d}{\mathrm dt}
\sum_{k=1}^{2}
\left\langle \mathbf 1,\mathbf u_{\mu}^{(k)}\right\rangle_H
=&-\sum_{k=1}^{2}
\left\langle
\mathbf 1,\mathbf F_{\mu}\left(\mathbf u^{(k)},\mathbf f(\mathbf u^{(k)}), D_1\right)
\right\rangle_H
+\sum_{k=1}^{2}
\mathbf 1^\top \mathbf{SAT}_{\mu}(\mathbf u^{(k)})\\
&+\mathbf 1^\top \mathbf{T}_{\mu}(\mathbf u)
+\sum_{k=1}^{2}
\left\langle \mathbf 1,\mathbf{OF}_{\mu}^{(k)}\right\rangle_H .
\end{align*}
The first three terms on the right-hand side correspond to the underlying SBP-SAT discretization~\eqref{generalform_nondp}. By Theorem~\ref{thm_cons_nodp} in Section~\ref{Sec.3.2}, their contribution is conservative. Therefore,
\[
-\sum_{k=1}^{2}
\left\langle 1,
\mathbf F_{\mu}\left(\mathbf u^{(k)},\mathbf f(\mathbf u^{(k)}), D_1\right)
\right\rangle_H
+\sum_{k=1}^{2}
\mathbf 1^\top \mathbf{SAT}_{\mu}(\mathbf u^{(k)})
+\mathbf 1^\top \mathbf{T}_{\mu}(\mathbf u)=0.
\]
For the last term, by the discrete orthogonality condition~\eqref{projection_system}, we have $\langle \mathbf 1,\mathbf r_{\mu}^{(k,l)}\rangle_H=0$.
Therefore,
\[
\begin{aligned}
\left\langle \mathbf 1,\mathbf{OF}_{\mu}^{(k)}\right\rangle_H
&=-\sum_{l=0}^{s}
\dfrac{\sigma_{\mu}^{(k,l)}}{|\Omega_k|}
\left\langle \mathbf 1,\mathbf r_{\mu}^{(k,l)}\right\rangle_H=0 .
\end{aligned}
\]
Summing over the two subdomains yields
\[
\sum_{k=1}^{2}
\left\langle \mathbf 1,\mathbf{OF}_{\mu}^{(k)}\right\rangle_H=0.
\]
Combining the above results yields
\[
\dfrac{\mathrm d}{\mathrm dt}
\sum_{k=1}^{2}
\left\langle \mathbf 1,\mathbf u_{\mu}^{(k)}\right\rangle_H=0,
\qquad \mu=1,\ldots,m.
\]
Thus, the general semi-discrete formulation~\eqref{generalform_system} preserves conservation for each conserved variable.
\end{proof}

Next, we consider the discrete entropy stability for the Eq.~\eqref{generalform_system}.
\begin{thm}\label{thm_L2}
Consider the general semi-discrete oscillation-free SBP-SAT
formulation~\eqref{generalform_system} with periodic boundary conditions.
Then the semi-discrete approximation is entropy-stable, that is,
\begin{equation}\label{general_L2}
\dfrac{\mathrm dE_h}{\mathrm dt}\leq 0,
\end{equation}
where $E_h$ denotes the discrete entropy defined in Section~\ref{Sec.3.2}.
\end{thm}

\begin{proof}
Taking the discrete inner product of each component of Eq.~\eqref{generalform_system} with the corresponding entropy variable gives
\[
\begin{aligned}
\dfrac{\mathrm dE_h}{\mathrm dt}
=&-\sum_{k=1}^{2}\sum_{\mu=1}^{m}
\left\langle
\mathbf g_{\mu}^{(k)},
\mathbf F_{\mu}\left(\mathbf u^{(k)},\mathbf f(\mathbf u^{(k)}), D_1\right)
\right\rangle_H
+\sum_{k=1}^{2}\sum_{\mu=1}^{m}
(\mathbf g_{\mu}^{(k)})^\top
\mathbf{SAT}_{\mu}(\mathbf u^{(k)})\\
&+\sum_{\mu=1}^{m}
(\mathbf g_{\mu})^\top
\mathbf T_{\mu}(\mathbf u)
+\sum_{k=1}^{2}\sum_{\mu=1}^{m}
\left\langle
\mathbf g_{\mu}^{(k)},\mathbf{OF}_{\mu}^{(k)}
\right\rangle_H .
\end{aligned}
\]
The first three terms on the right-hand side correspond to the underlying SBP-SAT discretization~\eqref{generalform_nondp}. By Theorem~\ref{thm_L2_nodp} in Section~\ref{Sec.3.2}, we have
\[\small
\begin{aligned}
-\sum_{k=1}^{2}\sum_{\mu=1}^{m}
\left(\left\langle
\mathbf g_{\mu}^{(k)},
\mathbf F_{\mu}\left(\mathbf u^{(k)},\mathbf f(\mathbf u^{(k)}), D_1\right)
\right\rangle_H
+(\mathbf g_{\mu}^{(k)})^\top
\mathbf{SAT}_{\mu}(\mathbf u^{(k)})\right)
+\sum_{\mu=1}^{m}(\mathbf g_{\mu})^\top
\mathbf T_{\mu}(\mathbf u)
\leq0.
\end{aligned}
\]
For the damping term, using $\mathbf g_{\mu}^{(k)}
=P^{\,l-1}\mathbf g_{\mu}^{(k)}+\mathbf r_{\mu}^{(k,l)}$ and the orthogonality of the projection residual gives
\[
\begin{aligned}
\left\langle
\mathbf g_{\mu}^{(k)},\mathbf r_{\mu}^{(k,l)}
\right\rangle_H=\langle P^{\,l-1}\mathbf g_{\mu}^{(k)},\mathbf r_{\mu}^{(k,l)}\rangle_H+\left\|\mathbf r_{\mu}^{(k,l)}
\right\|_H^2
=\left\|\mathbf r_{\mu}^{(k,l)}\right\|_H^2.
\end{aligned}
\]
Consequently,
\[
\begin{aligned}
\sum_{k=1}^{2}\sum_{\mu=1}^{m}
\left\langle
\mathbf g_{\mu}^{(k)},\mathbf{OF}_{\mu}^{(k)}
\right\rangle_H
=-\sum_{k=1}^{2}\sum_{\mu=1}^{m}\sum_{l=0}^{s}
\dfrac{\sigma_{\mu}^{(k,l)}}{|\Omega_k|}
\left\|\mathbf r_{\mu}^{(k,l)}\right\|_H^2\leq0.
\end{aligned}
\]
Combining the above estimates yields
\[
\dfrac{\mathrm dE_h}{\mathrm dt}\leq0.
\]
Thus, the general semi-discrete formulation~\eqref{generalform_system} is entropy-stable.
\end{proof}

\begin{rem}
We further clarify the differences between the proposed damping term $\mathbf{OF}$ in method~\eqref{generalform_system} and classical artificial viscosity methods. Classical artificial viscosity methods introduce numerical dissipation through additional diffusion operators, often involving high-order even spatial derivatives \cite{Mattsson2004,Nordstrom2006,Stiernstrom2021}. In contrast, the proposed damping term is constructed from discrete local polynomial projection residuals of the entropy variables, with damping coefficients determined by jumps of numerical derivatives across subdomain interfaces. This construction allows the damping term to target local projection residuals while adaptively adjusting its strength according to the interface derivative jumps. For sufficiently smooth solutions, the global damping term satisfies an order $p+1$ estimate in the discrete $\widetilde H$-norm, as follows from Lemma~\ref{lemma_OF_acc}, without reducing the order of accuracy for the underlying SBP-SAT scheme~\eqref{generalform_nondp}.
\end{rem}
\begin{rem}
The damping term $\mathbf{OF}$ is incorporated as an additional term in the underlying semi-discrete SBP-SAT formulation~\eqref{generalform_nondp}. Although the damping coefficients are determined by numerical derivative jumps across neighboring subdomains, the local projection residuals are computed independently within each subdomain. The inclusion of the OF damping does not require any modifications to the SAT terms or upwind penalties. Moreover, Theorems~\ref{thm_cons} and~\ref{thm_L2} establish that the resulting oscillation-free SBP-SAT scheme~\eqref{generalform_system} retains the discrete conservation and entropy stability properties, demonstrating that the additional damping term does not compromise these properties.
\end{rem}

\begin{rem}
We further examine whether the OF damping imposes an additional restriction on the explicit time-step. The maximum stable CFL number is first determined for the scheme without OF damping. Using the same CFL number, the scheme with OF damping remains stable throughout the simulation, even for long-time simulations. These results indicate that no additional restriction on the explicit time-step size is observed in the numerical experiments considered.
 \end{rem}

\subsection{Two-dimensional system}\label{4_2}
For the two-dimensional case, we first introduce the corresponding two dimensional operators. We define
$D_x=D_1\otimes I_{N+1}$, $D_y=I_{N+1}\otimes D_1$, $\hat H=H\otimes H$.  These operators are further extended to the two-subdomain and $m$-variable setting as $\widetilde D_x=I_m\otimes I_2\otimes D_x$, $\widetilde D_y=I_m\otimes I_2\otimes D_y$, $\widetilde H=I_m\otimes I_2\otimes \hat H$ and $\mathcal{B}_x=(\mathcal B_I+\mathcal B_n) \otimes H$.  Since there are no internal interfaces in the $y$-direction, we define  $\mathcal{B}_y=H\otimes B_p$, where $B_p=-(\boldsymbol{\delta}_0-\boldsymbol{\delta}_N)(\boldsymbol{\delta}_0-\boldsymbol{\delta}_N)^{\top}$ satisfies $\mathbf 1^{\top}B_p=\mathbf 0^{\top}$. The upwind interface coefficients are $\boldsymbol{\alpha}_x
=\dfrac12\operatorname{diag}\left(\alpha_{x,1}^{(1)}+\alpha_{x,1}^{(2)},\alpha_{x,2}^{(1)}+\alpha_{x,2}^{(2)},\ldots,
\alpha_{x,m}^{(1)}+\alpha_{x,m}^{(2)}\right)$ and
$\boldsymbol{\alpha}_y=\operatorname{diag}
\left(\alpha_{y,1},\alpha_{y,2},\ldots,\alpha_{y,m}\right)$ with $\alpha_{x,\mu}^{(k)}\ge0$ and $\alpha_{y,\mu}=\operatorname{diag}
\left(\alpha_{y,\mu}^{(1)}, \alpha_{y,\mu}^{(2)}\right)\ge0$. For simplicity, we also use uniform grid spacings in both spatial directions, with $\Delta x=\Delta y$.
The corresponding oscillation-free SBP-SAT formulations can then be written in the general form
\begin{equation}\label{generalform_2D}
   \partial_t \mathbf{u}= -\mathbf{F}(\widetilde{D},\mathbf{u}) + \widetilde{H}^{-1}\mathbf{SAT}(\mathbf{u})+\widetilde{H}^{-1}\mathbf{T}(\mathbf{u})+\mathbf{OF},  
\end{equation}
where $\widetilde{D}=(\widetilde{D}_x,\widetilde{D}_y)$. $\widetilde H^{-1}\mathbf T(\mathbf u)=\dfrac12\widetilde H^{-1}\left(\boldsymbol{\alpha}_x\otimes\mathcal B_x+
\boldsymbol{\alpha}_y\otimes\mathcal B_y\right)\mathbf g$ represents the upwind term. $\mathbf{OF}_{\mu}=\left[\mathbf {OF}_{\mu}^{(1)},\mathbf {OF}_{\mu}^{(2)}
\right]^\top$ is the damping term for the $\mu$-th conserved variable.

As in the one-dimensional case, the damping term for the $\mu$-th conserved variable on each subdomain is given by
\[
\mathbf {OF}_{\mu}^{(k)}=-\sum_{l=0}^{s}\dfrac{\sigma_{\mu}^{(k,l)}}{|\Omega_k|}(\mathbf g_{\mu}^{(k)}-P^{\,l-1}\mathbf g_{\mu}^{(k)}),
\qquad \mu=1,\ldots,m,\quad k=1,2.
\]
where $P^{\,l}$ denotes the discrete $L^2$ projection in 2D. In particular, $P^{\,-1}=P^{\,0}$. $|\Omega_k|$ denotes the area of the $k$-th rectangular subdomain. In addition, $\sigma_{\mu}^{(k,l)}\ge 0$ is a scalar damping coefficient determined by the jumps at the vertices of the $k$-th subdomain.

We extend the damping term to the two-dimensional case. 
The grid points $\left\{\left(x_i^{(k)},y_j^{(k)}\right)\right\}_{i,j=0}^{N}$ 
on the subdomain $\Omega_k=[a_x^{(k)},b_x^{(k)}]\times[a_y^{(k)},b_y^{(k)}]$,  are mapped to the reference square
$[-1,1]^2$ by
\[
\xi_i=\dfrac{2(x_i^{(k)}-a_x^{(k)})}{b_x^{(k)}-a_x^{(k)}}-1,
\quad
\eta_j=\dfrac{2(y_j^{(k)}-a_y^{(k)})}{b_y^{(k)}-a_y^{(k)}}-1.
\]

The two-dimensional Legendre basis matrix
$\widetilde{\Phi}\in\mathbb R^{(N+1)^2\times(s+1)^2}$
is constructed by the tensor product of the one-dimensional Legendre basis matrices, that is,
$\widetilde{\Phi}=\Phi_x\otimes\Phi_y$,
where $\Phi_x$ and $\Phi_y$ denote the Legendre basis matrices in the $x$- and $y$-directions, respectively.
In particular, when the same reference points are used in both directions, $\Phi_x=\Phi_y=\Phi$.
For each $l=0,\ldots,s$, let
$\widetilde{\Phi}_l=\Phi_{x,l}\otimes\Phi_{y,l}
\in\mathbb R^{(N+1)^2\times(l+1)^2}$,
where $\Phi_{x,l}$ and $\Phi_{y,l}$ consist of the first
$l+1$ columns of $\Phi_x$ and $\Phi_y$, respectively.

For the $\mu$-th conserved variable on the $k$-th subdomain,
the discrete $L^2$ projection is defined by the following
orthogonality condition:
\begin{equation*}
\langle \mathbf g_{\mu}^{(k)}
-\widetilde{\Phi}_l\mathbf c_{\mu,l}^{(k)},
\widetilde{\Phi}_l\mathbf v\rangle_{\hat H}=0,
\quad \forall\,\mathbf v\in\mathbb R^{(l+1)^2}.
\end{equation*}
It follows that $\mathbf c_{\mu,l}^{(k)}=\left(
\widetilde{\Phi}_l^\top\hat H \widetilde{\Phi}_l\right)^{-1} \widetilde{\Phi}_l^\top
\hat H\mathbf g_{\mu}^{(k)}$. Here,
$\widetilde{\Phi}_l^\top\hat H\widetilde{\Phi}_l=
(\Phi_{x,l}^\top H\Phi_{x,l})\otimes (\Phi_{y,l}^\top H\Phi_{y,l})$.
Therefore, the discrete $L^2$ projection is given by
\begin{equation*}
P^{\,l}\mathbf g_{\mu}^{(k)}=
\widetilde{\Phi}_l\mathbf c_{\mu,l}^{(k)}=
\sum_{i=0}^{l}\sum_{j=0}^{l}
c_{\mu,i,j}^{(k)}
\left( L_i(\boldsymbol{\xi})
\otimes
L_j(\bm{\eta})
\right),
\end{equation*}
where $\boldsymbol{\xi}=\left[\xi_0,\xi_1,\ldots,\xi_N\right]^{\top} $ and $\bm{\eta}=\left[\eta_0,\eta_1,\ldots,\eta_N\right]^{\top}$ are the reference vectors. In particular, we define $P^{-1}=P^0$.

Let $\mathcal V_k$ denote the set of the four vertices of $\Omega_k$. For $0\le p\le l$, the discrete mixed derivative of $\mathbf z_\mu^{(k)}$ at a vertex $v=(x_{\mathcal I_x}^{(k)},y_{\mathcal I_y}^{(k)})\in\mathcal V_k$, where $\mathcal I_x,\mathcal I_y \in \{0,N\}$, is evaluated by
\[
\left.
\mathcal D_x^p\mathcal D_y^{l-p}
\mathbf z_{\mu}^{(k)}
\right|_v
=\left(
\mathbf d_{\mathcal I_x}^{(k, p)}
\otimes
\mathbf d_{\mathcal I_y}^{(k, l-p)}
\right)^\top
\mathbf z_{\mu}^{(k)},
\]
where $\mathbf d_{\mathcal I_x}^{(k, p)}$ and $\mathbf d_{\mathcal I_y}^{(k, l-p)}$ are the same one-sided
finite difference coefficient vectors used in the
one-dimensional construction. In particular, $\mathbf d_{\mathcal I_x}^{(k, 0)}$ and $\mathbf d_{\mathcal I_y}^{(k, 0)}$ are defined as the endpoint selection vectors in each spatial direction.  The auxiliary variable $\mathbf z_\mu^{(k)}$ remains equal to $\mathbf u_\mu^{(k)}$ for the advection and Burgers' equations, while its definition for the shallow water equations is given in~\eqref{z_2D}.

The two-dimensional damping coefficient is defined by
\[
\sigma_{\mu}^{(k,l)}=C_l (\Delta x)^l\sum_{p=0}^{l}\left[
\dfrac{1}{4}\sum_{v\in\mathcal V_k}\left(\left.[\mkern-2.8mu[\mathcal D_x^p
\mathcal D_y^{\,l-p}\mathbf z_{\mu}^{(k)}
]\mkern-2.8mu]\right|_v\right)^2\right]^{1/2},
\]
To illustrate the computation of the vertex jumps, let \(v\) be the upper-right vertex of \(\Omega_1\). The squared jump of the mixed derivative at \(v\) consists of two contributions: one between \(\Omega_1\) and its neighboring subdomain \(\Omega_2\) in the \(x\)-direction, and the other between the upper-right and lower-right vertices of \(\Omega_1\) in the periodic \(y\)-direction, as follows:
\[
\begin{aligned}
\left(\left.[\mkern-2.8mu[\mathcal D_x^p\mathcal D_y^{l-p}
\mathbf z_\mu^{(1)}]\mkern-2.8mu]\right|_v\right)^2
=& \left(\left.\mathcal D_x^p\mathcal D_y^{l-p}
\mathbf z_\mu^{(2)}\right|_v-\left.\mathcal D_x^p\mathcal D_y^{l-p}\mathbf z_\mu^{(1)}\right|_v\right)^2 \\
+& \left(\left.\mathcal D_x^p\mathcal D_y^{l-p}
\mathbf z_\mu^{(1)}\right|_v-\left.\mathcal D_x^p\mathcal D_y^{l-p}\mathbf z_\mu^{(1)}\right|_{v_{\mathrm{per}}}\right)^2.
\end{aligned}
\]
Here, \(v_{\mathrm{per}}\) denotes the lower-right vertex of \(\Omega_1\), periodically connected to \(v\) in the \(y\)-direction. The jumps at the remaining vertices are computed analogously, with periodic connections included when determining the neighboring vertices.
\begin{rem}
  The proposed damping term $\mathbf{OF}$ can be extended to $d$ spatial dimensions following the same approach as in the two-dimensional case. The discrete operators are extended through tensor products. For the damping coefficients, the mixed derivative is defined by $\mathcal D^{\boldsymbol{\rho}}
=\mathcal D_{x_1}^{\rho_1}\cdots\mathcal D_{x_d}^{\rho_d}$,
where the multi-index $\boldsymbol{\rho}=(\rho_1,\ldots,\rho_d)$ satisfies $|\boldsymbol{\rho}|=\rho_1+\cdots+\rho_d=l$. Otherwise,
when computing the derivative jumps at the vertices, only neighboring subdomains sharing a $(d-1)$-dimensional face with the current subdomain are considered.
\end{rem}

Indeed, the two-dimensional SBP operators, norm matrices and discrete
$L^2$ projection operators are constructed by tensor products of their one-dimensional counterparts. 
The SBP property and orthogonality of the discrete $L^2$ projection are preserved in 2D. 
Therefore, the conservation and entropy stability for the two-dimensional system follow directly from the corresponding one-dimensional results through the tensor-product construction.

\begin{thm}\label{thm_2D}
Consider the two-dimensional semi-discrete oscillation-free SBP-SAT
formulation~\eqref{generalform_2D} with periodic boundary conditions.
Then the following properties hold:
\begin{enumerate}
\item  The semi-discrete approximation~\eqref{generalform_2D} is conservative for each conserved variable, that is,
\[
\dfrac{\mathrm d}{\mathrm dt}
\sum_{k=1}^{2} \langle \mathbf 1,\mathbf u_{\mu}^{(k)}\rangle_{\hat H}=0,
\qquad \mu=1,\ldots,m.
\]

\item  The semi-discrete approximation~\eqref{generalform_2D} is entropy-stable, that is,
\[
\dfrac{\mathrm dE_h}{\mathrm dt}\leq 0,
\]
where $E_h$ denotes the discrete entropy of the system.
\end{enumerate}
\end{thm}

\begin{proof}
For the conservation property, the underlying two-dimensional
SBP-SAT discretization without the damping term is conservative.
Therefore, it only remains to consider the contribution of the
damping term. By the definition
$\mathbf r_{\mu}^{(k,l)}
=\mathbf g_{\mu}^{(k)}-P^{\,l-1}\mathbf g_{\mu}^{(k)}$
and the orthogonality of the discrete projection, we have
$\langle\mathbf 1,\mathbf r_{\mu}^{(k,l)}
\rangle_{\hat H}=0$ for $l=0,\ldots,s$. Therefore,
\[
\begin{aligned}
\langle
\mathbf 1,\mathbf {OF}_{\mu}^{(k)}\rangle_{\hat H}
&=-\sum_{l=0}^{s}
\dfrac{\sigma_{\mu}^{(k,l)}}{|\Omega_k|}
\langle \mathbf 1,\mathbf r_{\mu}^{(k,l)}
\rangle_{\hat H}=0.
\end{aligned}
\]
After summing over the two subdomains,
\[
\sum_{k=1}^{2}
\langle
\mathbf 1,\mathbf {OF}_{\mu}^{(k)}
\rangle_{\hat H}=0.
\]
Together with the conservation property of the underlying
SBP-SAT discretization, this yields
\[
\dfrac{\mathrm d}{\mathrm dt}
\sum_{k=1}^{2}
\langle\mathbf 1,\mathbf u_{\mu}^{(k)}
\rangle_{\hat H}=0,
\quad \mu=1,\ldots,m.
\]
Thus, part~(i) follows.

For the entropy-stable property, the underlying
two-dimensional SBP-SAT discretization without the damping
term is entropy-stable. It therefore remains to
consider the contribution of the damping term.
Using
$\mathbf g_{\mu}^{(k)}
=P^{\,l-1}\mathbf g_{\mu}^{(k)}
+\mathbf r_{\mu}^{(k,l)}$
and the orthogonality of the discrete projection gives
\[
\begin{aligned}
\langle\mathbf g_{\mu}^{(k)},\mathbf r_{\mu}^{(k,l)}
\rangle_{\hat H}
=\langle P^{\,l-1}\mathbf g_{\mu}^{(k)},
\mathbf r_{\mu}^{(k,l)}
\rangle_{\hat H}
+\left\|\mathbf r_{\mu}^{(k,l)}\right\|_{\hat H}^2=\left\|\mathbf r_{\mu}^{(k,l)}\right\|_{\hat H}^2.
\end{aligned}
\]
Therefore,
\[
\begin{aligned}
\sum_{k=1}^{2}\sum_{\mu=1}^{m}
\langle
\mathbf g_{\mu}^{(k)},\mathbf {OF}_{\mu}^{(k)}
\rangle_{\hat H}
=-\sum_{k=1}^{2}\sum_{\mu=1}^{m}\sum_{l=0}^{s}
\dfrac{\sigma_{\mu}^{(k,l)}}{|\Omega_k|}
\left\|\mathbf r_{\mu}^{(k,l)}
\right\|_{\hat H}^2
\leq0.
\end{aligned}
\]
Together with the entropy stability of the
underlying SBP-SAT discretization, we obtain
\[
\dfrac{\mathrm dE_h}{\mathrm dt}\leq0.
\]
Thus, part~(ii) follows.
\end{proof}

\begin{rem}
Although the difference operators used for multi-dimensional problems and systems of conservation laws may take different forms, they can all be obtained from the one-dimensional SBP operator $D_1$ through tensor-product extensions. Throughout this paper, $\widetilde D$ is used as a generic notation for the corresponding extended difference operator in the setting under consideration. Similarly, $\widetilde H$ denotes the corresponding extended norm matrix. Their precise forms depend on the dimension and on whether the problem is scalar or a system. 
\end{rem}

\section{Numerical experiments}\label{Sec.5}
In this section, we present numerical experiments to validate the theoretical results obtained in the previous sections and to assess the performance of the proposed numerical method. We adopt the third-order strong stability-preserving Runge-Kutta (SSP-RK3) method for time integration \cite{Shu1988}, employ a uniform grid with spatial step denoted by $\Delta x$ and use periodic boundary conditions in all cases. The numerical tests are presented in order of increasing complexity. 
We first consider one-dimensional tests for the linear advection equation, Burgers' equation, and the shallow water equations, followed by the corresponding two-dimensional tests. Throughout the numerical experiments, we compare the errors with and without the damping term, and examine the ability of the damping term to control spurious oscillations near discontinuities. In addition, for the shallow water equations, we carry out several additional benchmark tests and physically motivated examples to further examine the robustness and stability of the proposed method.

\subsection{One-dimensional problems }
\begin{figure}[!htb]     
   \centering 
  \subfigure[ Without $\mathbf{OF}$]
  {
     \includegraphics[width=0.45\textwidth]{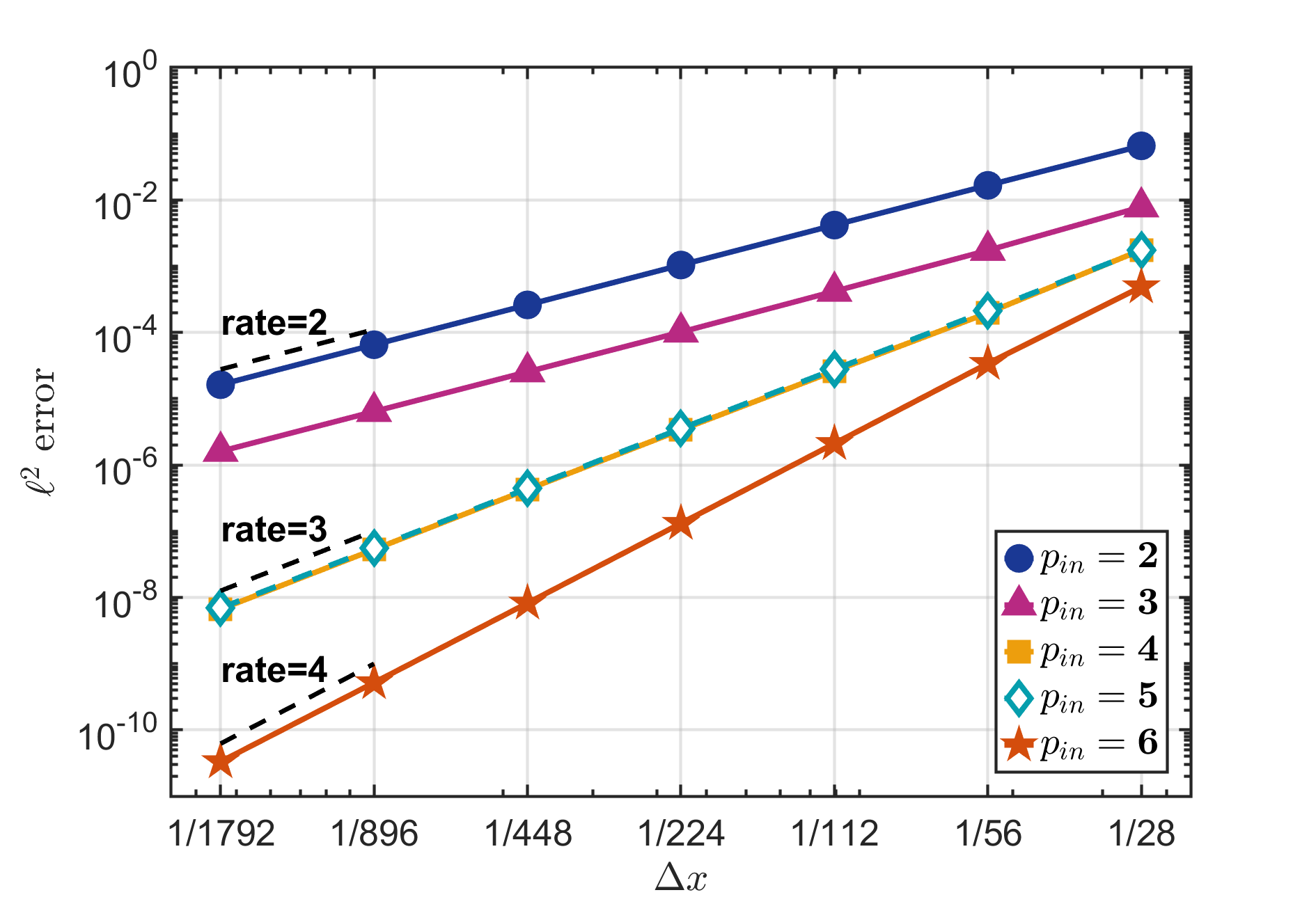}
  }\hspace{-2.2em}
  \subfigure[  With $\mathbf{OF}$]
  {
      \includegraphics[width=0.45\textwidth]{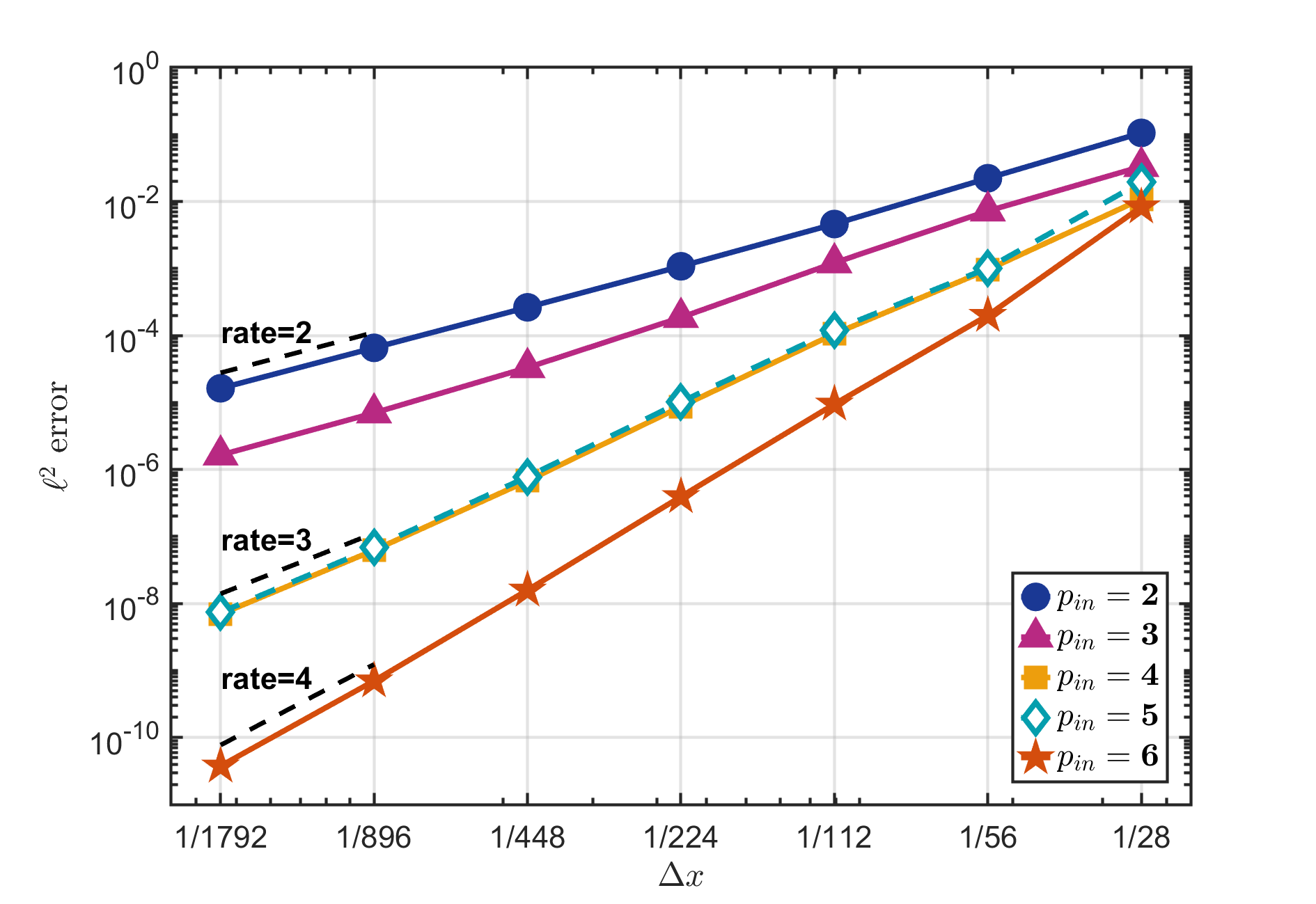}
  }
 \captionsetup{name={Fig.},labelsep=period}
  \caption{ $\ell^2$ errors and convergence rates initial value (a) in  Example~\ref{ex1}.}   
  \label{adv_conv_rate}
\end{figure}

\begin{figure}[!htb]     
   \centering 
  \subfigure[$K=2^7$]
  {
     \includegraphics[width=0.33\textwidth]{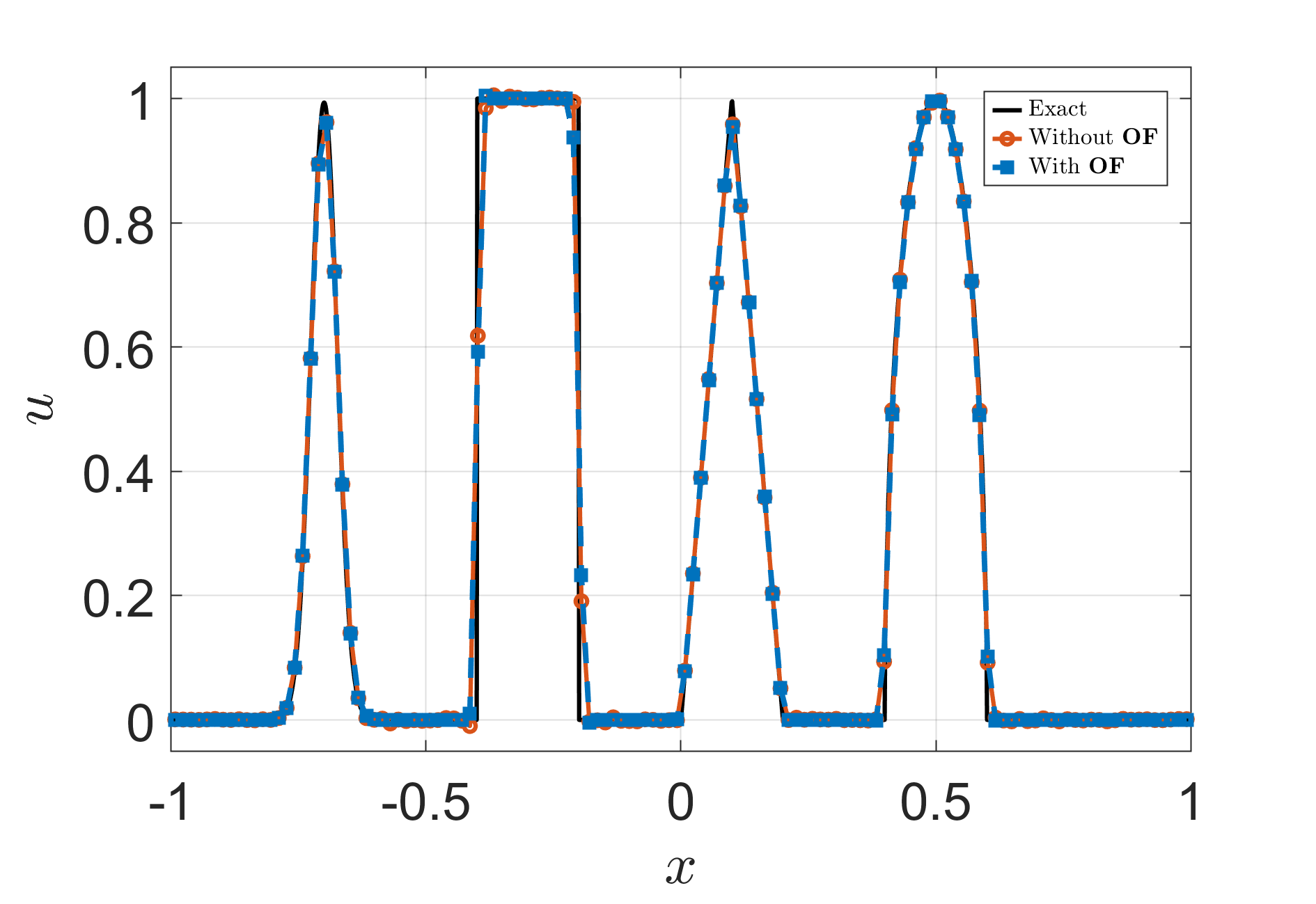} \label{adv_Dis_ic_a}
  }\hspace{-2.2em}
  \subfigure[  $|\mathbf {OF}|$]
  {
      \includegraphics[width=0.33\textwidth]{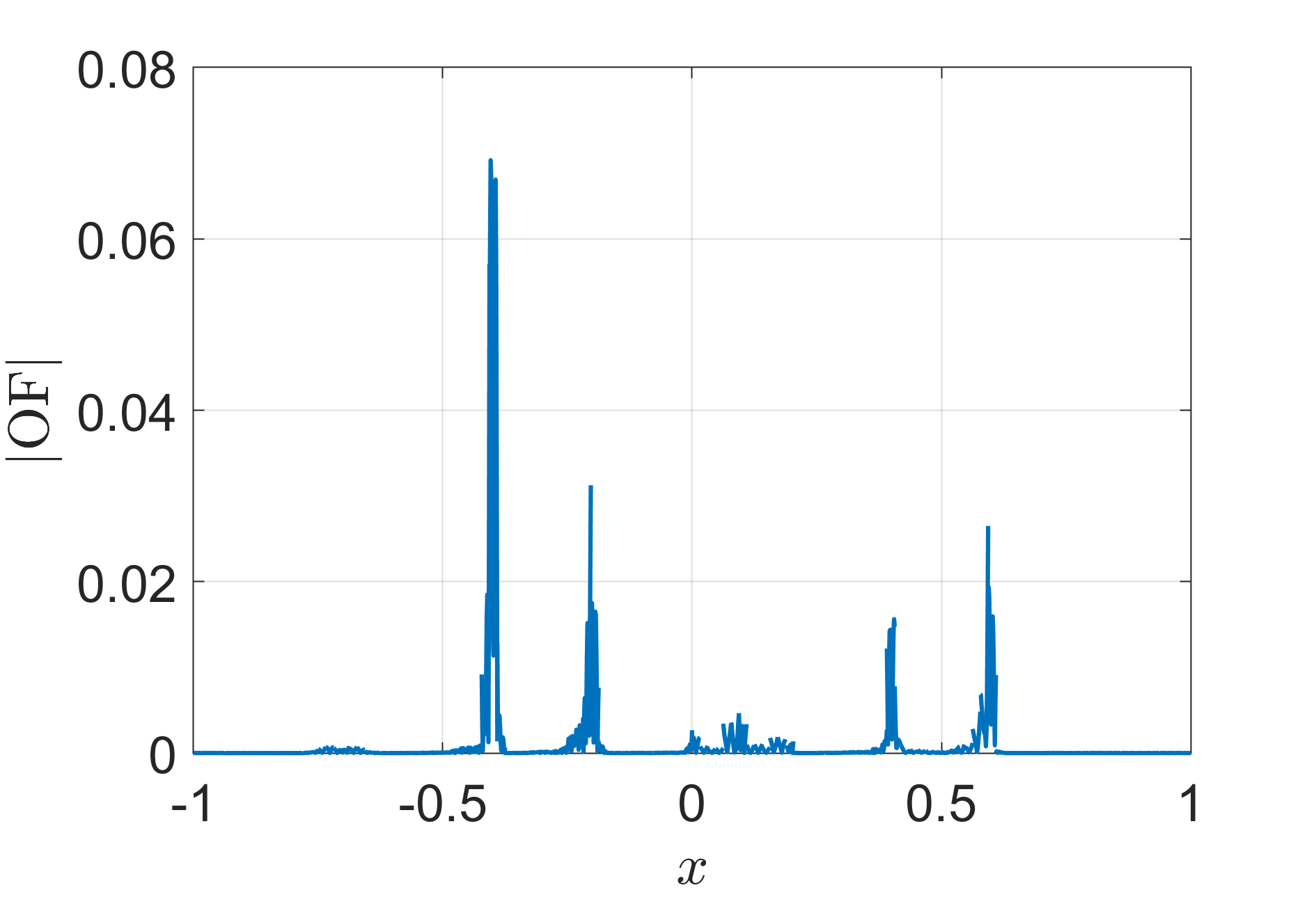} \label{adv_Dis_ic_b}
  }\hspace{-2.2em}
  \subfigure[  $\sigma^{(k,l)}$]
  {
      \includegraphics[width=0.33\textwidth]{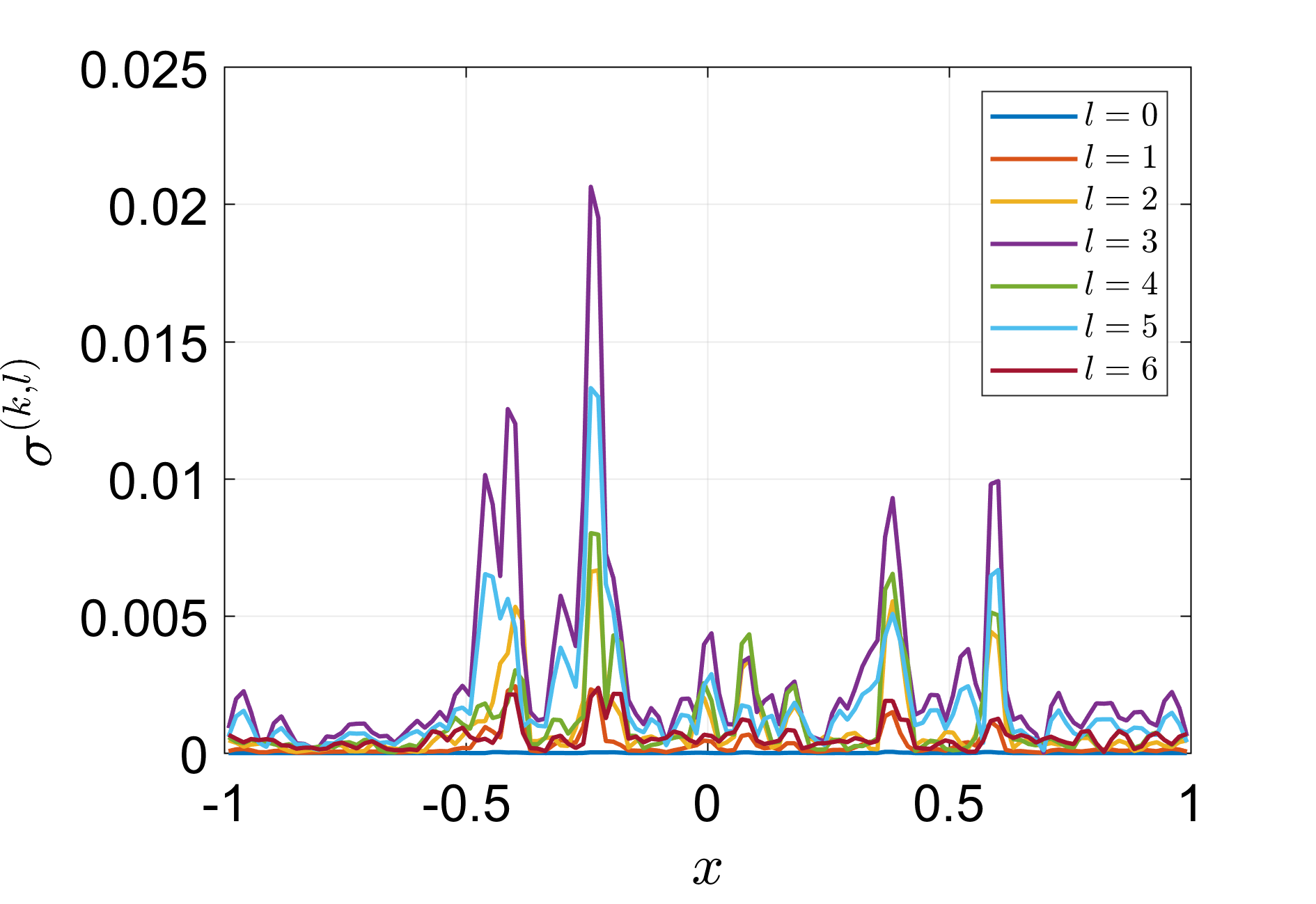} \label{adv_Dis_ic_c}
  }
 \captionsetup{name={Fig.},labelsep=period}
  \caption{The numerical solution and the corresponding  $|\mathbf {OF}|$, $\sigma^{(k,l)}$ for the $6$th order method with initial value (b) in  Example~\ref{ex1}.}   \label{adv_Dis_ic}
\end{figure}

\begin{figure}[!htb]     
   \centering 
  \subfigure[Without $\mathbf{OF}$]
  {
     \includegraphics[width=0.45\textwidth]{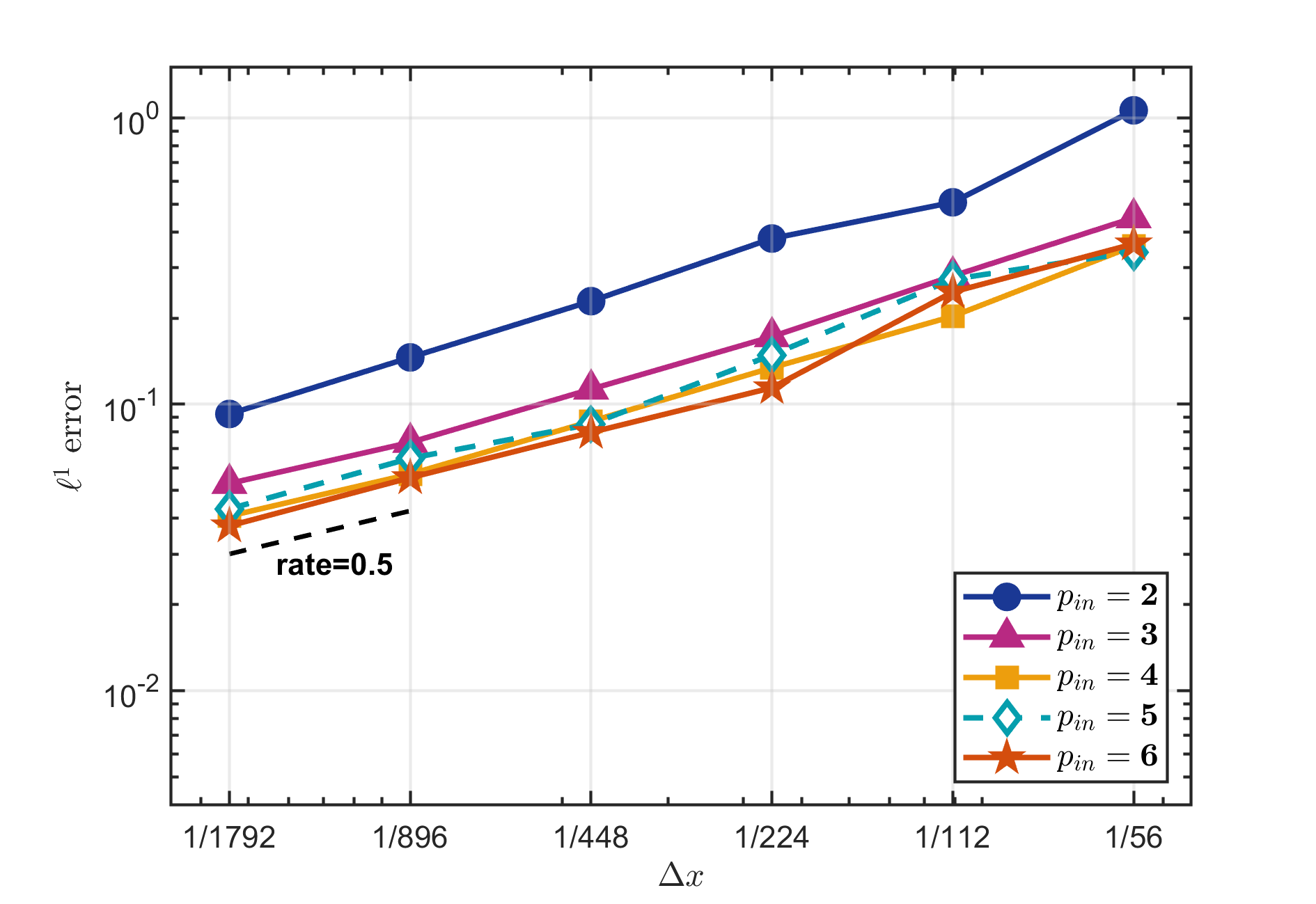}
  }\hspace{-2.2em}
  \subfigure[ With $\mathbf{OF}$]
  {
      \includegraphics[width=0.45\textwidth]{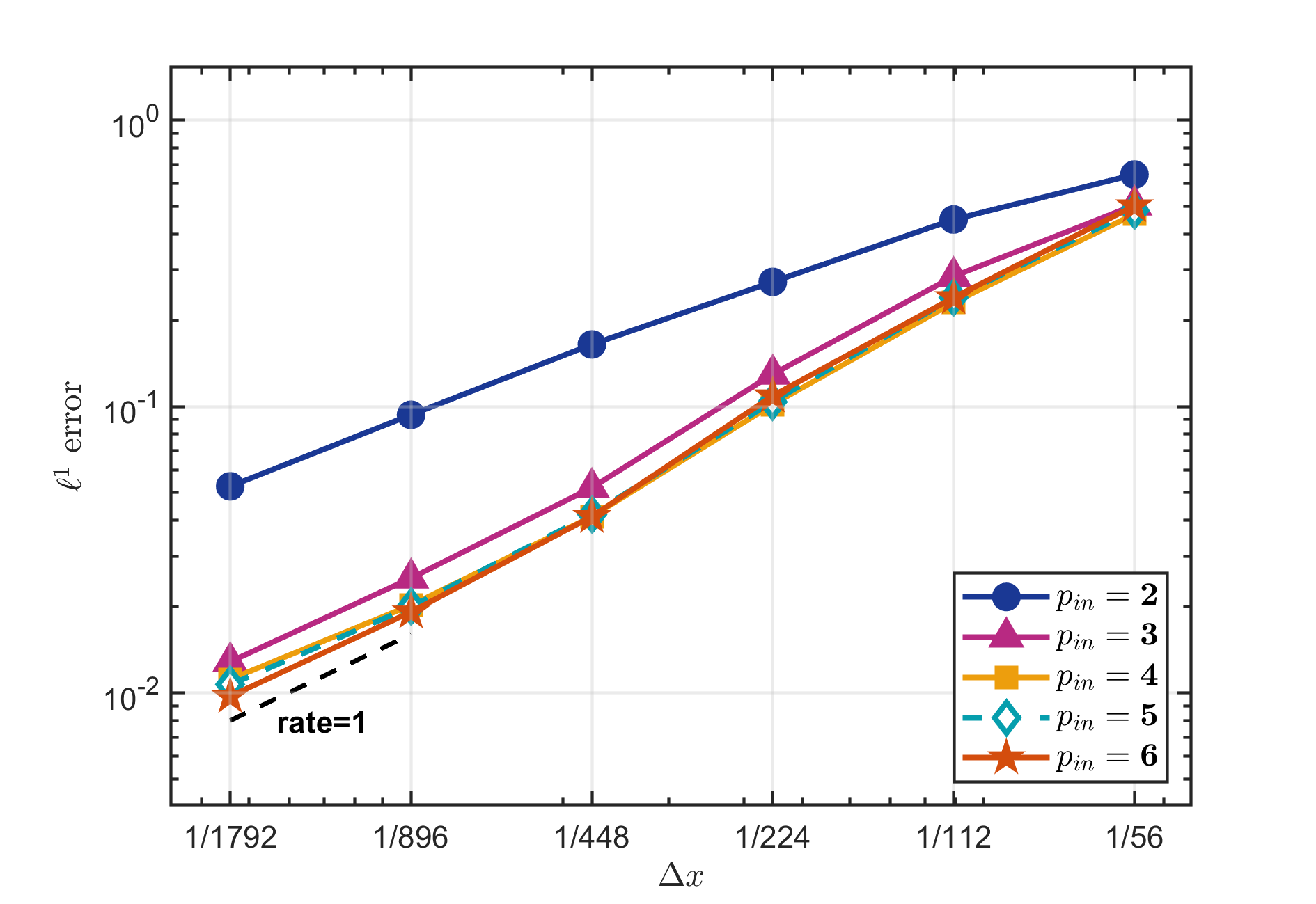}
  }
 \captionsetup{name={Fig.},labelsep=period}
  \caption{$\ell^1$ errors and convergence rates initial value (b) in  Example~\ref{ex1}.}   \label{adv_Dis_conver}
\end{figure}

In this subsection, we present numerical examples for the one-dimensional problems. In all cases, each subdomain is discretized by $N+1$ grid points with $N=14$. 
The time step is chosen as $\Delta t=\mathrm{CFL}\,\Delta x/\lambda$, where $\mathrm{CFL}$ is the Courant number and $\lambda$ is the maximum characteristic wave speed. For linear advection, $\lambda=|a|$. For Burgers' and shallow water equations, $\lambda$ is evaluated from the initial conditions as $\max |u|$ and $\max (|u|+\sqrt{gh})$ with $g$ denoting the gravitational acceleration.
To assess the spatial accuracy of the proposed scheme, we conduct convergence tests using SBP operator $D_1$. For each order, the mesh is refined by increasing the number of subdomains as $K=2^{\gamma}$, $\gamma=1,\ldots,7$, while keeping the number of grid points per subdomain fixed. In particular, the convergence rates are determined from the discrete $\ell^2$ errors. For the global error vector $\mathbf e=\mathbf u-\mathbf U$, the discrete error norm is defined by $\|\mathbf e\|^2_{\widetilde H}=\mathbf e^T\widetilde H\mathbf e$.
\begin{example}[Advection equation in 1D]\label{ex1}
 For the one-dimensional linear advection equation~\eqref{adv}, the oscillation-free SBP-SAT scheme is obtained from the generic formulation~\eqref{generalform_system} by setting $m=1$ and using the corresponding flux function. In the tests, we take the coefficient of the upwind interface term $\alpha=\dfrac{1}{2}\left(\alpha^{(k)}+\alpha^{(k+1)}\right), 
 \alpha^{(k)}=1$, for $k\ge 1$.  The following two initial conditions are considered.
 
 (a) The smooth initial condition:
\begin{equation*}
U_0(x) = \sin(2\pi x+0.3).
\end{equation*}
(b) The non-smooth initial condition :
\begin{equation*}
\scalebox{0.9}{$\displaystyle
U_0(x)=
\begin{cases}
\dfrac{1}{6}\left(
G(x,z-\delta)+G(x,z+\delta)+4G(x,z)
\right),
& -0.8\le x\le -0.6,\\
1,& -0.4\le x\le -0.2,\\
1-\left|10(x-0.1)\right|,
& 0\le x\le 0.2,\\
\dfrac{1}{6}\left(
F(x,a-\delta)+F(x,a+\delta)+4F(x,a)\right), & 0.4\le x\le 0.6,\\
0,
& \text{otherwise},
\end{cases}
$}
\end{equation*}
where $G(x)=e^{-\dfrac{\log 2}{36\delta^2}(x-z)^2}$ and $F(x)=\sqrt{\max\left(1-100(x-a)^2,0\right)}$. The parameters are chosen as $\delta=10$, $z=-0.7$ and $a=0.5$.

For the smooth case (a), we take $\Omega=[0,1]$ and $T=1$. For the non-smooth case (b), we take $\Omega=[-1,1]$ and $T=2$. The time step is set to $\Delta t=0.1\Delta x$ in both cases.

 Fig.~\ref{adv_conv_rate} shows the convergence results for case (a), with and without the $\mathbf{OF}$ damping terms. In both cases, the observed convergence rates are consistent with the expected convergence rates listed in Table~\ref{accuracy}. The results also show that the global convergence rate depends on the boundary accuracy: the SBP operator $D_1$ with the same interior order $p_1$ but different boundary orders $p_b$ exhibit different convergence rates, while those with the same boundary order exhibit the same rate. A comparison of the errors also shows that the results with and without the damping term become increasingly close as the mesh is refined, indicating that the damping term has little effect on the accuracy for smooth exact solutions.

Figures~\ref{adv_Dis_ic} and~\ref{adv_Dis_conver} correspond to case (b), non-smooth initial conditions. With all other parameters fixed, we use the 6th-order SBP operator and compare the exact solution with two numerical solutions (with $\mathbf{OF}$ and without $\mathbf{OF}$) for $K=2^7$.  As shown in Fig.~\ref{adv_Dis_ic_a}, the oscillation-free SBP-SAT method accurately captures the sharp features and discontinuities, whereas spurious oscillations are observed without the damping term. We have obtained similar results using 2nd-, 3rd-, 4th-, and 5th-order SBP operators, with the exception that the 2nd-order operator requires doubling the damping coefficient. Figures~\ref{adv_Dis_ic_b} and ~\ref{adv_Dis_ic_c} show the corresponding absolute value of the damping term $\mathbf{OF}$ and the damping coefficients $\sigma^{(k,l)}$, respectively. Both exhibit similar patterns, with their magnitudes increasing near the discontinuities in Fig.~\ref{adv_Dis_ic_a} and approaching zero in smooth regions. 
For discontinuous solutions, classical first-order finite difference methods may achieve $0.5$-order accuracy in the $\ell^1$ norm~\cite{LeVeque1992}. 
Figure~\ref{adv_Dis_conver} shows the convergence results in the discrete $\ell^1$ norm for the discontinuous initial condition. All methods exhibit convergence rates higher than $0.5$, both with and without OF damping. Moreover, the proposed scheme with the damping term preserves the convergence behavior and generally reduces the numerical errors on finer meshes.

Finally, we also repeated the above two tests using the upwind SBP operators developed in~\cite{Mattsson2017}. The results point to the results and data show no significant differences from those obtained using the SBP operator $D_1$.
\end{example}

\begin{figure}[!htb]     
   \centering 
  \subfigure[ Without $\mathbf{OF}$]
  {
     \includegraphics[width=0.45\textwidth]{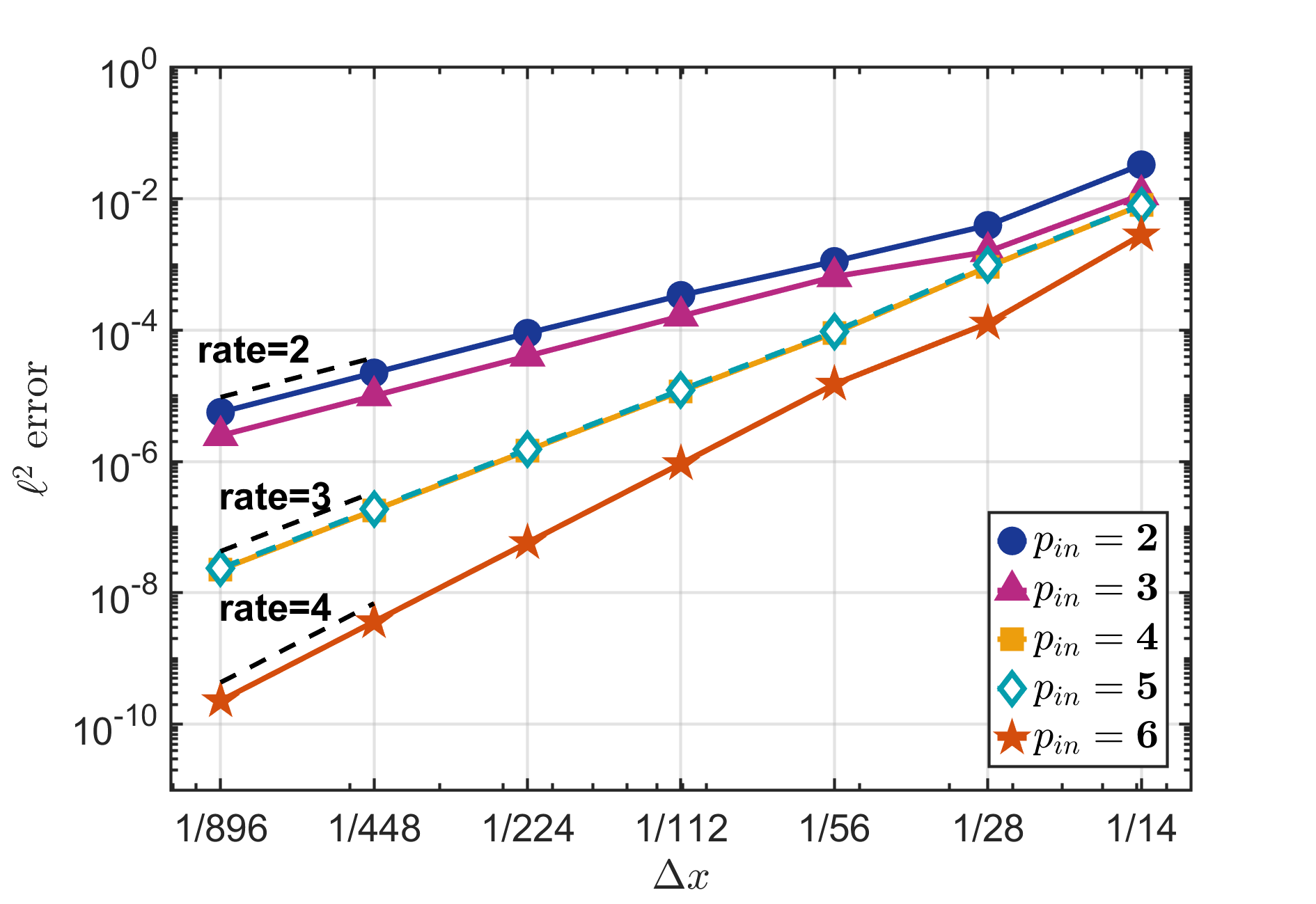}
  }\hspace{-2.2em}
  \subfigure[  With $\mathbf{OF}$]
  {
      \includegraphics[width=0.45\textwidth]{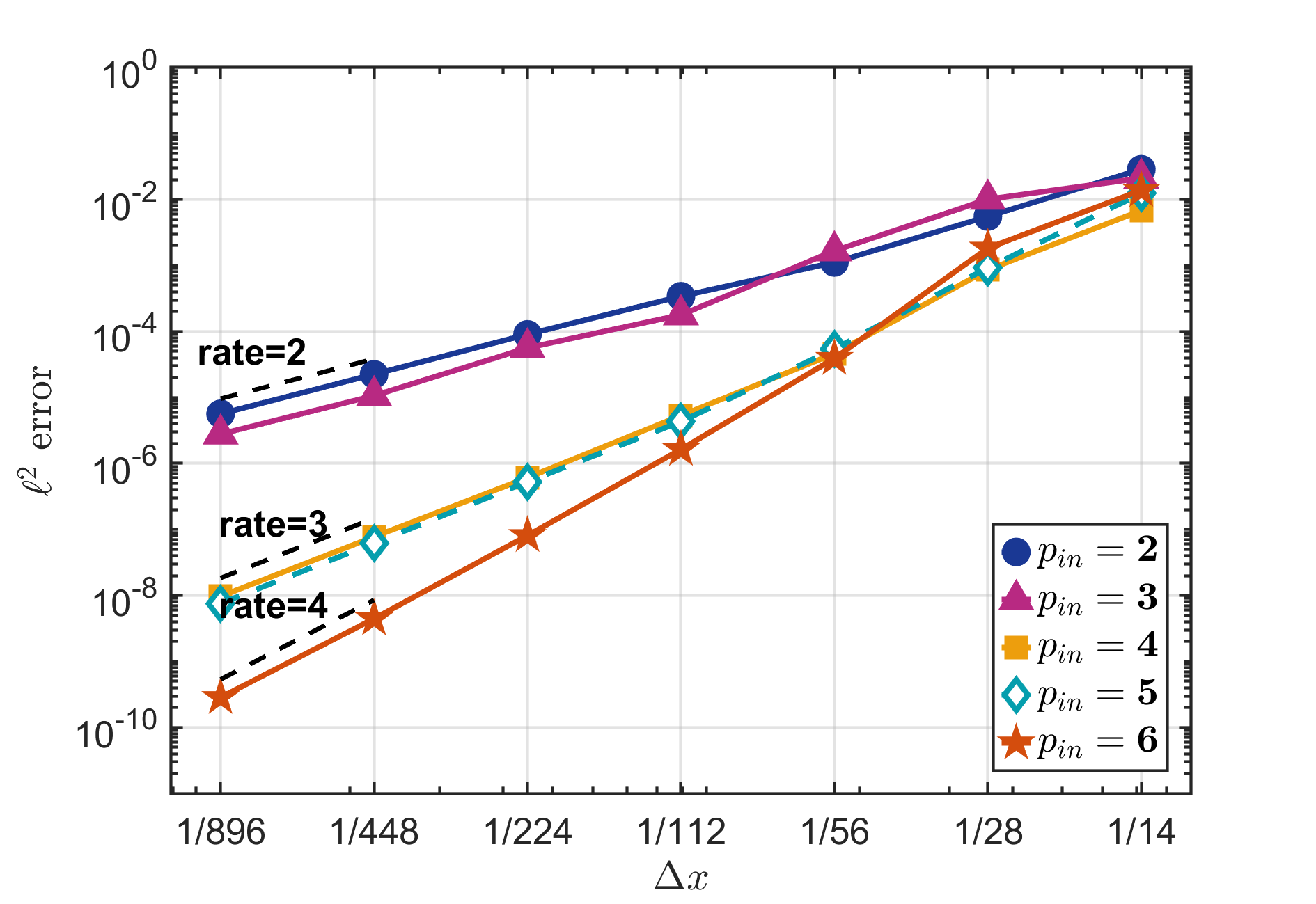}
  }
 \captionsetup{name={Fig.},labelsep=period}
  \caption{ $\ell^2$ errors and convergence rates with manufactured solutions $U(x,t)$ in  Example~\ref{ex2}.}   
  \label{bur_conv_rate}
\end{figure}

\begin{figure}[!htb]     
   \centering 
  \subfigure
  {
     \includegraphics[width=0.33\textwidth]{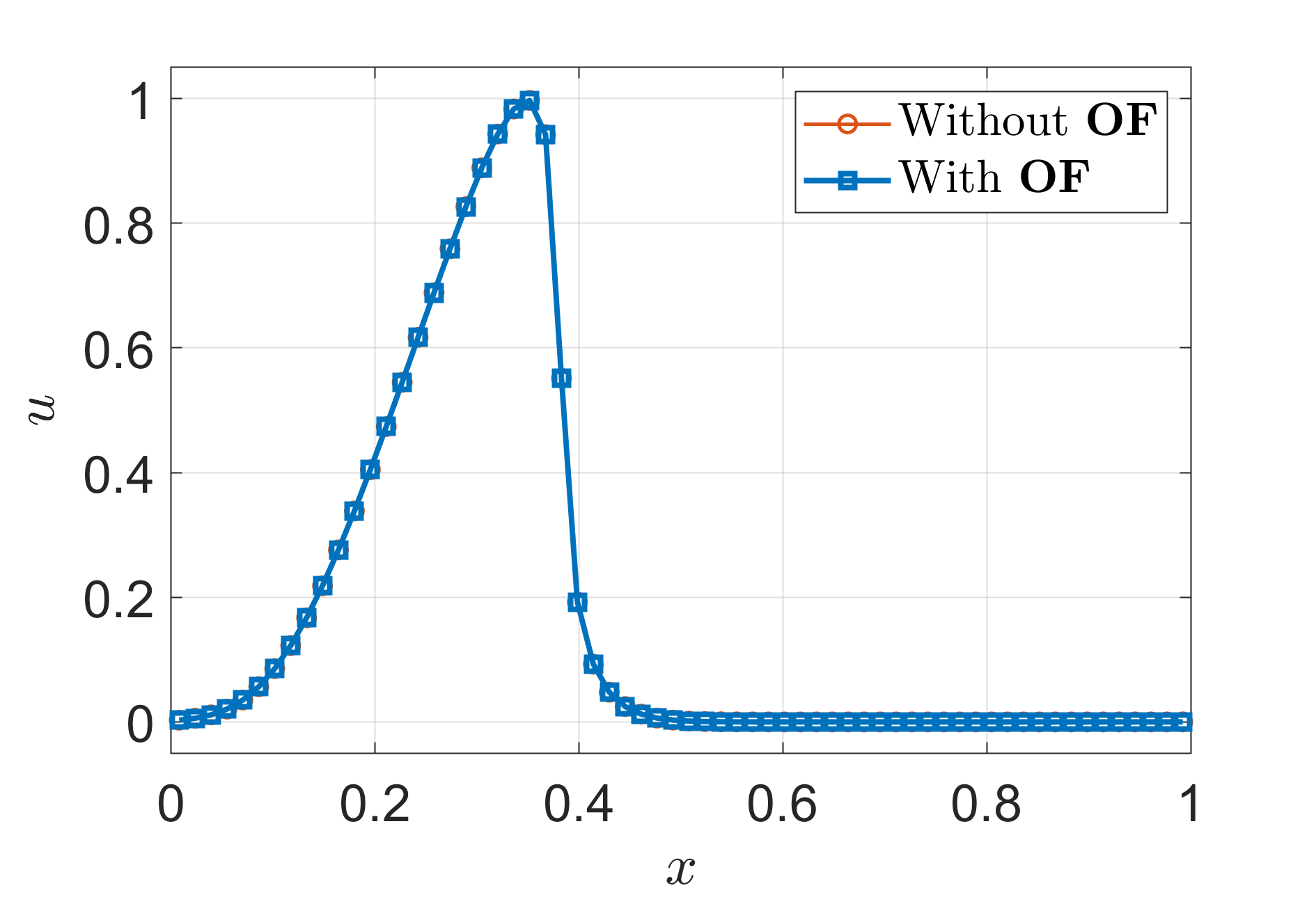}
  }\hspace{-2em}
  \subfigure 
  {
      \includegraphics[width=0.33\textwidth]{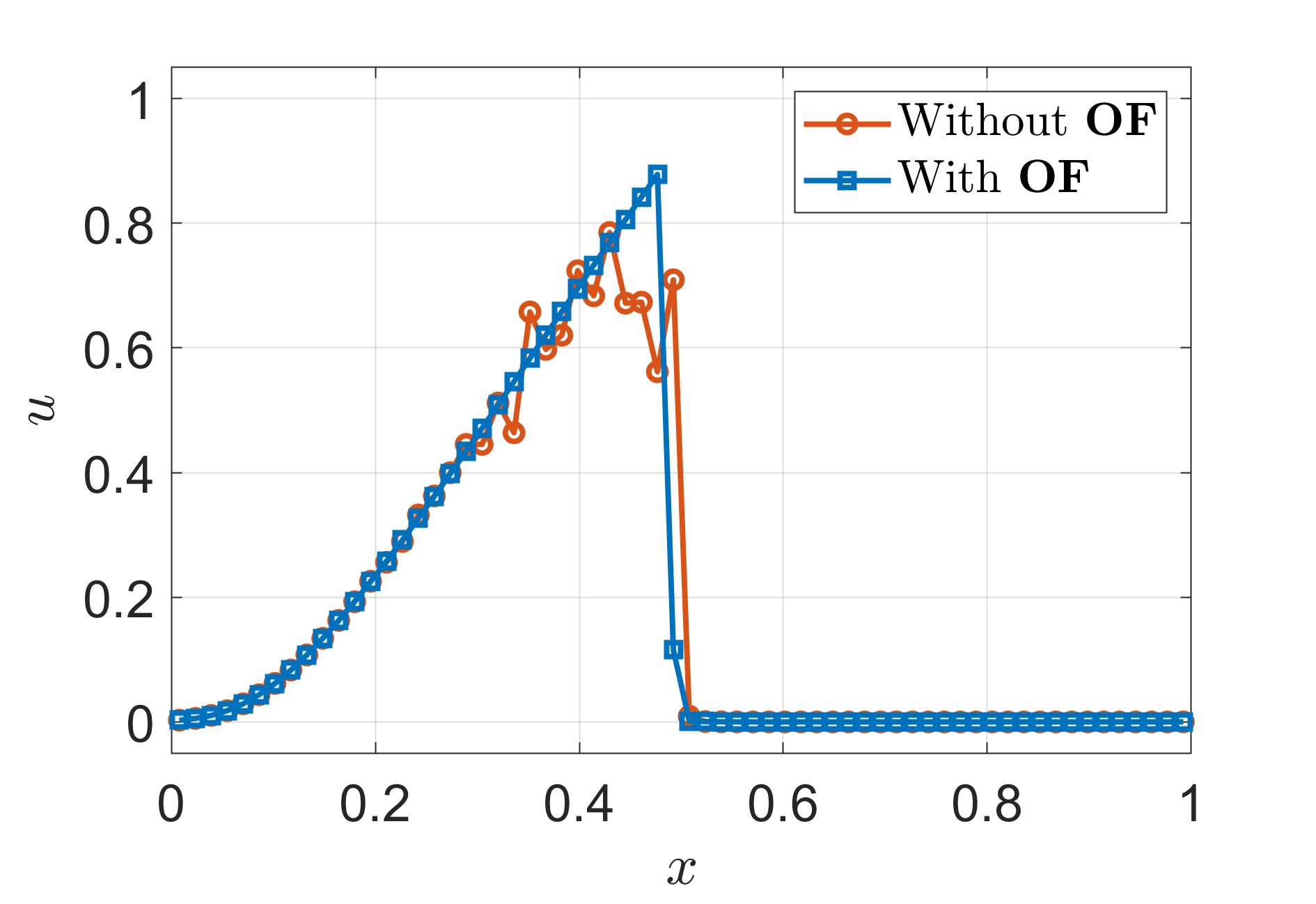}
  }\hspace{-2em}
    \subfigure 
  {
      \includegraphics[width=0.33\textwidth]{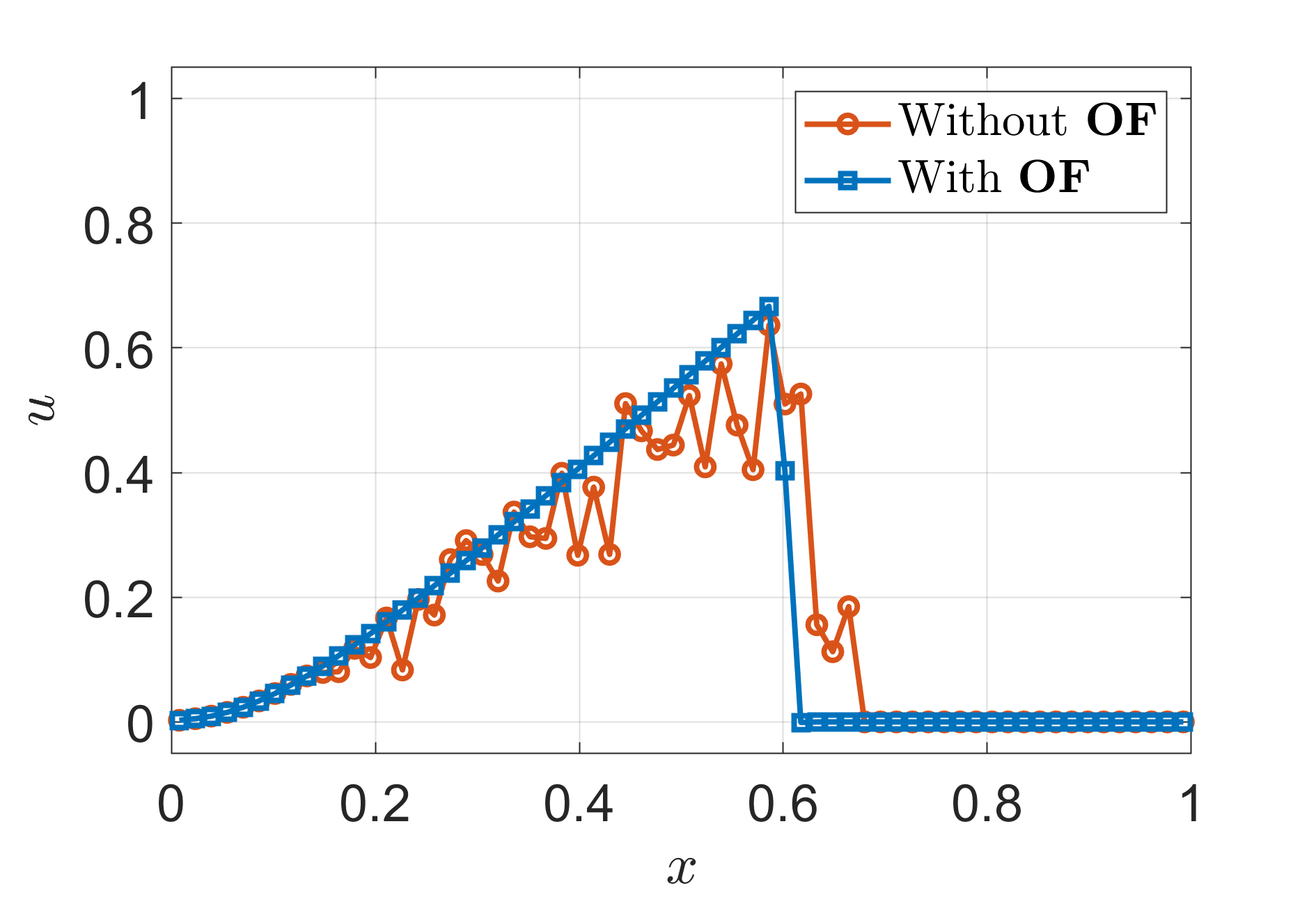}
  } \\[-1em]
   \centering
    \subfigure[$t=0.1$]  
  {
      \includegraphics[width=0.33\textwidth]{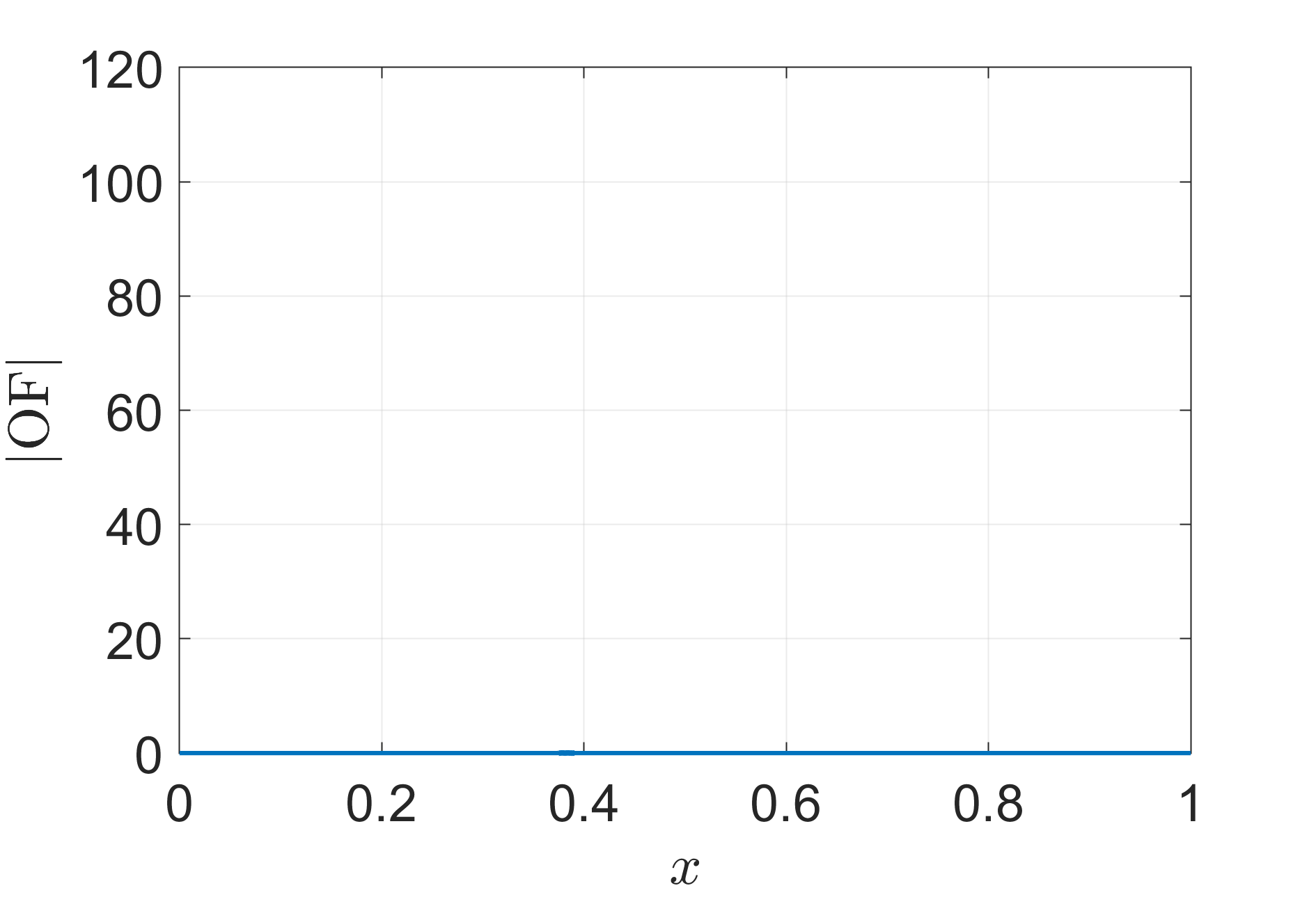}
  }\hspace{-2em}
    \subfigure[$t=0.3$]  
  {
      \includegraphics[width=0.33\textwidth]{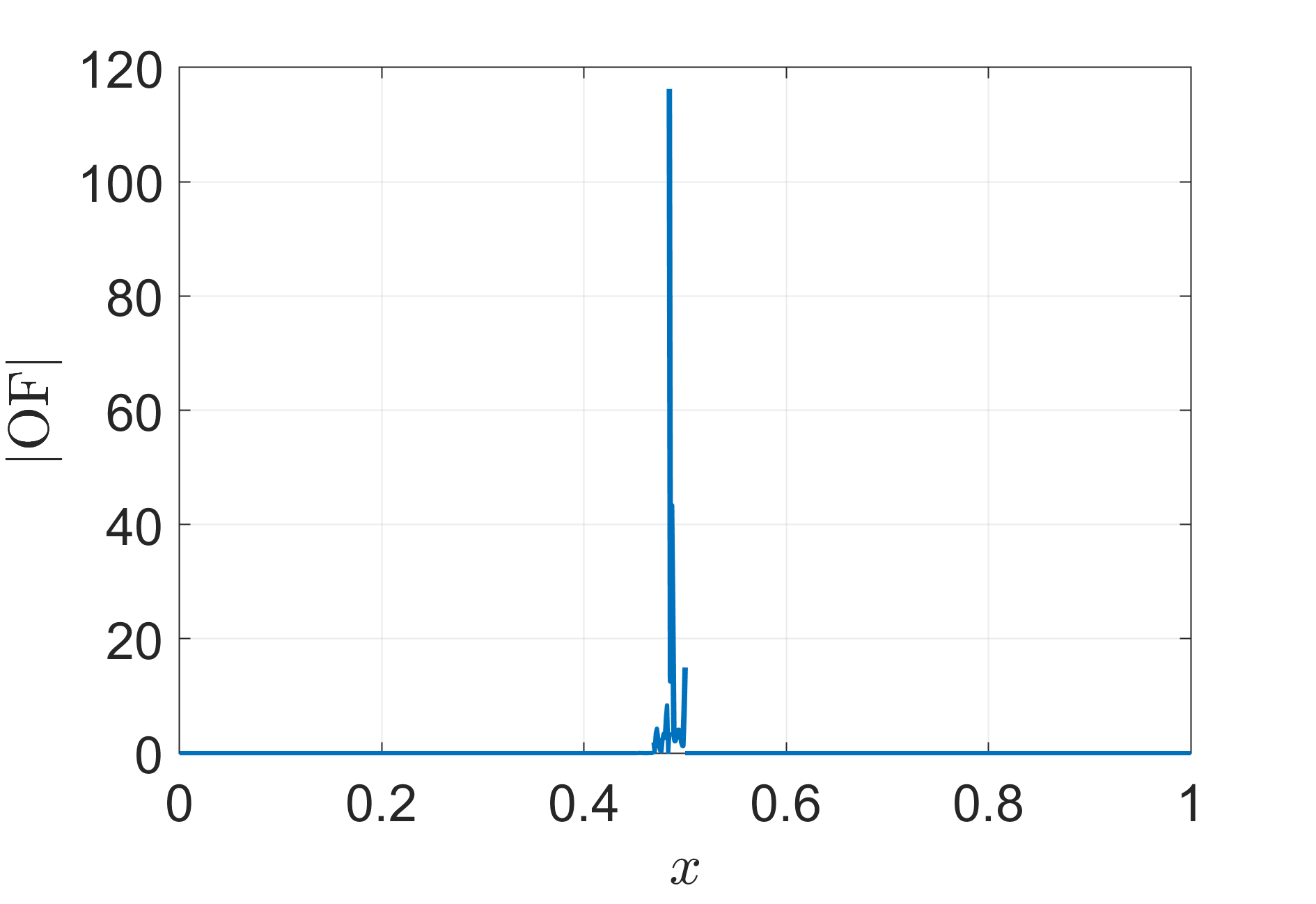}
   }\hspace{-2em}
    \subfigure[$t=0.6$]  
  {
      \includegraphics[width=0.33\textwidth]{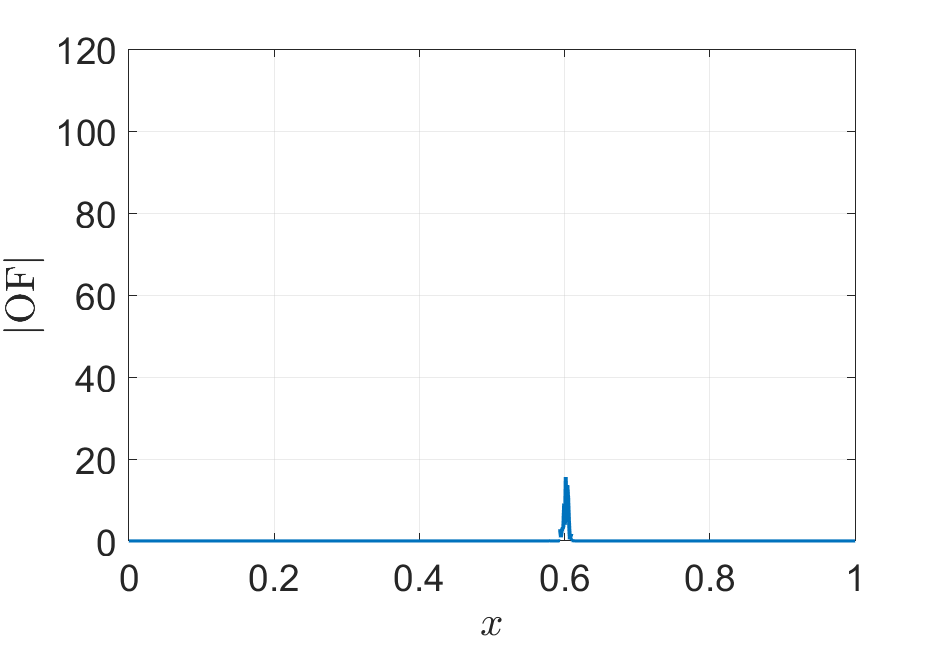}
  } 
 \captionsetup{name={Fig.},labelsep=period}
          \caption{ The numerical solution and $|\mathbf{OF}|$ with the $6$th order method on $2^6$ subdomains for the Gaussian initial value in  Example~\ref{ex2}.}   
    \label{bur_shock}
\end{figure}

\begin{rem}
    Although the third-order SSP Runge-Kutta method is used with $\Delta t=0.1\Delta x$, the sixth-order SBP stencil is observed to achieve the expected fourth-order convergence over the grid resolutions. To check the influence of the temporal discretization error, we repeat the computations with smaller time steps. The errors and convergence rates remain essentially unchanged, indicating that the temporal discretization error is not dominant in the present convergence tests. 
\end{rem}
\begin{example}[Burgers' equation in 1D]\label{ex2}
In this example, we specialize the inviscid Burgers' equation in Example~\ref{bur} to one dimension. The numerical scheme is obtained from the general semi-discrete formulation \eqref{generalform_system}. For each subdomain, the coefficient of the upwind interface term is $\alpha=\dfrac{1}{2}\left(\alpha^{(k)}+\alpha^{(k+1)}\right)$, with
 $\alpha^{(k)}=\max_{0\le i\le N}|u_i^{(k)}|$, for $k\ge 1$. 

To verify the accuracy of the proposed scheme, we use the method of manufactured solutions (MMS) to assess the errors and convergence rates. A source term is added such that the exact solution
$U(x,t)=2+0.3\sin\bigl(2\pi(x-t)\bigr)$, $ x\in[-1,1]$ satisfies the forced Burgers' equation. The time step is set to $\Delta t=0.1\Delta x$ and the $\ell^2$ errors are computed by comparing the numerical solution with the exact solution at the final time $T=2$. Fig~\ref{bur_conv_rate} shows the errors and convergence rates without and with the damping term, respectively. The observed convergence rates are as expected in both cases. In addition, the corresponding errors reported in the two subfigures become gradually close under mesh refinement. This shows that the damping term has little effect on refined meshes when the exact solution remains smooth.

 To investigate the performance of the proposed schemes after shock time, we consider the Gaussian initial condition
\begin{equation*}
U_0(x) = \exp\!\left(-\dfrac{(x-0.25)^2}{0.01}\right).
\end{equation*}
 In the simulation, the spatial domain is $\Omega=[0,1]$ and the final time is $T=0.6$. The 6th-order SBP operator is used and the time step is $\Delta t=0.1\Delta x$. Note that the exact solution remains smooth before approximately $t=0.12$, while a shock develops during the subsequent evolution. 
 We present snapshots of the numerical solutions in Fig.~\ref{bur_shock} for two cases: with both the damping and upwind interface terms, and without either term.
 At $t=0$ and $t=0.1$, the numerical solutions obtained with and without these two terms show similar evolution patterns. From the snapshots at $t=0.3$ and $t=0.6$, we can see that the numerical solution without either term develops spurious oscillations after shock time, while there are no oscillations near the discontinuity when the damping term is used. This demonstrates the good performance of the proposed algorithm.
\end{example}

\begin{figure}[!htb]     
   \centering 
  \subfigure[ Without $\mathbf{OF}$]
  {
     \includegraphics[width=0.45\textwidth]{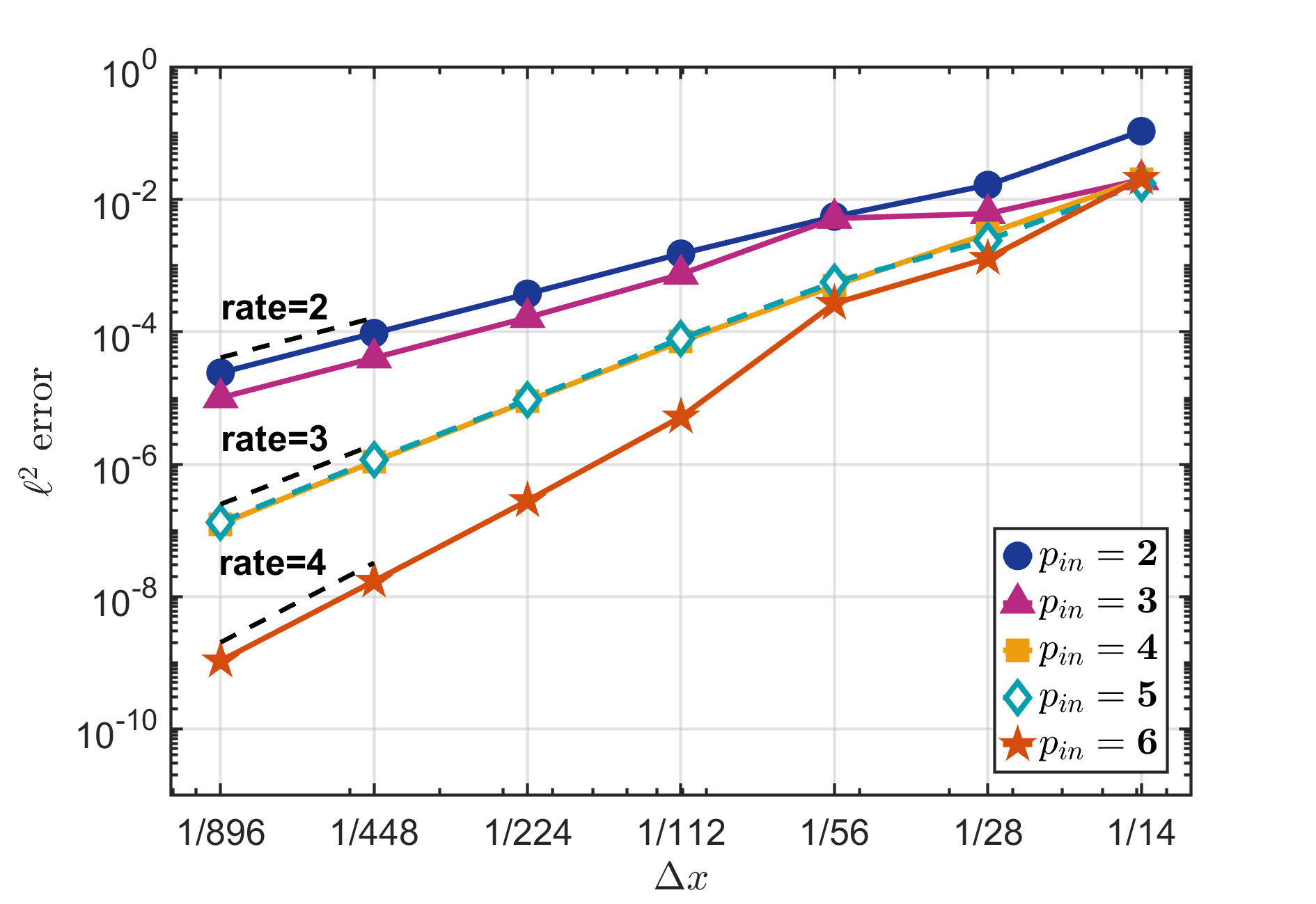}
  }\hspace{-2.2em}
  \subfigure[  With $\mathbf{OF}$]
  {
      \includegraphics[width=0.45\textwidth]{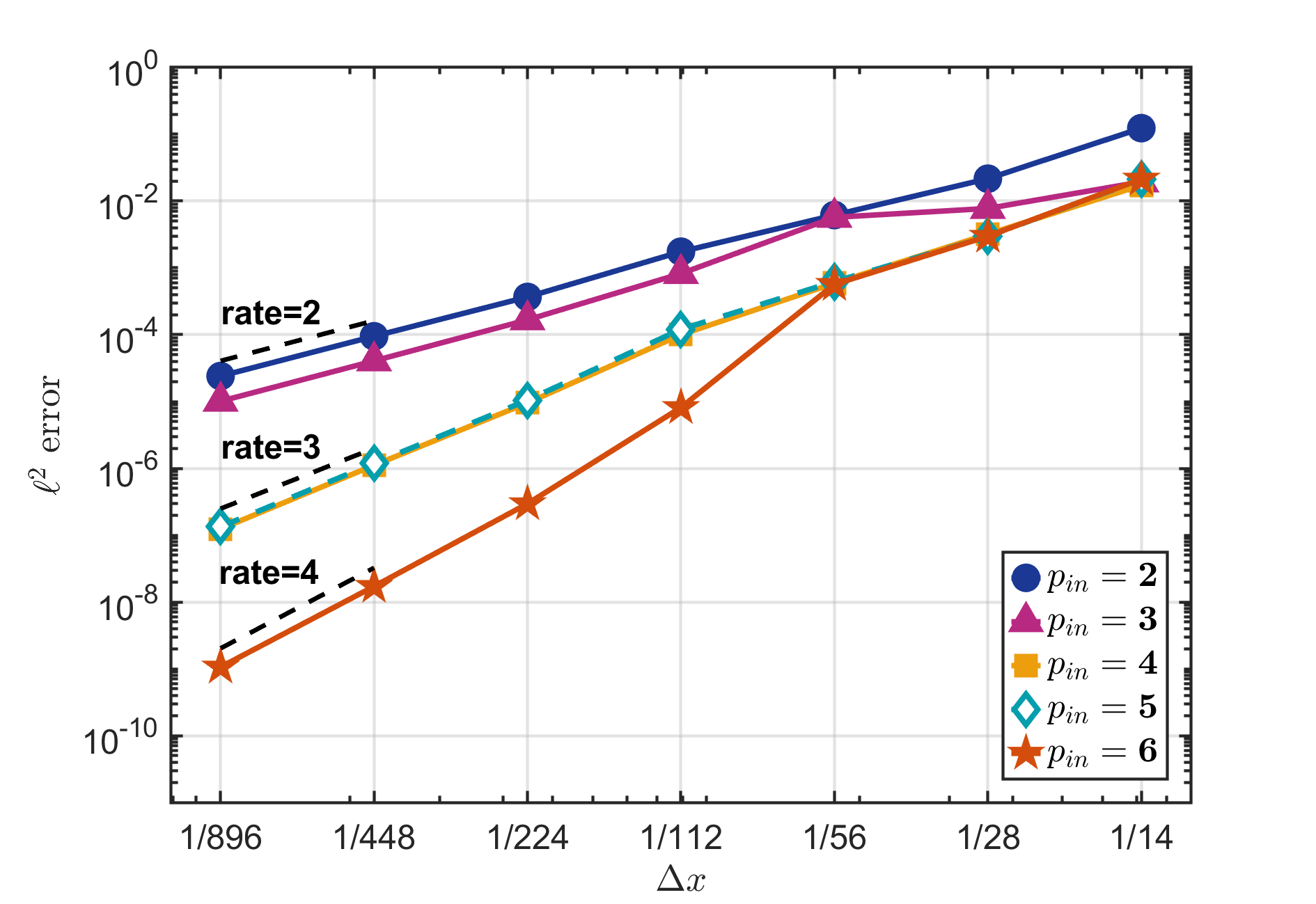}
  }
 \captionsetup{name={Fig.},labelsep=period}
  \caption{ $\ell^2$ errors and convergence rates with manufactured solutions in  Example~\ref{ex3}.}   
  \label{sw_conv_rate}
\end{figure}

\begin{figure}[!htb]     
   \centering 
   \subfigure[$\mathbf h$ at $t=2.243$ ]
  {
     \includegraphics[width=0.45\textwidth]{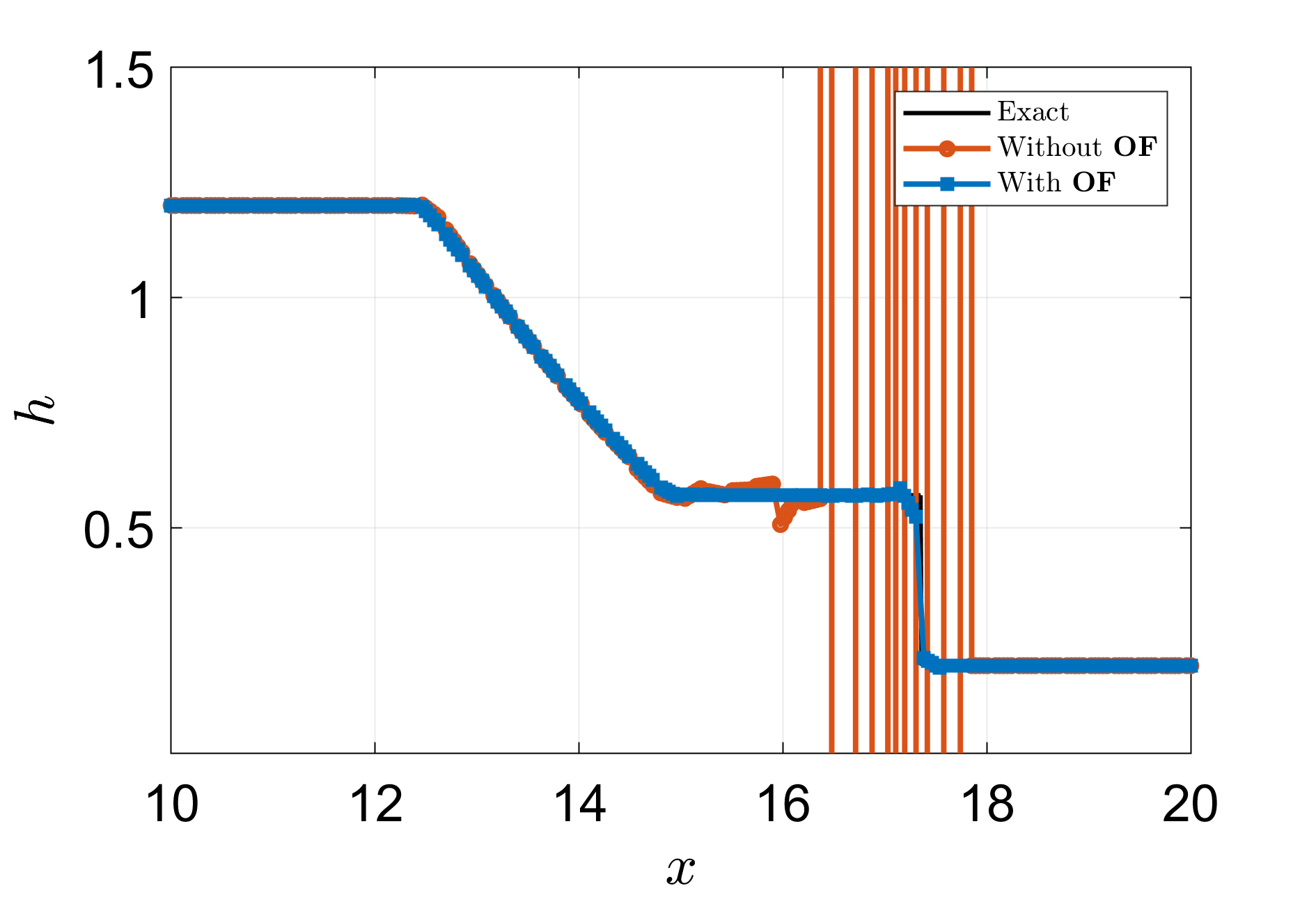}
  }\hspace{-2em}
  \subfigure[$\mathbf {hu}$ at $t=2.243$]
  {
     \includegraphics[width=0.45\textwidth]{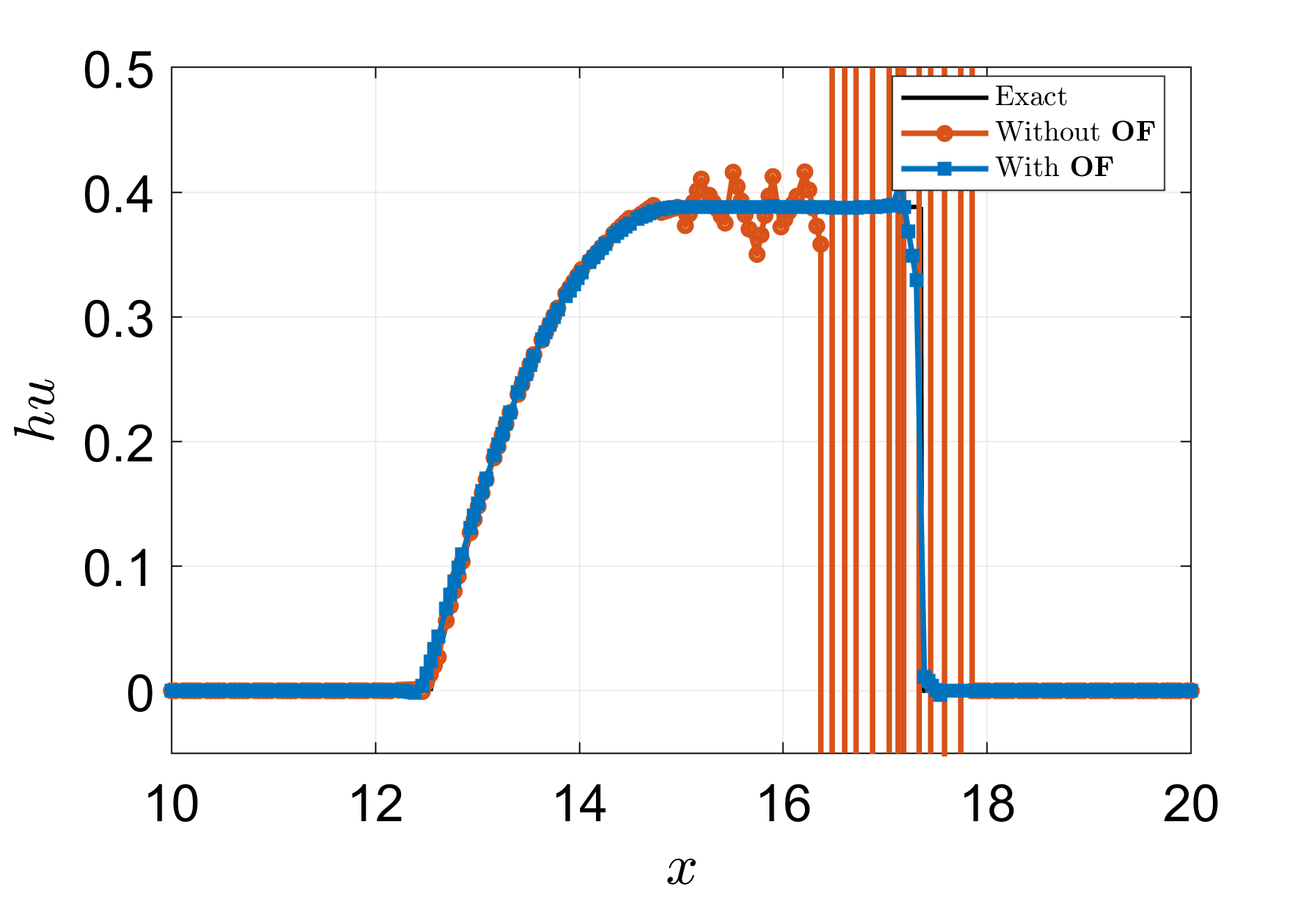}
  }
  \\[-1em]
   \centering
    \subfigure[$\mathbf h$ at $t=3.917$]  
  {
      \includegraphics[width=0.45\textwidth]{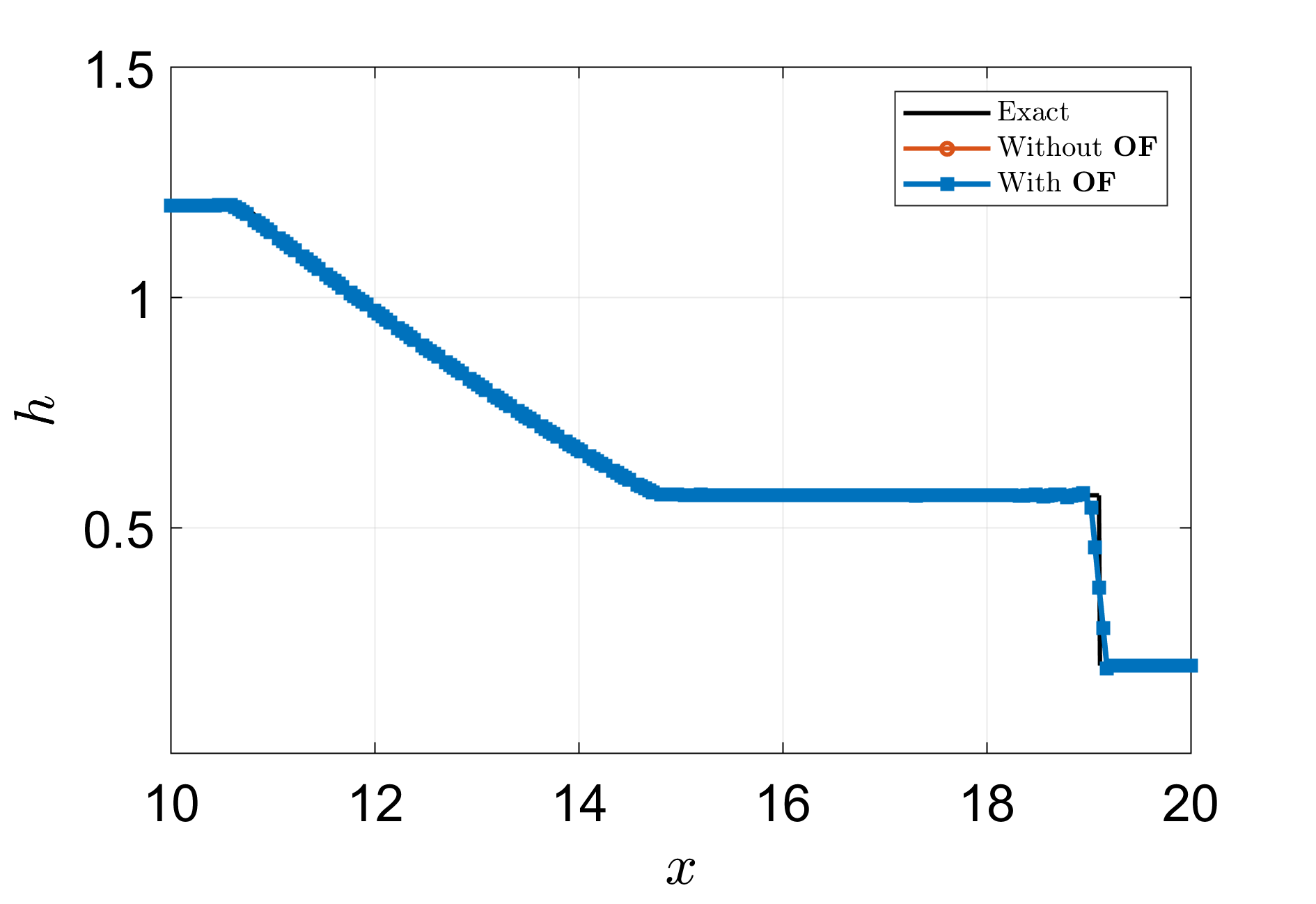}
  }\hspace{-2em}
    \subfigure[$\mathbf {hu}$ at $t=3.917$]  
  {
      \includegraphics[width=0.45\textwidth]{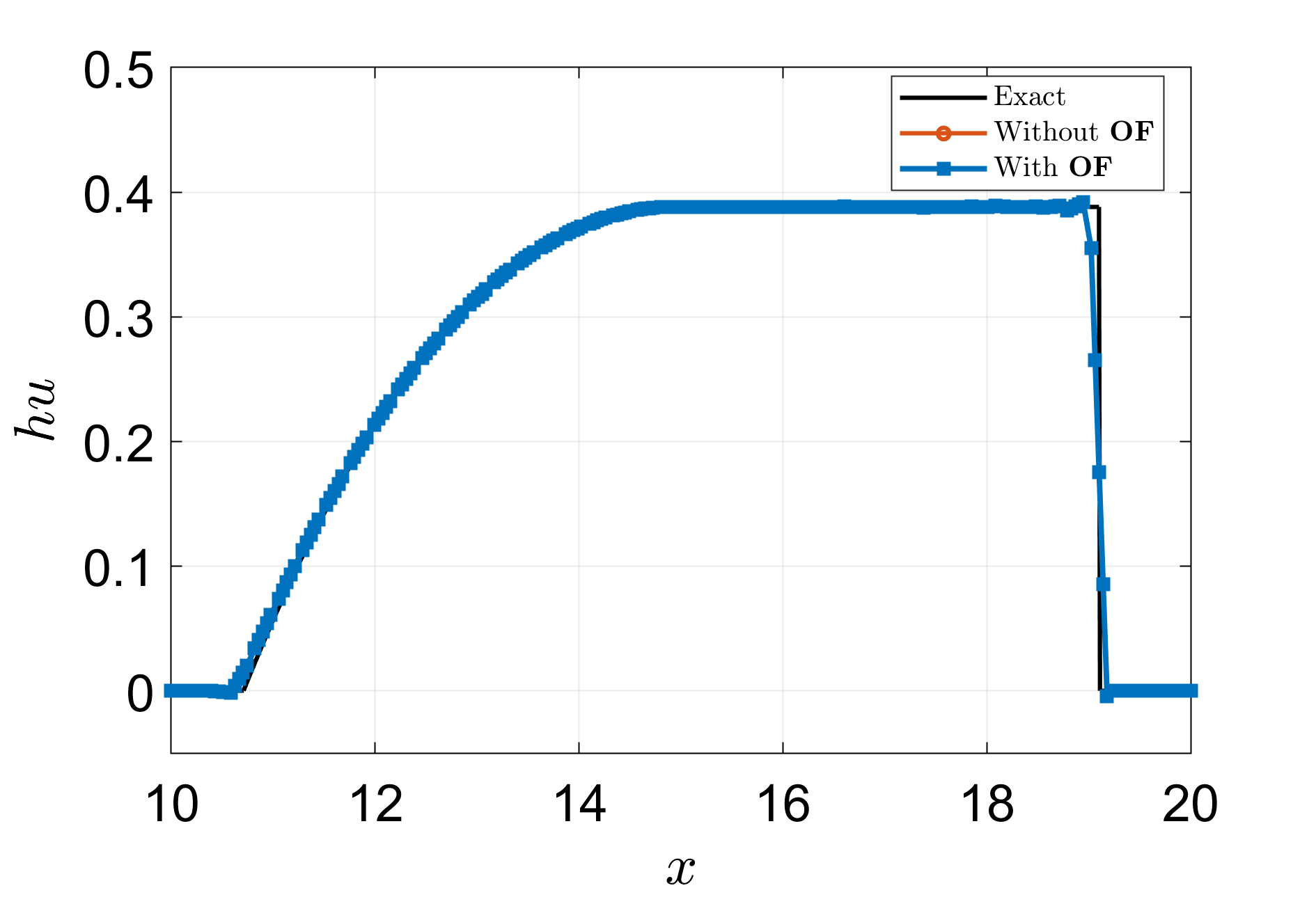}
  } 
 \captionsetup{name={Fig.},labelsep=period}
  \caption{The numerical solution with the $6$th order method for the dam  break problem with $2^8$ subdomains in  Example~\ref{ex3}.}   
  \label{sw_db}
\end{figure}

\begin{figure}[!htb]     
   \centering 
  \subfigure[Entropy]
  {
     \includegraphics[width=0.33\textwidth]{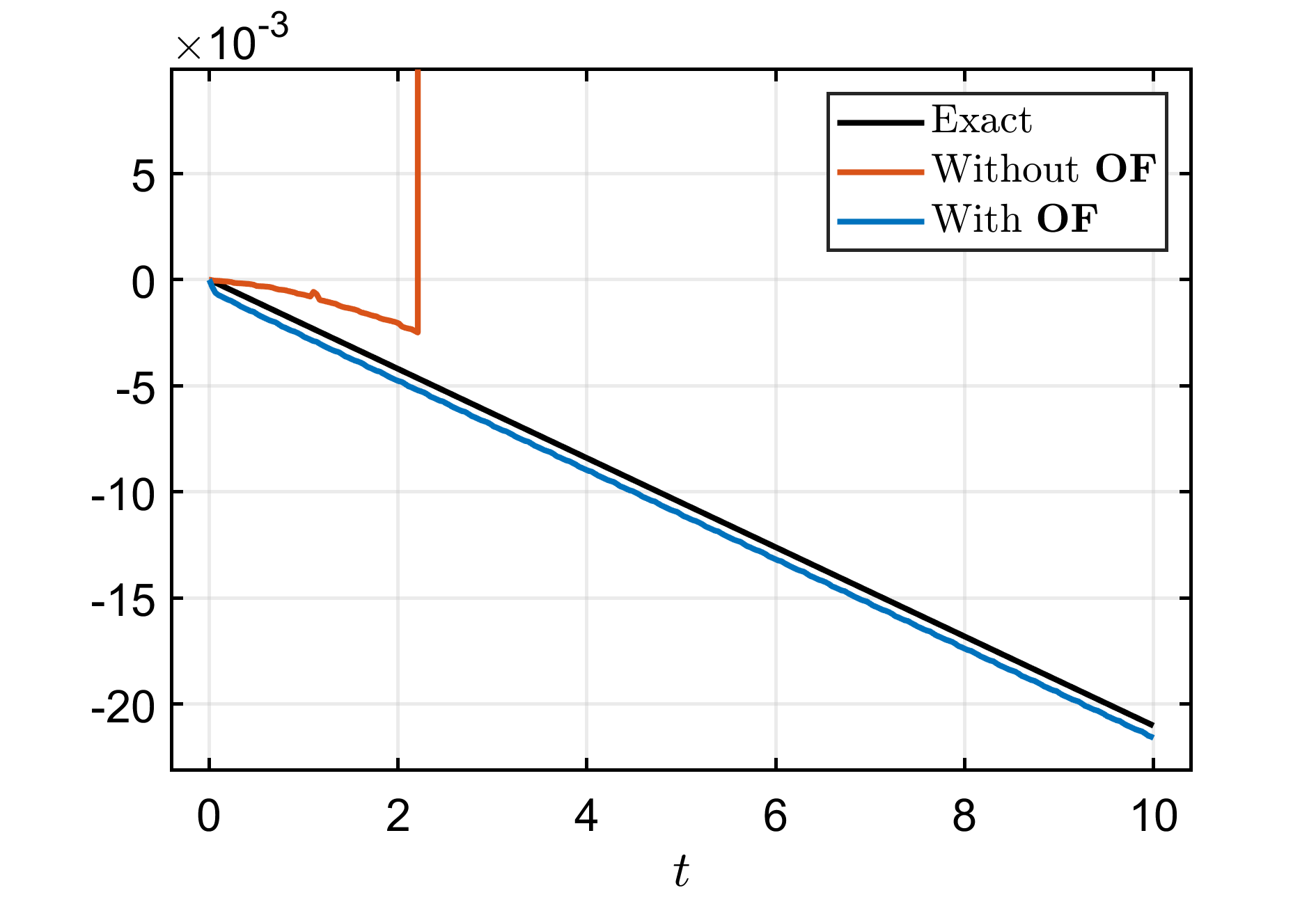}
  }\hspace{-2.2em}
  \subfigure[  Total mass]
  {
      \includegraphics[width=0.33\textwidth]{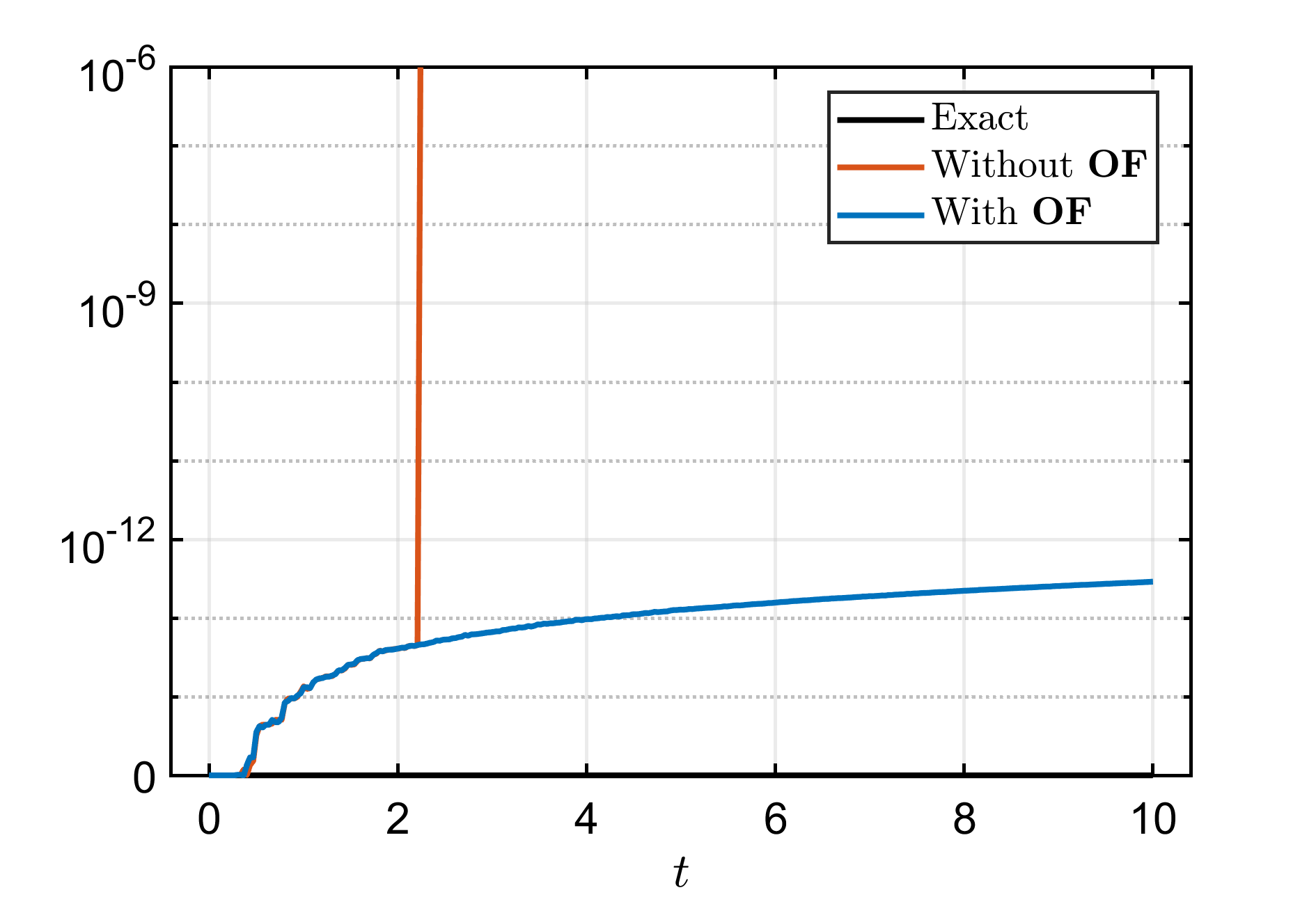} 
  }\hspace{-2.2em}
  \subfigure[ Total momentum ]
  {
      \includegraphics[width=0.33\textwidth]{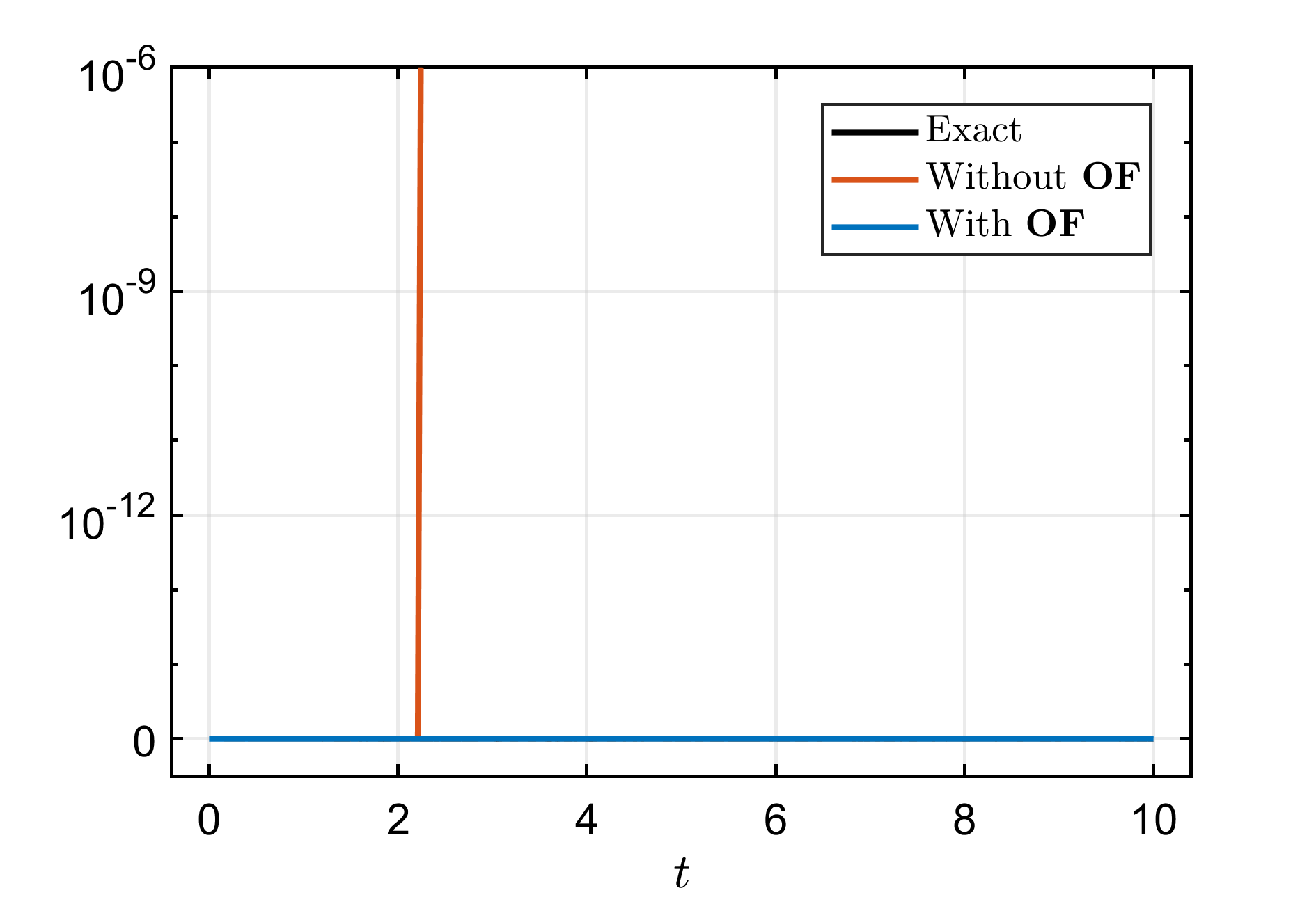}
  }
 \captionsetup{name={Fig.},labelsep=period}
  \caption{The relative change in total entropy,  mass and momentum for dam break problem in Example~\ref{ex3}.}   \label{sw_reltive_en}
\end{figure}

\begin{example}[Shallow water equations in 1D]\label{ex3}
For the shallow water equations in Example~\ref{sw}, we consider the one-dimensional problem. The corresponding semi-discrete formulation is given by Eq.~\eqref{generalform_system}. For each subdomain, the coefficient of the upwind interface term is 
$\boldsymbol{ \alpha}=\dfrac12\operatorname{diag}
([\alpha_1^{(k)}+\alpha_1^{(k+1)},\alpha_2^{(k)}+\alpha_2^{(k+1)}])$ with 
\[
\alpha_1^{(k)}=\max_{0\le i\le N}
\dfrac{h_i^{(k)}\left(|u|+\sqrt{gh_i^{(k)}}\right)}
{h_i^{(k)} g+\left(u_i^{(k)}\right)^2}
\ge 0, \quad \alpha_2^{(k)}=\max_{0\le i\le N}
\left|h_i^{(k)} u_i^{(k)}\right|\ge 0,  \quad k\ge 1.
\]
Two auxiliary variables for the adaptive damping coefficients are defined as 
\begin{equation}\label{z_1D}
    \mathbf z_1^{(k)}=\dfrac{|\mathbf u^{(k)}|+\sqrt{g\mathbf h^{(k)}}}{g},\quad \mathbf{z}_2^{(k)} =|\mathbf h^{(k)}\mathbf u^{(k)}|.
\end{equation}
To verify the accuracy of the one-dimensional shallow water equations, an appropriate forcing term is constructed to yield the exact solution
\begin{equation*}
h = 2 + 0.3\sin(2\pi(x-t)), \quad U = 2 + 0.3\sin(2\pi(x+t)), \quad
x\in [-1,1], \quad
t\in [0,2].
\end{equation*}
The time step is fixed as $\Delta t=0.01\Delta x$. 
In Fig~\ref{sw_conv_rate}, the observed convergence rates are consistent with the expected results. For smooth exact solutions, the errors with the damping term are the same order of magnitude as those without the damping term. Additionally, with all other parameters unchanged, we also test the case $\Delta t=0.1\Delta x$. The second-order SBP operator results in blow-up, whereas the other SBP operators yield the expected convergence results.

To investigate the robustness, we consider the 1D dam break problem with the initial conditions
\begin{equation*}
h_0 =
\begin{cases}
1.2, & \text{if } |x| < 15,\\
0.2, & \text{if } |x|> 15,
\end{cases}
\quad
U_0 = 0,
\quad
x \in [-30,30],
\quad
t \in [0,10].
\end{equation*}
The gravitational acceleration is set to $g=1$, the 6th-order SBP operator is used and the time step is $\Delta t=0.1\Delta x$. 
Fig.~\ref{sw_db} shows snapshots of the numerical water depth $\mathbf h$ and momentum $\mathbf {hu}$ at $t=2.243$ and $t=3.917$. 
Under the same parameter setting, the numerical solution without both upwind interface and the damping term blows up at about $t= 0.6$. Although the upwind interface term is added, the numerical solution without the damping term still blows up at around $t=2.2$. In contrast, the oscillation-free SBP-SAT scheme remains stable and can be evolved to the final time. This demonstrates the robustness of the proposed scheme. Furthermore, the SBP-SAT scheme with the OF damping is conservative and entropy consistent, as illustrated in Fig.~\ref{sw_reltive_en}.

\end{example}

\begin{example}[Lake at rest with immersed bump]\label{ex4} 
We consider the canonical lake at rest problem in 1D, modeled by the shallow water equations in Example~\ref{sw} with bathymetry, following \cite{Hew2025}. The length of the domain is $L=25$ and the bathymetry is defined as
\begin{equation}
b(x)=
\begin{cases}
0.2-0.05(x-10)^2, & \text{if } 8<x<12,\\
0, & x\in \Omega\backslash[8,12].
\end{cases}
\label{eq:bathymetry}
\end{equation}
To verify the well-balanced property of the numerical method, we prescribe the initial condition $h_0(x)+b=0.5$, $U_0(x)=0$. The final time is $T=5$ and the time step is $\Delta t=0.1\Delta x$. As described in \cite{Hew2025}, the initial condition is an exact steady-state solution. A well-balanced numerical scheme should preserve this prescribed initial state until the final time. Although the function $b(x)$ in this test are non-smooth, no damping term is used here. $h+b$ preserves at the analytical value in Fig.~\ref{wb} and the numerical errors are within the machine-precision rounding errors in Table~\ref{sw_wb}. This shows that the proposed scheme preserves the lake-at-rest steady state and verifies its well-balanced property.
\end{example}

\begin{center}
\begin{minipage}[c]{0.55\textwidth}
\centering
\refstepcounter{table}
\label{sw_wb}

{\small
Table~\thetable. Accuracy test: errors at $T=5$ for the lake at rest problem in Example~\ref{ex4}.
}

\vspace{0.4em}

\resizebox{\linewidth}{!}{
\begin{tabular}{ccccccccccc}
\hline
$K$ & $\log \ell^2_{(3rd)}$ &   &
$\log \ell^2_{(4th)}$ &  & $\log \ell^2_{(5th)}$ &   &
$\log \ell^2_{(6th)}$ &   \\
\hline
4 & 6.0489e-15 & & 2.6547e-15 & & 3.3445e-15 & & 3.2147e-15 & \\
8 & 5.3465e-15 & & 3.4693e-15 & & 4.0058e-15 & & 3.9429e-15 & \\
16 & 6.5566e-15 & & 5.6257e-15 & & 5.1499e-15 & & 4.8505e-15 & \\
32 & 8.4754e-15 & & 5.3983e-15 & & 7.4160e-15 & & 5.8078e-15 & \\
64 & 8.0712e-15 & & 8.3004e-15 & & 9.5482e-15 & & 8.0266e-15 & \\
128 & 1.0562e-14 & & 1.3236e-14 & & 1.2281e-14 & & 1.0164e-14 & \\
\hline
\end{tabular}
}

\end{minipage}
\hfill
\begin{minipage}[c]{0.44\textwidth}
\centering

\includegraphics[width=0.98\linewidth]{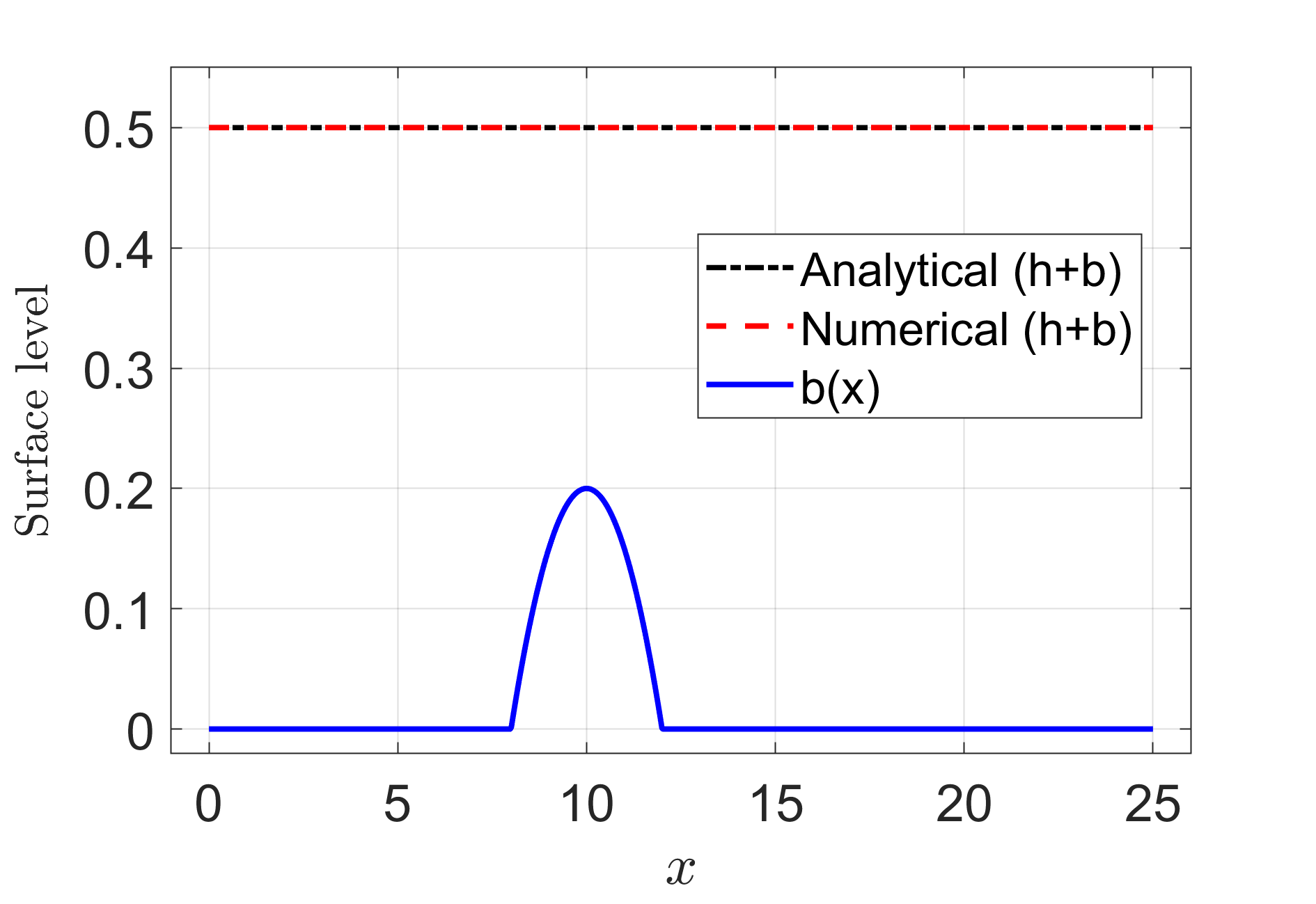}

\vspace{-0.6em}

\refstepcounter{figure}
\label{wb}
{\small
Fig.~\thefigure. Height profile for the lake at rest problem at $T=5$ with bathymetry $b(x)$ in Example~\ref{ex4}.
}

\end{minipage}

\end{center}
\subsection{Multi-dimensional problems }
In this subsection, we present numerical examples for the two-dimensional problems. For all tests, each subdomain is still divided into $N=14$ intervals in each spatial direction.
The time step is chosen as $\Delta t=\mathrm{CFL}\Delta x/(\sqrt{2}\lambda)$, where $\Delta x=\Delta y$ is the uniform grid spacing.
Here, $\lambda$ is the maximum wave speed in the medium.

\begin{figure}[!htb]     
   \centering 
  \subfigure[ Without $\mathbf{OF}$]
  {
     \includegraphics[width=0.45\textwidth]{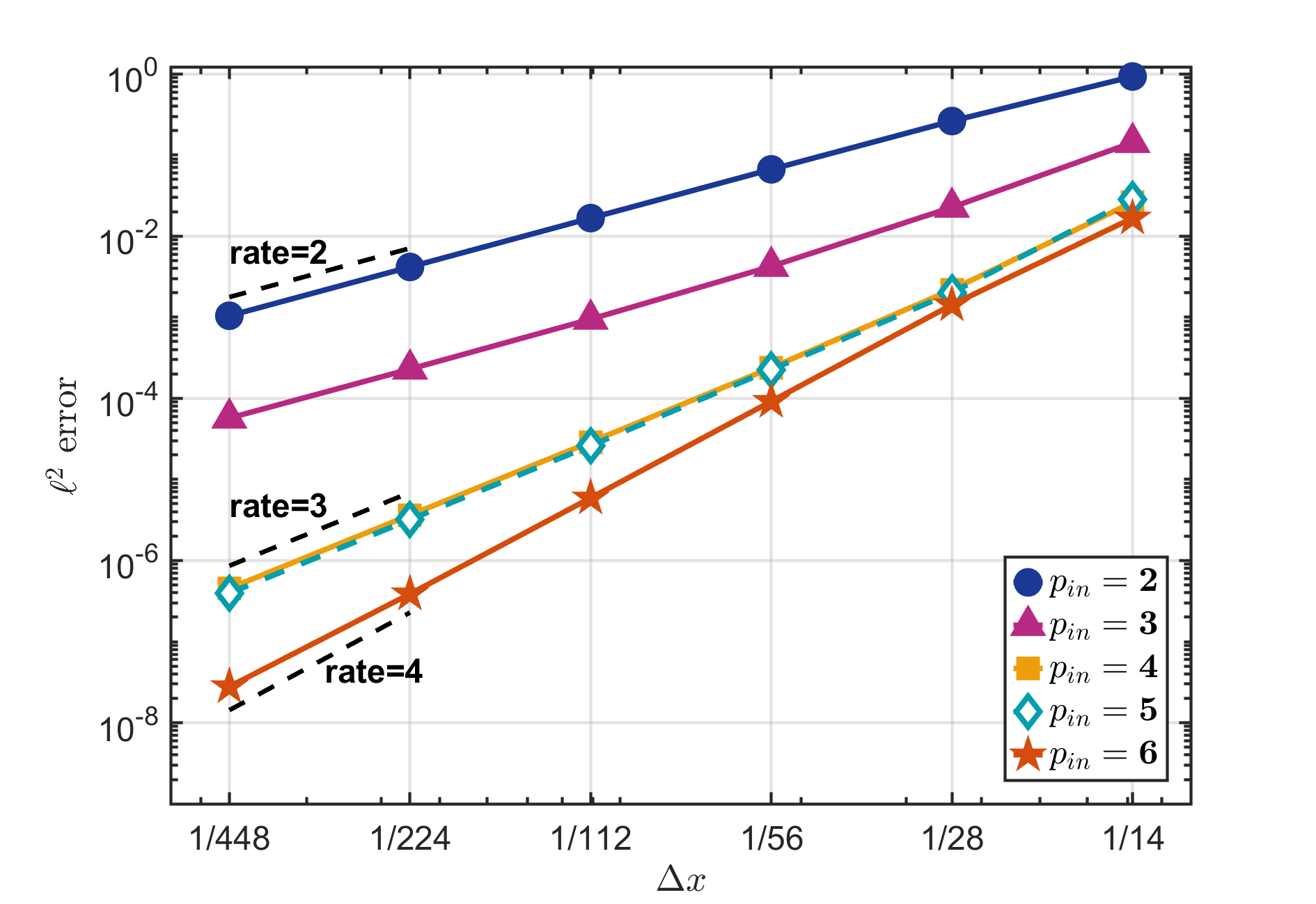}
  }\hspace{-2.2em}
  \subfigure[  With $\mathbf{OF}$]
  {
      \includegraphics[width=0.45\textwidth]{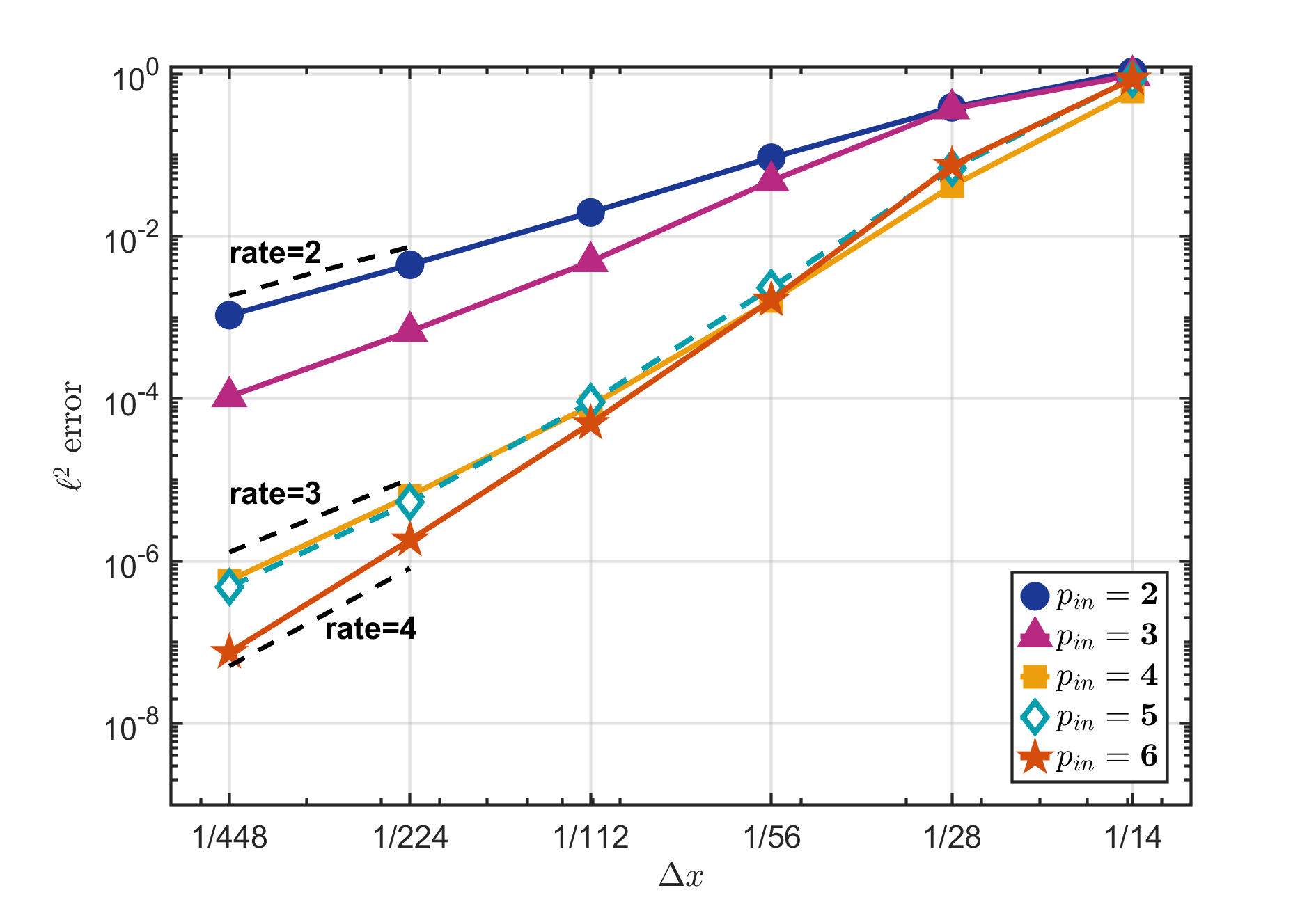}
  }
 \captionsetup{name={Fig.},labelsep=period}
  \caption{ $\ell^2$ errors and convergence rates with  initial condition (a) in  Example~\ref{ex5}.}   
  \label{adv_2D_conv_rate}
\end{figure}

\begin{figure}[!htb] 
  \centering
  \subfigure[$t=0.147$]{
    \includegraphics[width=0.45\textwidth] {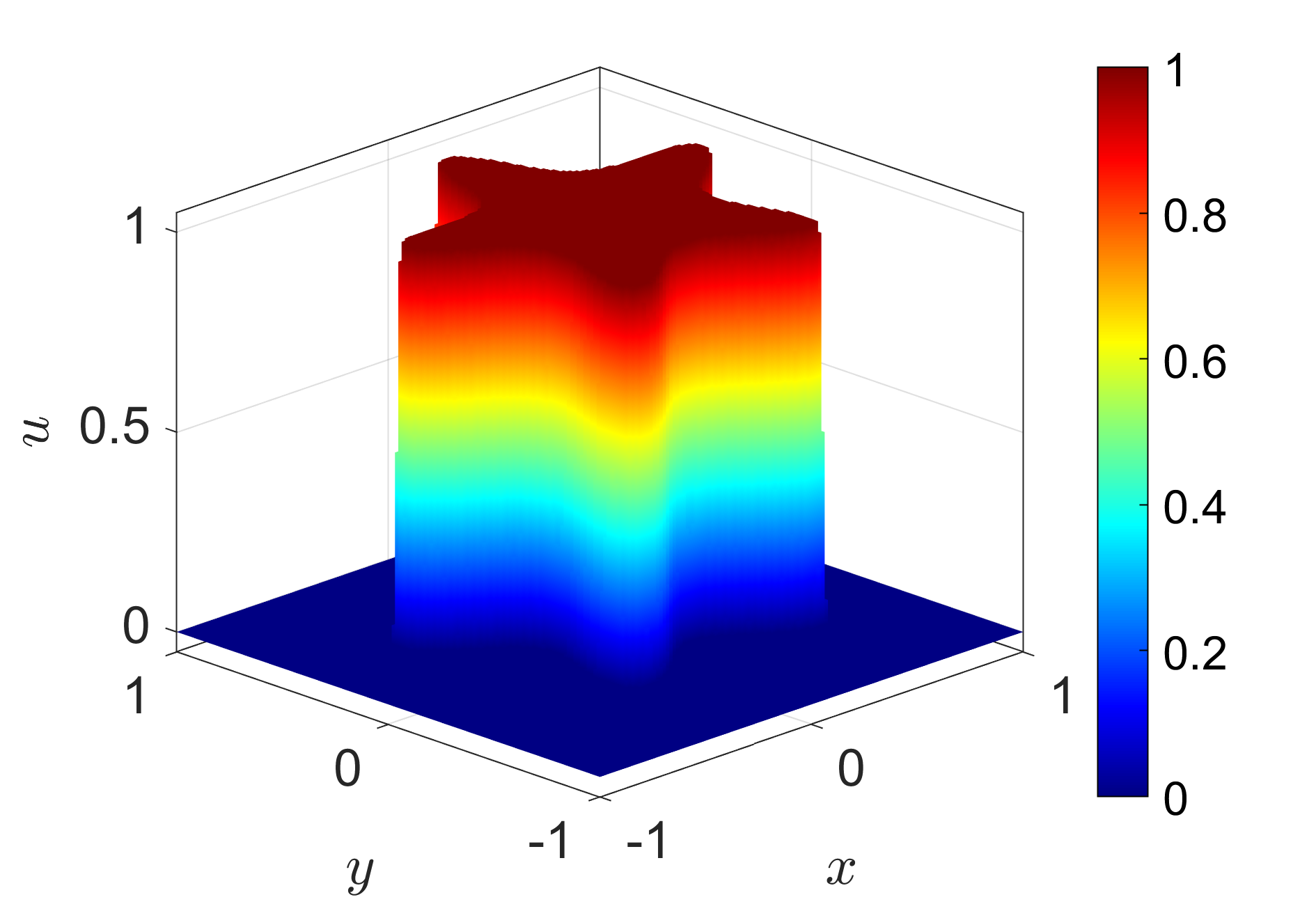}
  }
  \hspace{-2em}
  \subfigure[The cut at $y=-1/3$ for $t=0.147$]{
    \includegraphics[width=0.45\textwidth]
    {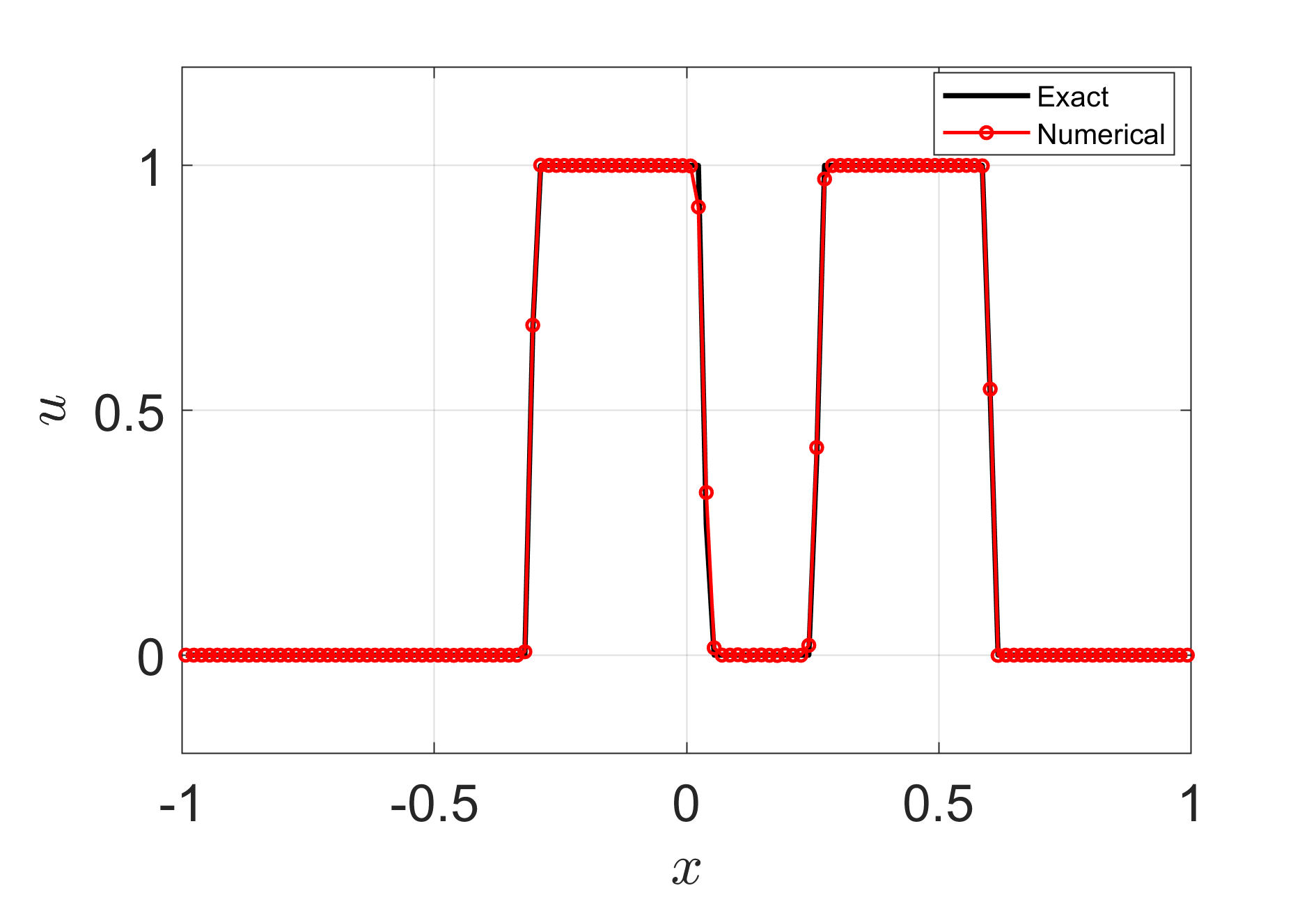}
  }
\\[-0.5em]

  \subfigure[$t=1.8$]{
    \includegraphics[width=0.45\textwidth]
    {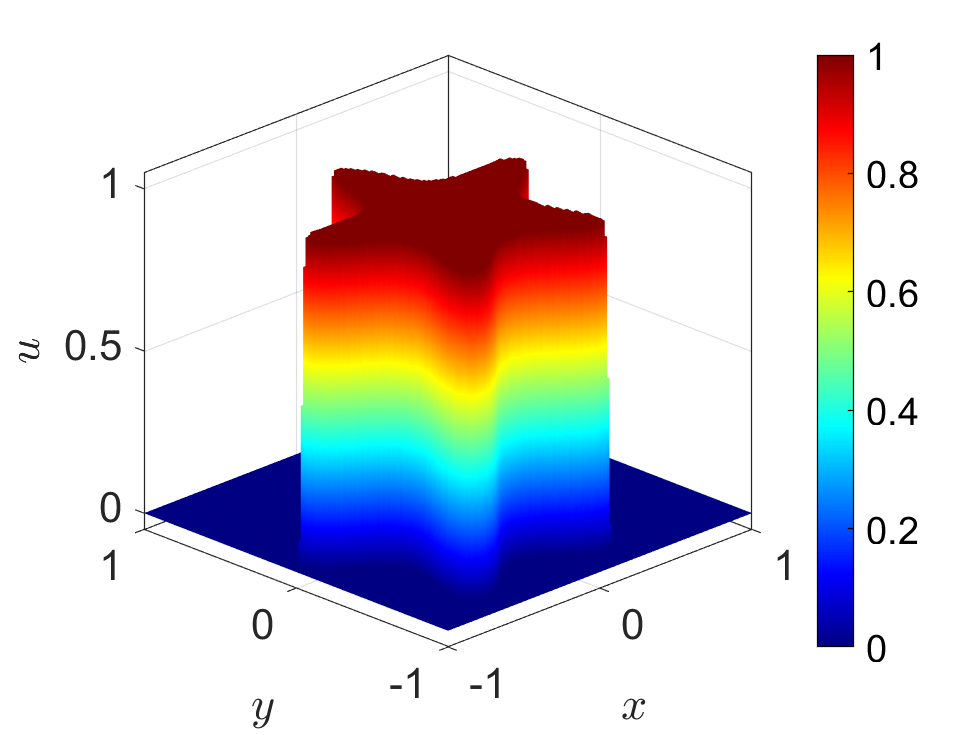}
  }
  \hspace{-2em}
  \subfigure[The cut at $y=-1/3$ for $t=1.8$]{
    \includegraphics[width=0.45\textwidth]
    {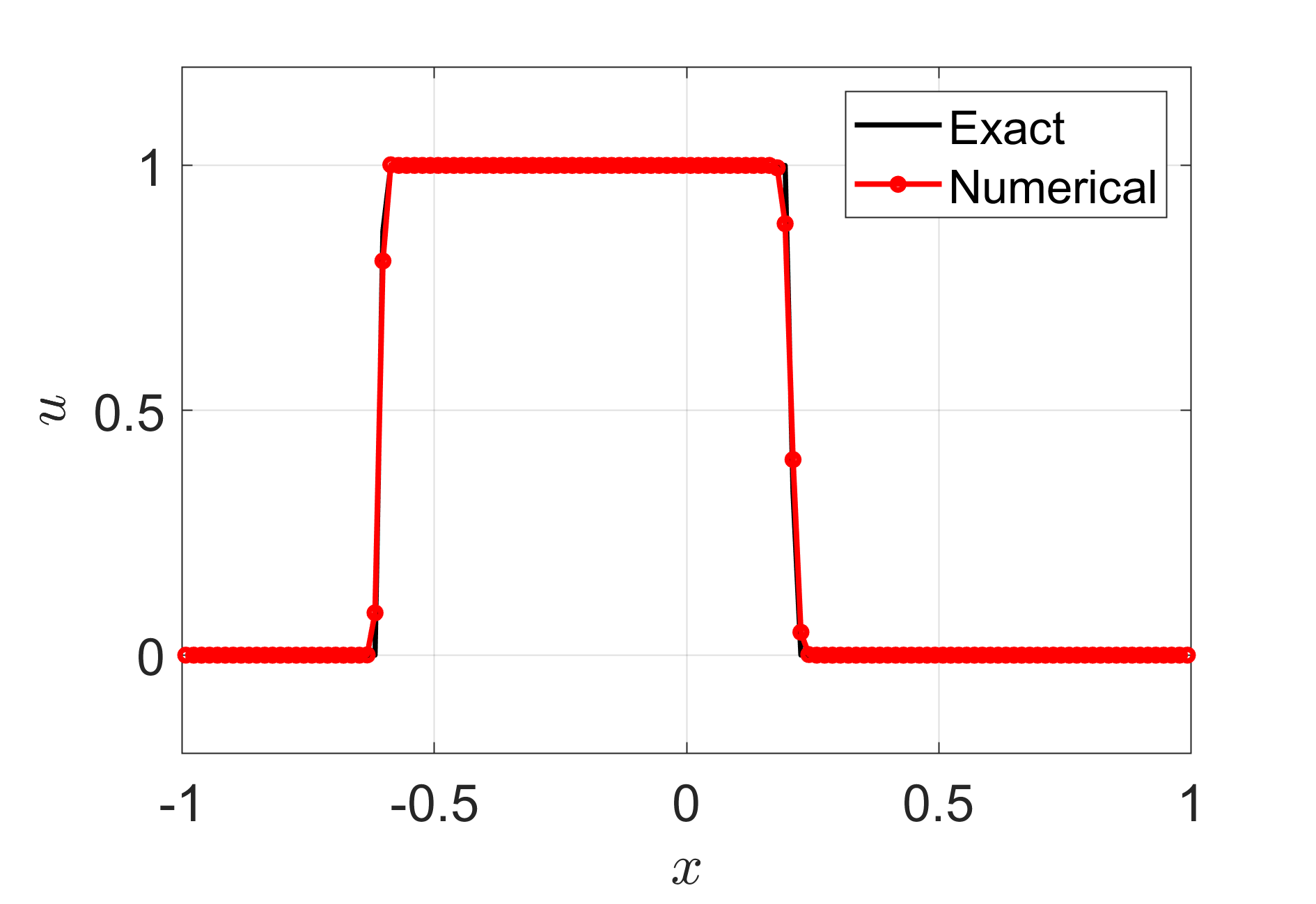}
  }

  \captionsetup{name={Fig.},labelsep=period}
  \caption{$6$th order for the initial value $U_0(x,y)$ with $2^7\times2^7$ subdomains
  in Example~\ref{ex5}.}
  \label{adv_2D_non_IC}
\end{figure}

The mesh is refined by taking the number of subdomains in each direction as $K=2^{\gamma}\times2^{\gamma}$, $\gamma=1,\ldots,6$. 
\begin{example}[Advection equation in 2D]\label{ex5}
In the advection equation \eqref{adv} in the 2D case, the numerical scheme is obtained from the generic formulation \eqref{generalform_2D}. For the upwind interface term, the coefficient is taken as in the one-dimensional case.  The following two initial conditions are considered.

(a) The smooth case:
\begin{equation*}
U_0(x,y) = \sin\big(2\pi(x+y)+0.3\big), \quad (x,y)\in[-1,1]^2.
\end{equation*}
(b) The non-smooth case :
\begin{equation*}
U_0(x,y)=
\begin{cases}
1, & (x^2+y^2)^{\dfrac12} \le \dfrac{1}{8}\left(3+3\sin(5\theta)\right),\\[6pt]
0, & \text{elsewhere},
\end{cases}
\quad (x,y)\in[-1,1]^2.
\end{equation*}
where $\theta$ is given as
\begin{equation*}
\theta =
\begin{cases}
\arccos\!\left(\dfrac{x}{\sqrt{x^2+y^2}}\right), & y\ge 0, \\[8pt]
2\pi - \arccos\!\left(\dfrac{x}{\sqrt{x^2+y^2}}\right), & y<0 .
\end{cases}
\end{equation*}
We take the final time $T=1$ for case (a) and $T=1.8$ for case (b). The time step is $\Delta t=0.1\Delta x$ in both cases.

The errors reported in Fig~\ref{adv_2D_conv_rate} are for case (a), where both the cases with and without the damping term exhibit the expected convergence rates. The results in Fig.~\ref{adv_2D_non_IC} are for case (b). We plot the contours of the numerical solution and the cuts at $y=-\dfrac{1}{3}$ with 6th-order SBP operator and $K=2^7\times 2^7 $ subdomains. In the snapshots, the contours are well captured and the spurious numerical oscillations are hardly visible in the figures.
Moreover, we have also checked the corresponding results using upwind SBP operators defined in Definition~\cite{Mattsson2017} and no significant difference was observed.
\end{example}

\begin{figure}[!htb]     
   \centering 
  \subfigure[ Without $\mathbf{OF}$]
  {
     \includegraphics[width=0.45\textwidth]{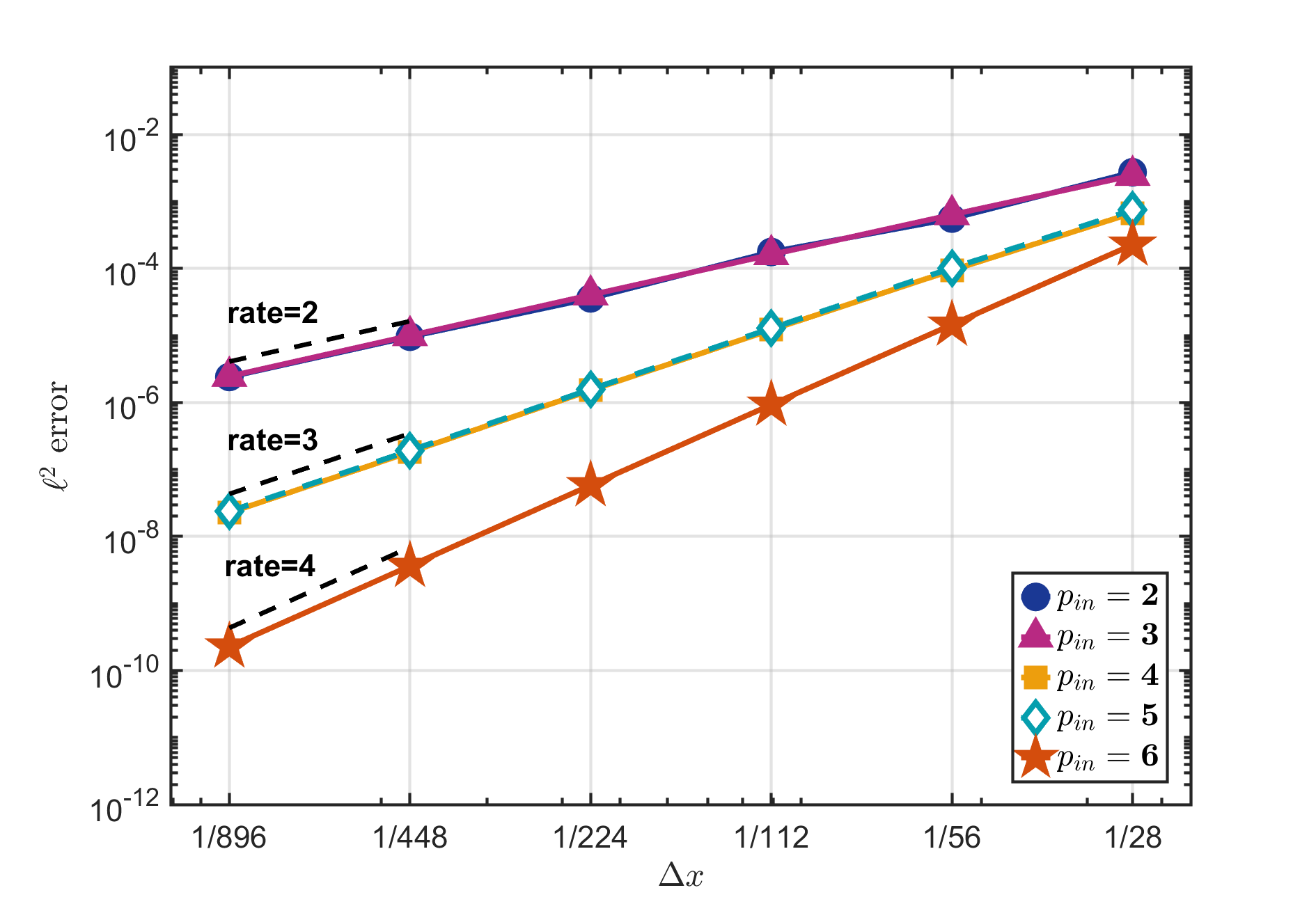}
  }\hspace{-2.2em}
  \subfigure[  With $\mathbf{OF}$]
  {
      \includegraphics[width=0.45\textwidth]{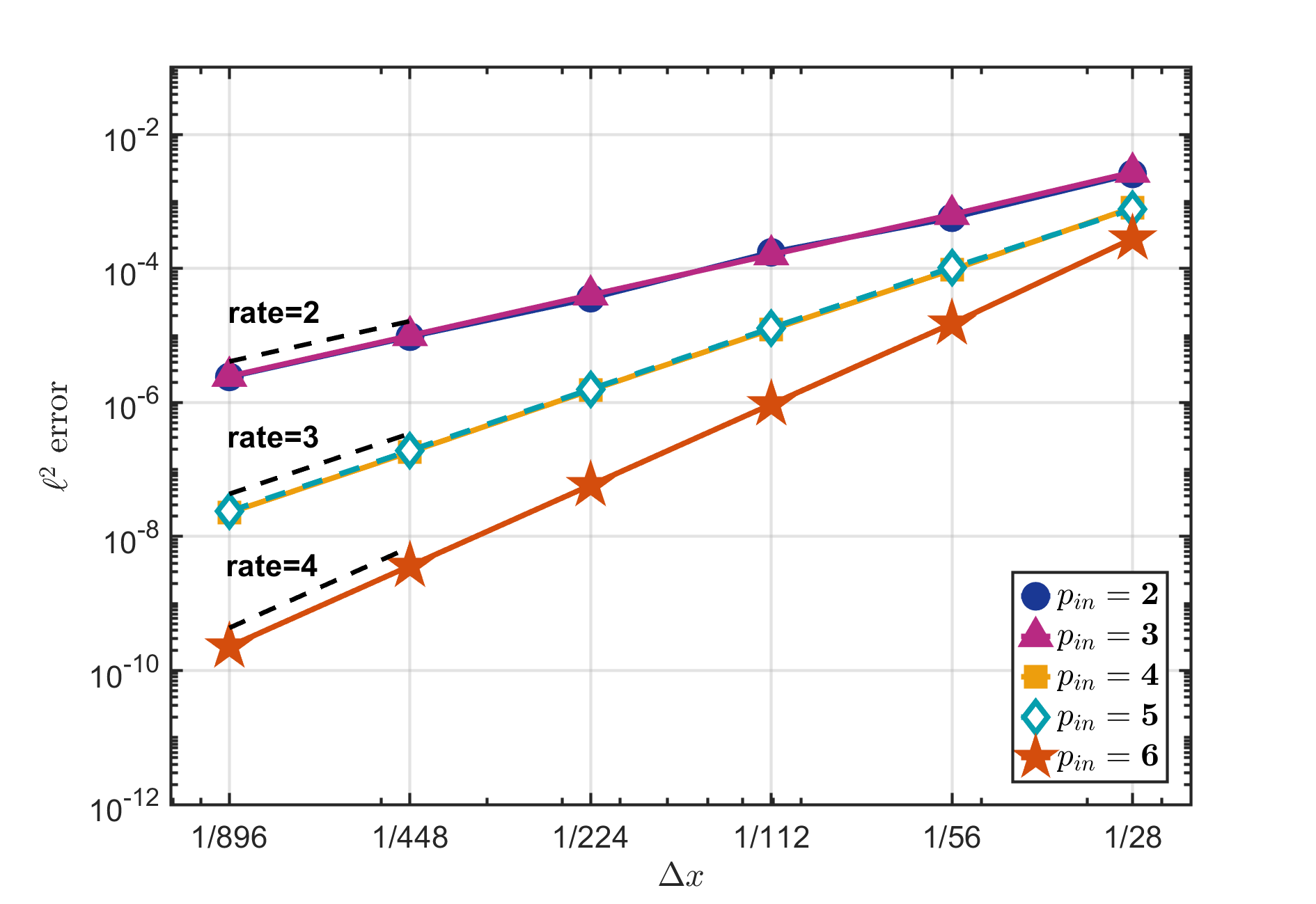}
  }
 \captionsetup{name={Fig.},labelsep=period}
  \caption{ $\ell^2$ errors and convergence rates with  manufactured solution $U(x,y,t)$ in  Example~\ref{ex6}.}   
  \label{Bur_2D_conv_rate}
\end{figure}

\begin{figure}[!htb]  
  \centering
  \subfigure[$t=0.6$]{
    \includegraphics[width=0.45\textwidth] {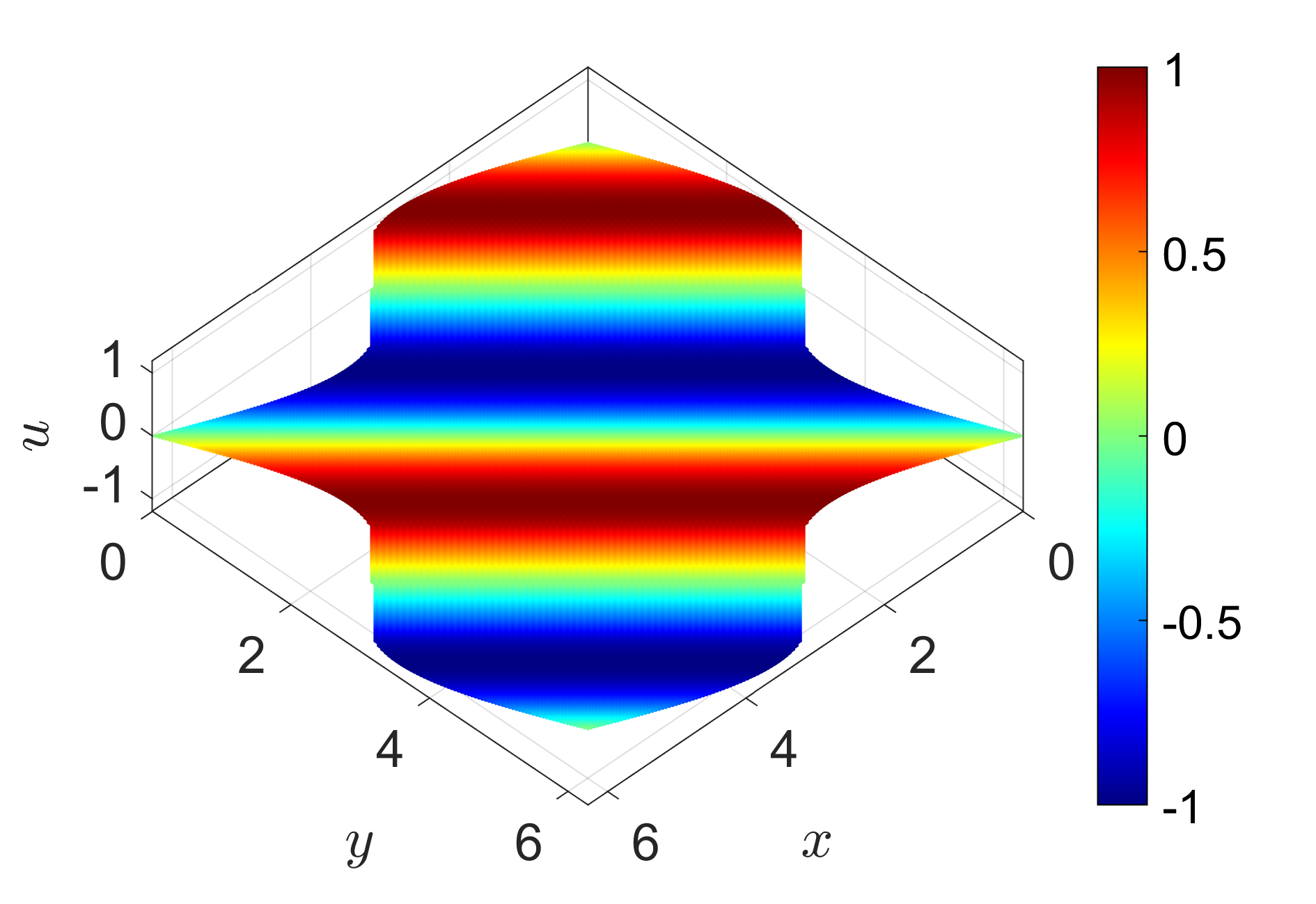}
  }
  \hspace{-2em}
  \subfigure[The cut at $y=-3\pi/4$ for $t=0.6$]{
    \includegraphics[width=0.45\textwidth]
    {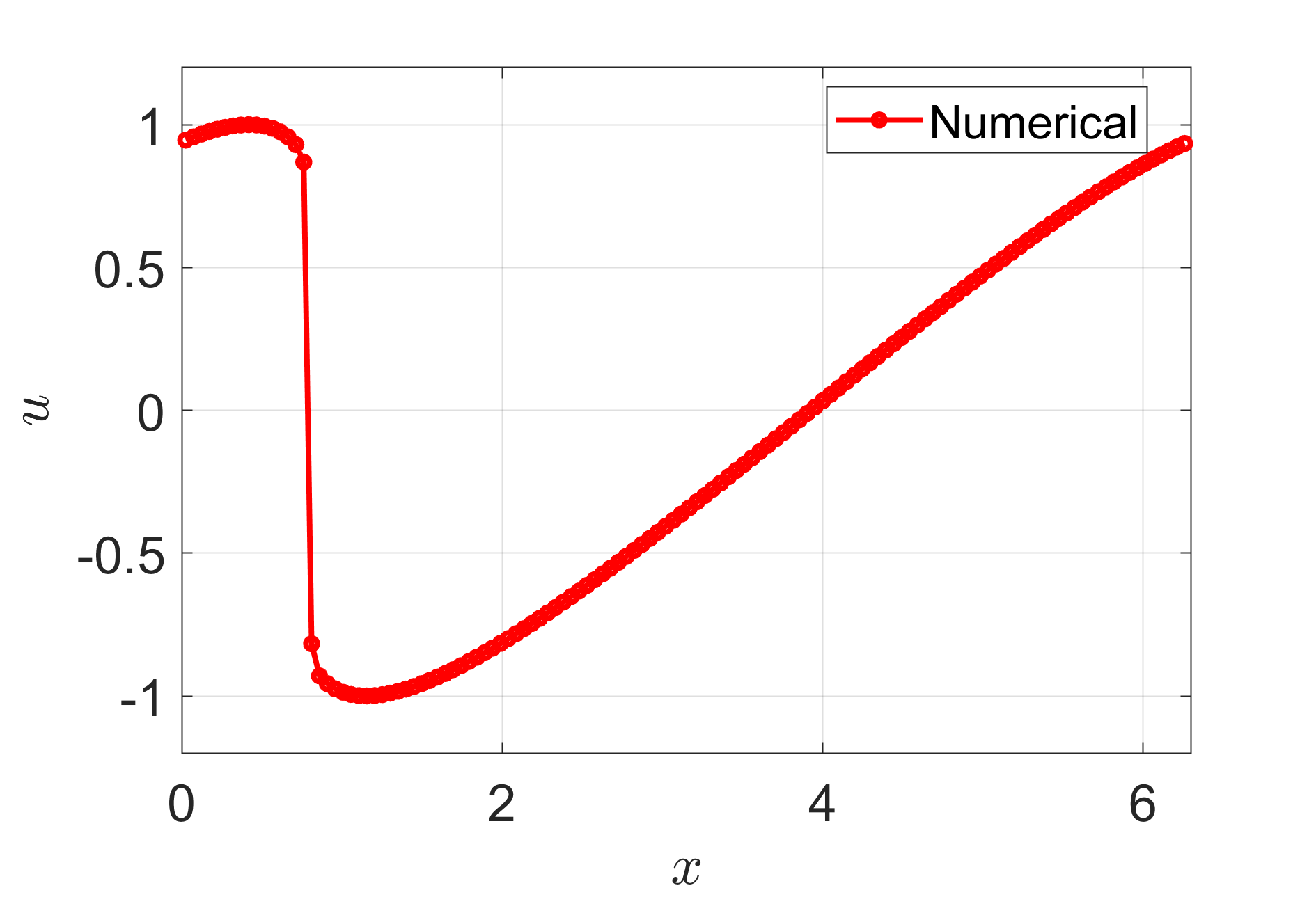}
  }
\\

  \subfigure[$t=0.8$]{
    \includegraphics[width=0.45\textwidth]
    {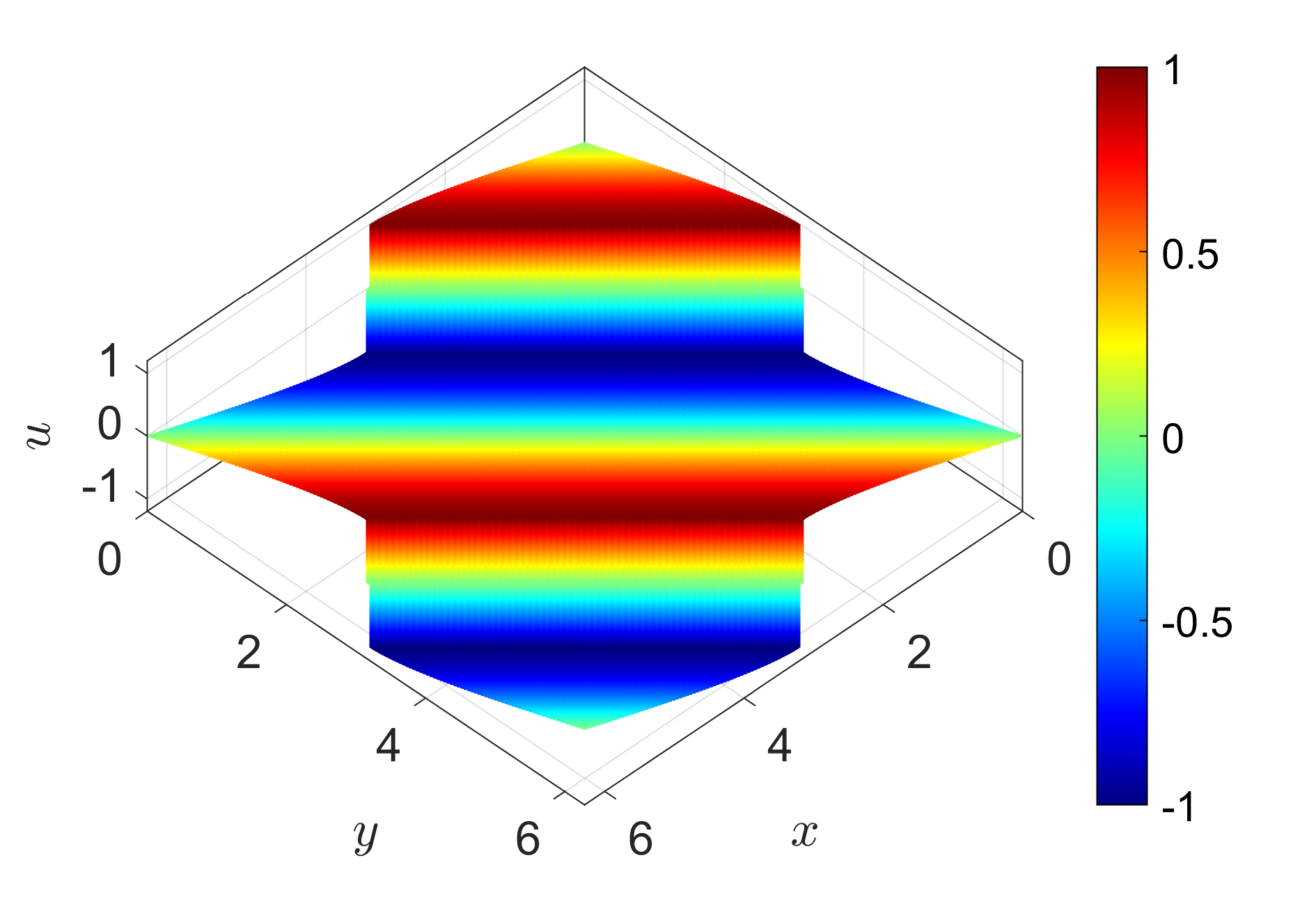}
  }
  \hspace{-2em}
  \subfigure[The cut at $y=-3\pi/4$ for $t=0.8$]{
    \includegraphics[width=0.45\textwidth]
    {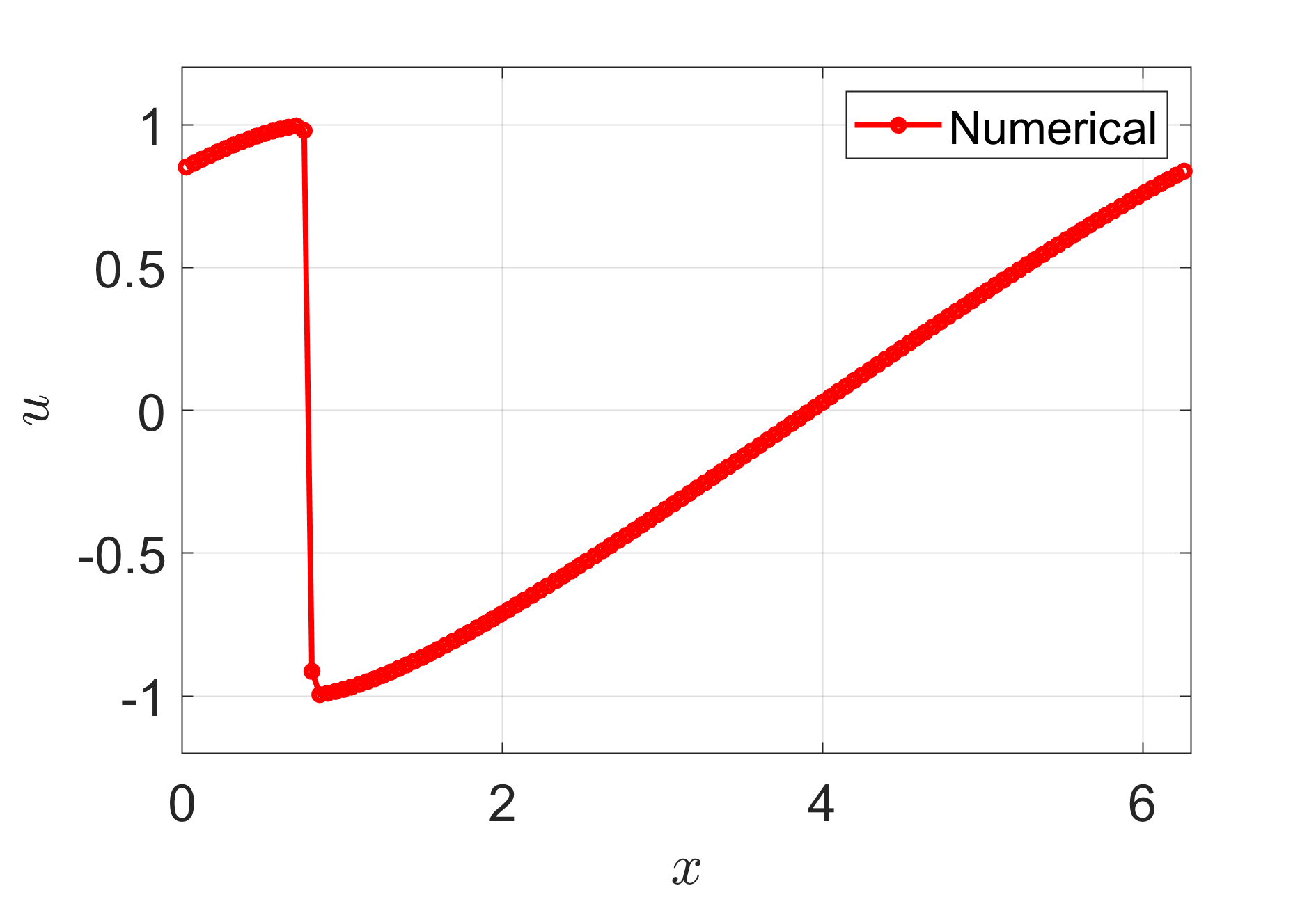}
  }

  \captionsetup{name={Fig.},labelsep=period}
  \caption{$6$th order for the  initial value (b) with $2^7\times2^7 $ subdomains in  Example~\ref{ex6}.}   
  \label{bur_2D_non}
\end{figure}

\begin{example}[Burgers' equation in 2D]\label{ex6}
In this example, we specialize the inviscid Burgers' equation in Example~\ref{bur} to 2D. The numerical scheme is obtained from the generic formulation \eqref{generalform_2D}. The coefficient of the upwind interface term is taken as in the 1D case.
To assess the convergence results, we employ MMS approach for the two-dimensional inviscid Burgers' equation. A source term is introduced such that it admits the exact solution
\begin{equation*}
U(x,y,t) = 2 + 0.3 \sin\bigl(2\pi(x + y - t)\bigr), \quad \Omega=[0,1]^2,
\end{equation*}
from which the corresponding source term is derived. 

 To investigate the performance of the proposed schemes, we consider the initial condition taken from~\cite{Lu2021}:
\[ U_0(x,y) = \sin(x+y), \quad \Omega \in [0, 2\pi]^2.\]
In convergence experiments, the time step is fixed as $\Delta t = 0.1\Delta x$ and the final time is $T = 0.8$. In Fig~\ref{Bur_2D_conv_rate}, the observed convergence rates are as expected. 
In this simulation, we use the 6th-order SBP operator with $T=0.6$, $\Delta t=0.1\Delta x$ and $K=2^7\times2^7$. As shown in Fig.~\ref{bur_2D_non}, the shock is well resolved by the proposed scheme and no spurious oscillations are observed in the cut plots at $y=-\dfrac{3\pi}{4}$.
\end{example}

\begin{figure}[!htb]     
   \centering 
  \subfigure[ Without $\mathbf{OF}$]
  {
     \includegraphics[width=0.45\textwidth]{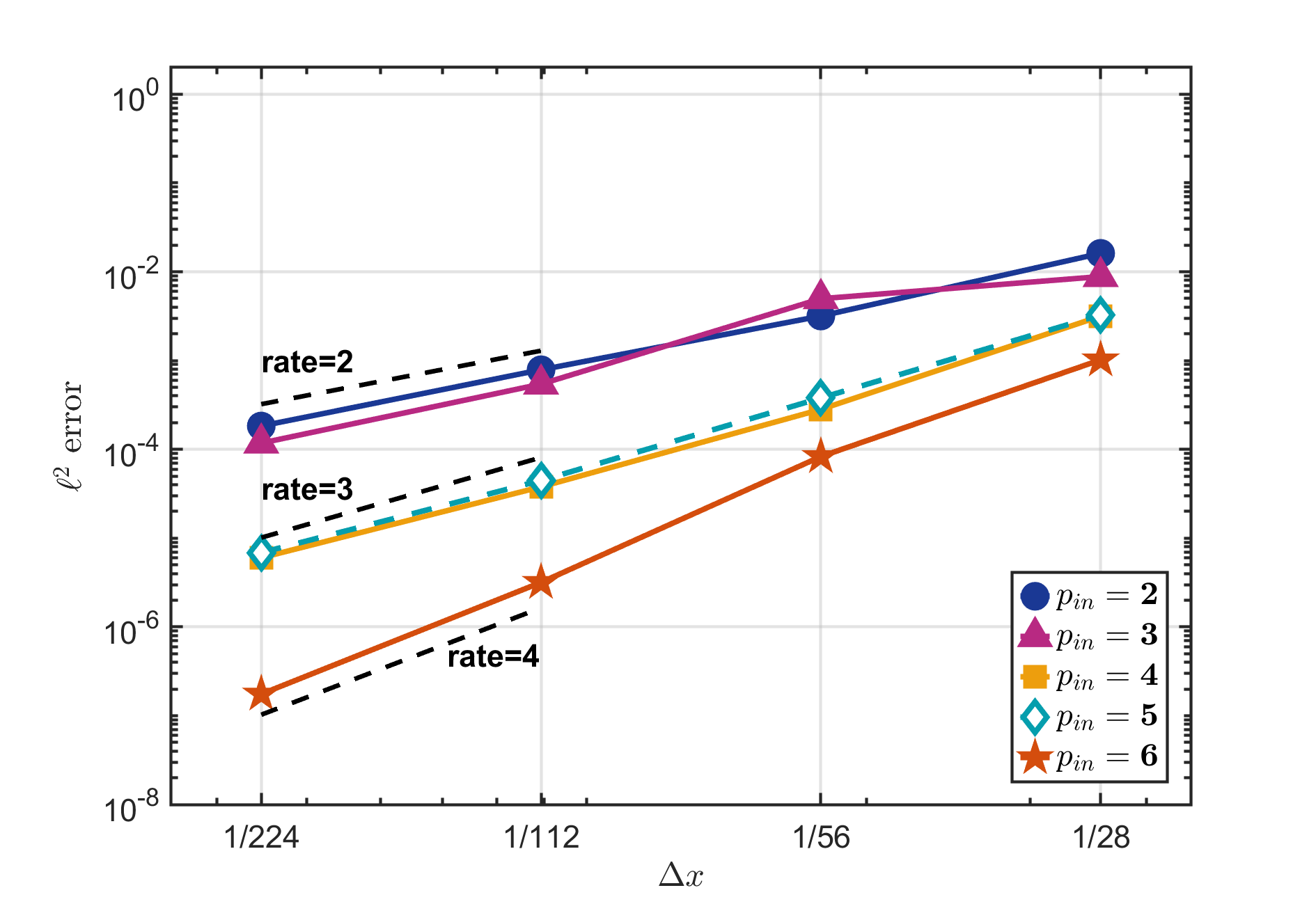}
  }\hspace{-2.2em}
  \subfigure[  With $\mathbf{OF}$]
  {
      \includegraphics[width=0.45\textwidth]{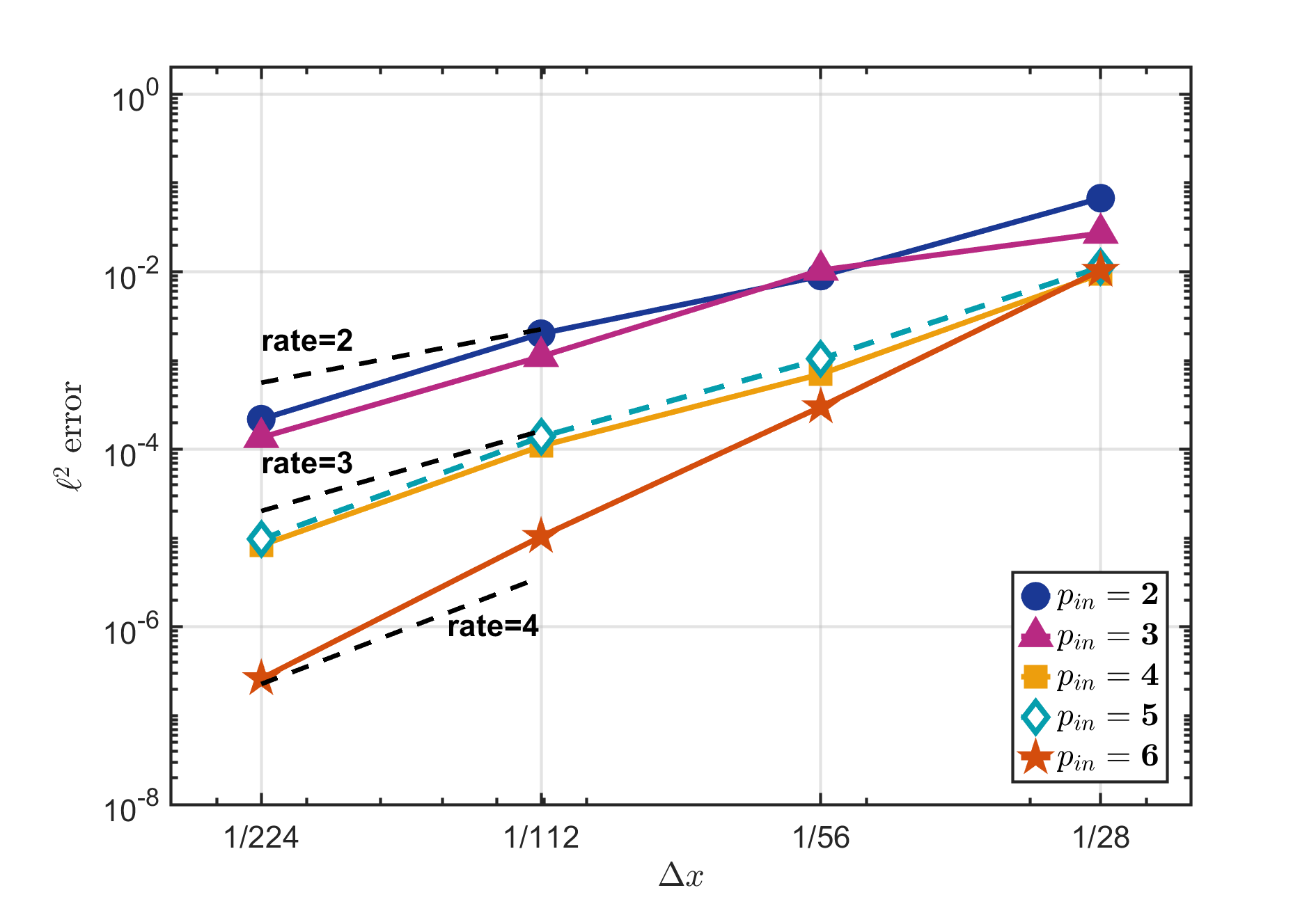}
  }
 \captionsetup{name={Fig.},labelsep=period}
  \caption{ $\ell^2$ errors and convergence rates with  manufactured solution in  Example~\ref{ex7}.}   
  \label{sw_2D_conv_rate}
\end{figure}

\begin{example}[Shallow water equations in 2D]\label{ex7}
For two-dimensional shallow water equations in Example~\ref{sw}, we consider the oscillation-free SBP-SAT scheme obtained 
 from the generic formulation~\eqref{generalform_2D}. 

For each subdomain $(x,y) \in \Omega_k$, we  define
\small{
\begin{align*}
\boldsymbol \alpha_x^{(k)}=\operatorname{diag}\left(\left[
2\max_{x,y}\dfrac{\mathbf h^{(k)}}{c}, 4\max_{ x,y}(\mathbf h^{(k)}(c-\dfrac12\sqrt{g\mathbf h^{(k)}})),
4\max_{ x,y} \mathbf h^{(k)}\sqrt{|\mathbf u^{(k)} \mathbf v^{(k)}|}
\right]\right),\\
\boldsymbol \alpha_y^{(k)}=\operatorname{diag}\left(\left[
2\max_{x,y}\dfrac{ \mathbf h^{(k)}}{c},
4\max_{x,y} \mathbf h^{(k)}\sqrt{|\mathbf u^{(k)} \mathbf v^{(k)}|},
4\max_{x,y}(\mathbf h^{(k)}(c-\dfrac12\sqrt{g \mathbf h^{(k)}}))
\right]\right),
\end{align*}}
where $c = \sqrt{(\mathbf{u}^{(k)})^2+(\mathbf{v}^{(k)})^2} + \sqrt{g\mathbf h^{(k)}}$.

For the damping term per subdomain, the entropy-variable vector is
\begin{equation*}
\mathbf g^{(k)}
=\begin{bmatrix}
\mathbf g_1^{(k)}\\
\mathbf g_2^{(k)}\\
\mathbf g_3^{(k)}
\end{bmatrix}
=\begin{bmatrix}
g\mathbf h^{(k)}
-\dfrac12\left((\mathbf u^{(k)})^2+(\mathbf v^{(k)})^2\right)\\
\mathbf u^{(k)}\\
\mathbf v^{(k)}
\end{bmatrix},
\end{equation*}
and the variables used in the damping coefficients are chosen as
\begin{equation}\label{z_2D}
\mathbf z_1^{(k)}
=\dfrac{\sqrt{(\mathbf u^{(k)})^2+(\mathbf v^{(k)})^2}
+\sqrt{g\mathbf h^{(k)}}}{g},\quad
\mathbf z_2^{(k)}=|\mathbf h^{(k)}\mathbf u^{(k)}|,\quad
\mathbf z_3^{(k)}=|\mathbf h^{(k)}\mathbf v^{(k)}|.
\end{equation}

To verify the convergence rates, we employ the MMS for the two-dimensional shallow water equations and the exact solution should satisfy
\begin{align*}
&h = 2 + 0.2\sin\bigl(2\pi(x - t)\bigr)\sin\bigl(2\pi(y - t)\bigr),\\
&U=V = 2 + 0.2\sin\bigl(2\pi(x + t)\bigr)\sin\bigl(2\pi(y + t)\bigr).
\end{align*}

We take the time step $\Delta t=0.01\Delta x$.
The errors are evaluated at the final time $T=1$ by comparing the numerical and exact solutions on domain $\Omega=[0,1]^2$. Fig~\ref{sw_2D_conv_rate} reports the corresponding convergence rates. The convergence results are again  consistent with the expected orders.
\end{example}

\begin{example}[Merging vortices]\label{ex8}
Next, we consider the merging vortices test for the two-dimensional nonlinear shallow water equations, following \cite{Hew2025, McRae2014}. The initial condition consists of a pair of Gaussian vortices function with the incompressible stream function
\[
\psi=\psi_{+}+\psi_{-},
\]
where
\[
\psi_{\pm}(x,y)=\exp\bigl(-2.5[(
x-\dfrac{3.05\pm0.45}{3}\pi)^2+(y-\pi)^2]
\bigr).
\]
 To ensure that the initial state is geostrophically balanced, we take
\[
h_0=H+\dfrac{f}{g}\psi(x,y), \quad U_0=\partial_y\psi(x,y),\quad
V_0=\partial_x\psi(x,y).
\]
The parameters are chosen as $f=g=5$ and $H=8$.  

 For each subdomain, we choose the upwind coefficients as $\frac{2}{5}\boldsymbol{\alpha}_x^{(k)}$ and $\frac{2}{5}\boldsymbol{\alpha}_y^{(k)}$ in the $x$- and $y$-directions, respectively.
The computational domain is $\Omega=[0,2\pi]^2$
 and divided into $K=2^5\times2^5$ uniform subdomains. The sixth-order SBP operator is used for the spatial discretization and the numerical solution is evolved until the final time $T=20$. The time step is $\Delta t=0.01\Delta x$.
We compare the numerical solutions obtained without the damping term and with the proposed damping term.
To examine the evolution of the vortical structures, we plot the potential vorticity $\dfrac{\omega}{h}$
with the absolute vorticity $\omega=\partial_xV-\partial_yU+f$.
We also monitor the relative changes in the total entropy and total enstrophy during the simulation.
The total enstrophy is defined by $E_{\mathrm{n}}=
\int_{\Omega}\frac{\omega^2}{h}\,\mathrm{d}x\mathrm{d}y $.
\begin{figure}[!htb]  
  \centering
  \subfigure{
    \includegraphics[width=0.45\textwidth] {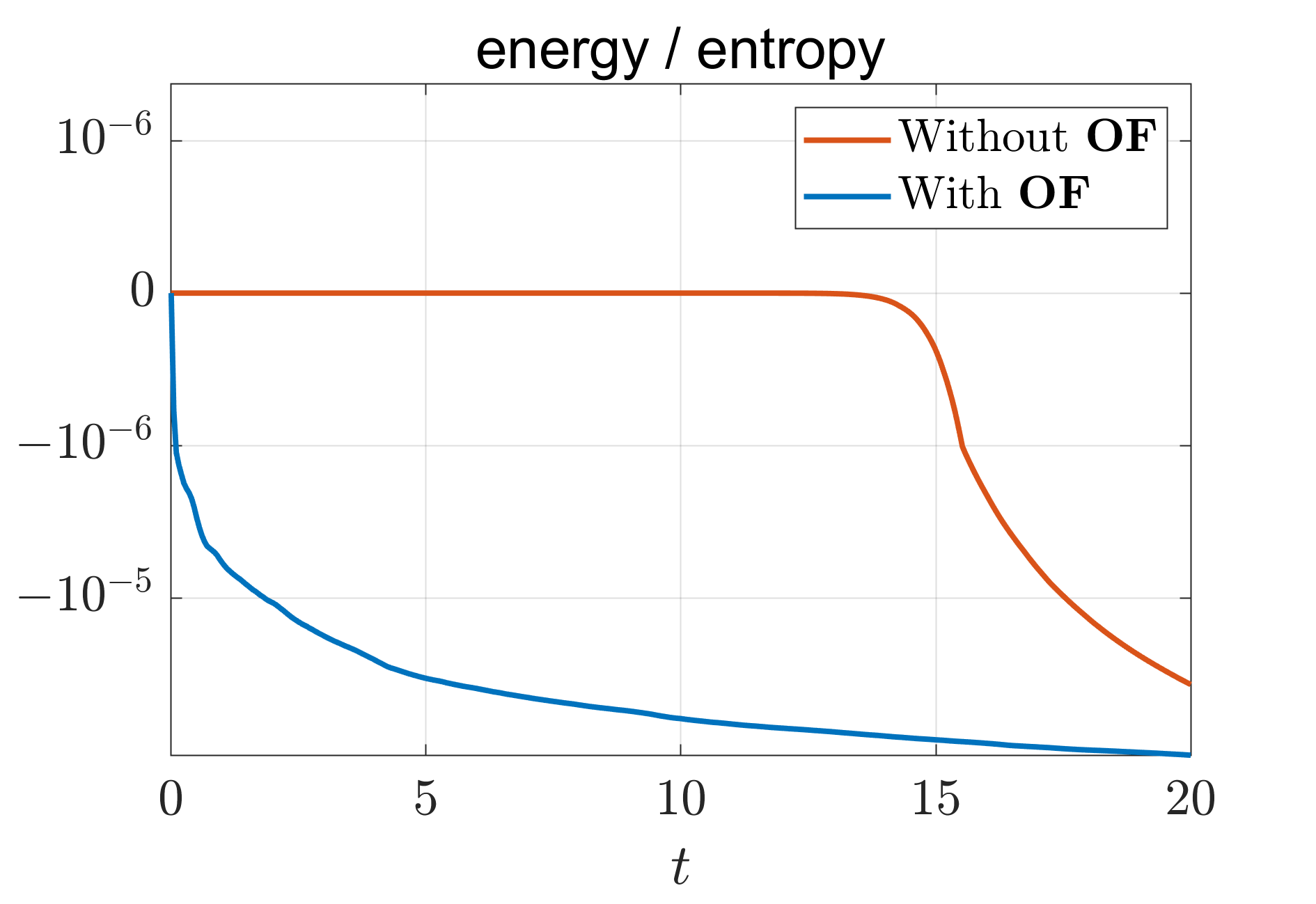} }
  \hspace{-1.5em}
  \subfigure{
    \includegraphics[width=0.45\textwidth]
    {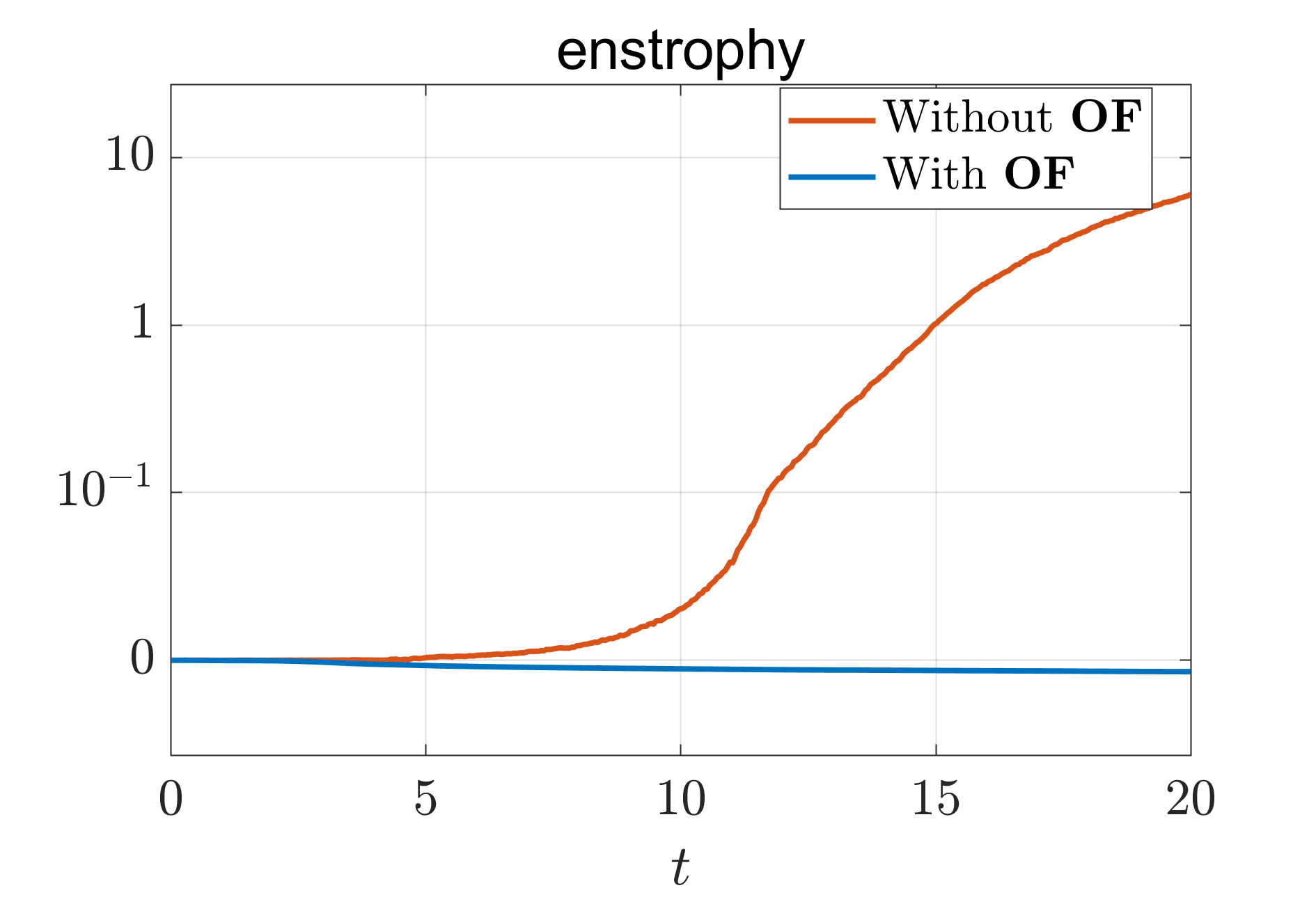}
  }
  \captionsetup{name={Fig.},labelsep=period}
  \caption{The relative change of Entropy and Enstrophy in Example~\ref{ex8}.}   
  \label{sw_2D_enes1}
\end{figure}

\begin{figure}[!htb]  
\centering
\begin{minipage}[c]{1\textwidth}
\centering

\subfigure[ $t=2$ without $\mathbf{OF}$]{
  \includegraphics[width=0.35\linewidth, trim=3cm 0.1cm 1.5cm 1.5cm, clip]
  {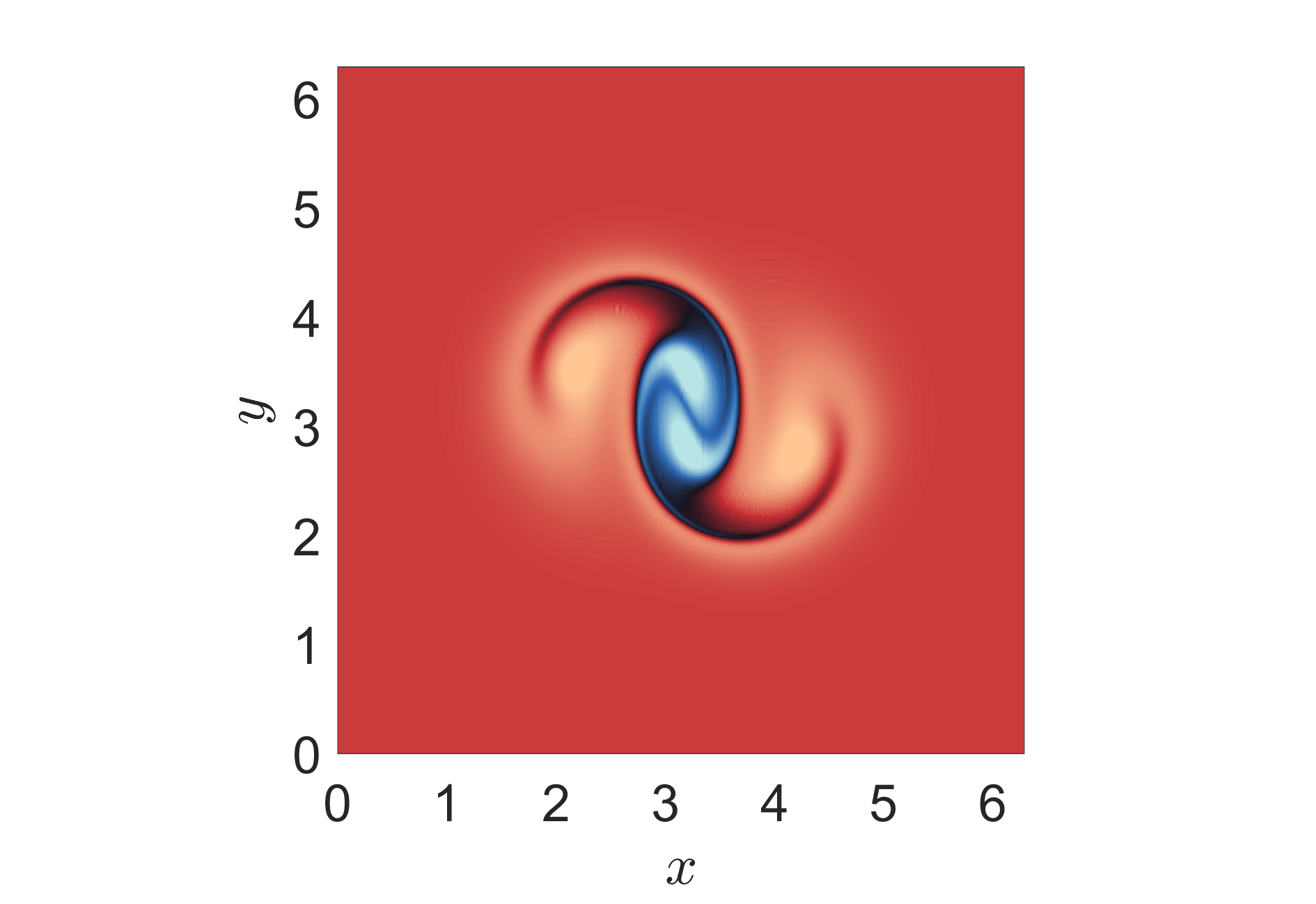}
}
\hspace{-3em}
\subfigure[ $t=2$ with $\mathbf{OF}$]{
  \includegraphics[width=0.35\linewidth,trim=3cm 0.1cm 1.5cm 1.5cm, clip]
  {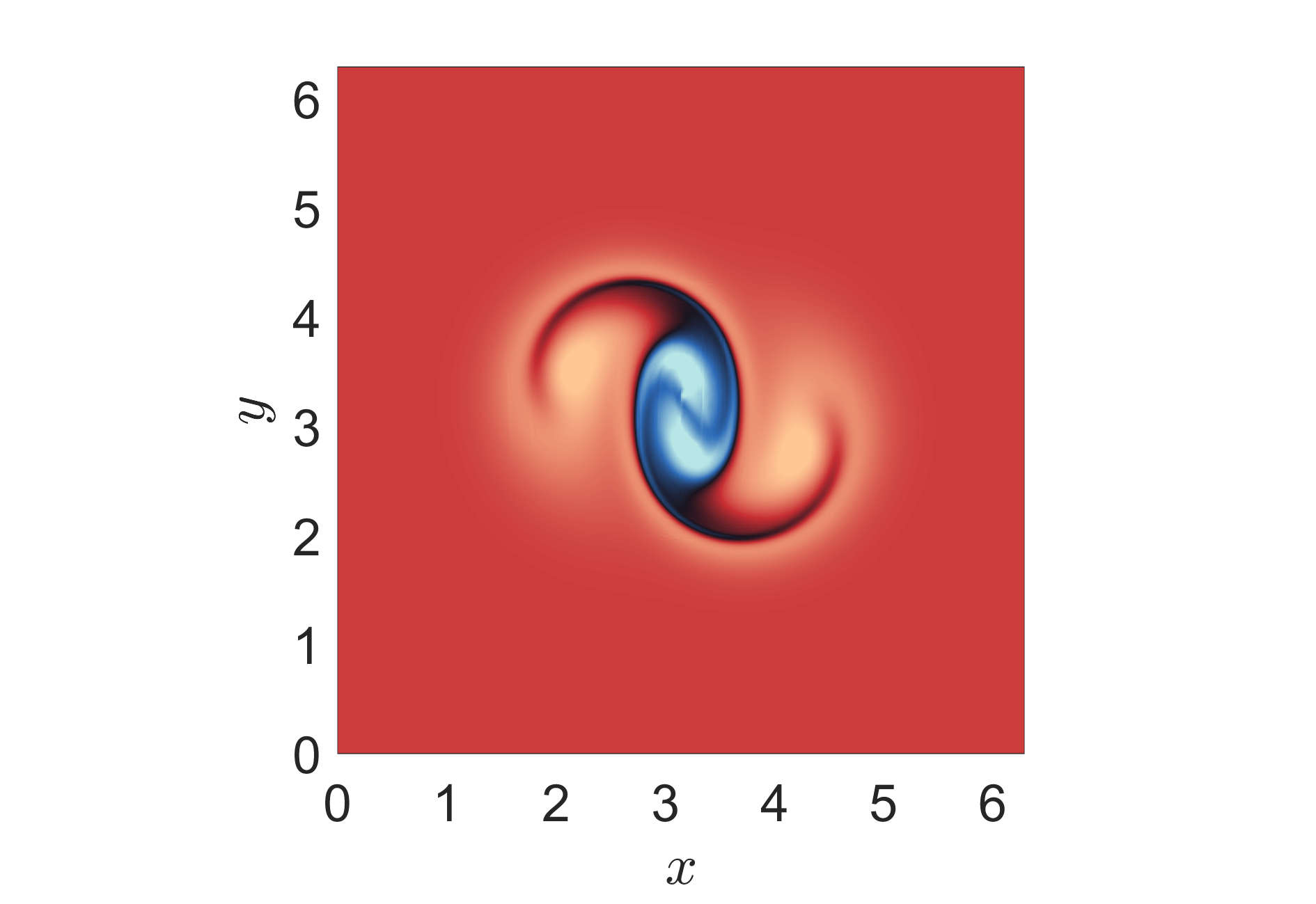}
}\\[-0.5em]
\subfigure[ $t=10$ without $\mathbf{OF}$]{
  \includegraphics[width=0.35\linewidth,trim=3cm 0.1cm 1.5cm 1.5cm, clip]
  {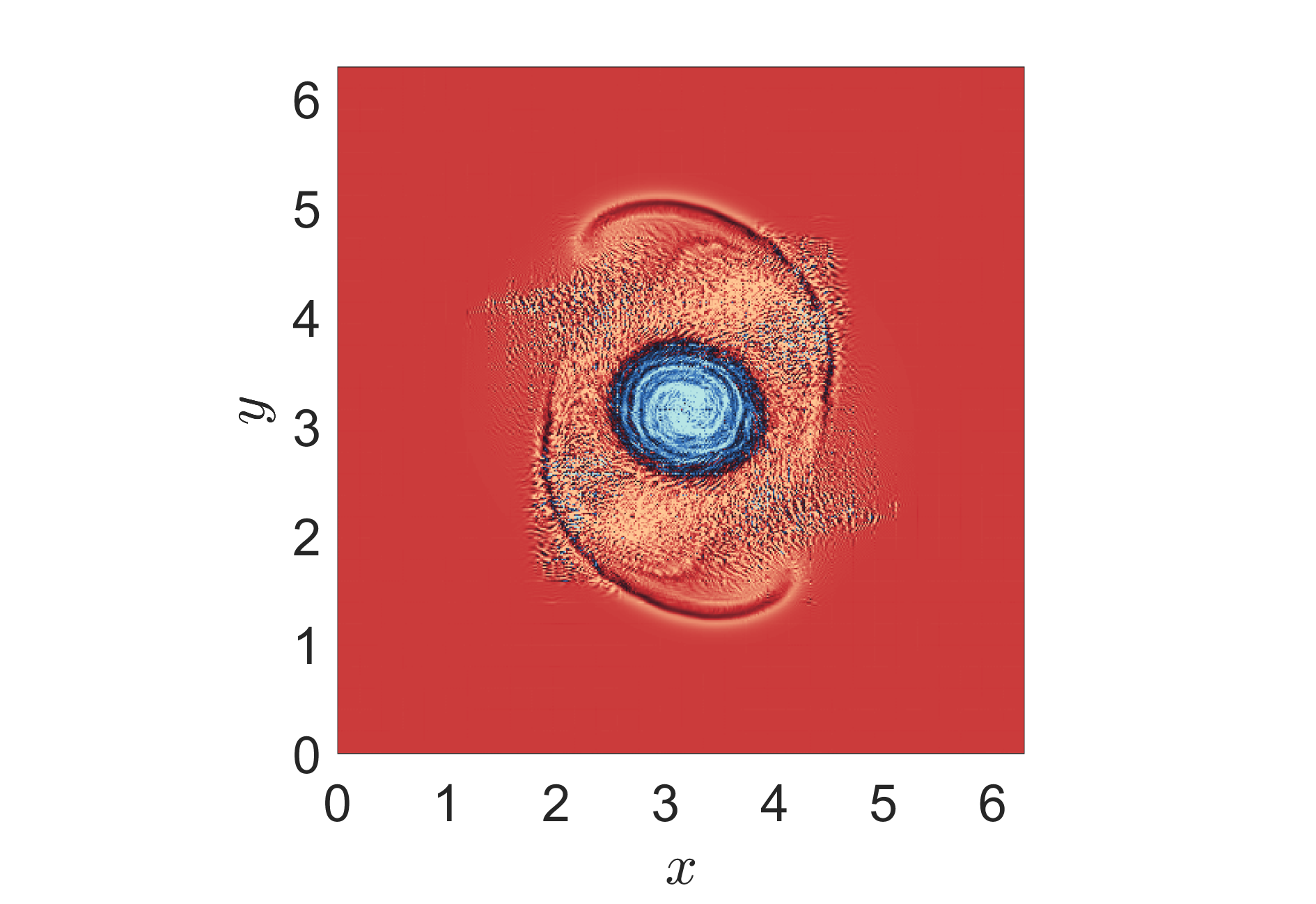}
}
\hspace{-3em}
\subfigure[ $t=10$ with $\mathbf{OF}$]{
  \includegraphics[width=0.35\linewidth,trim=3cm 0.1cm 1.5cm 1.5cm, clip]
  {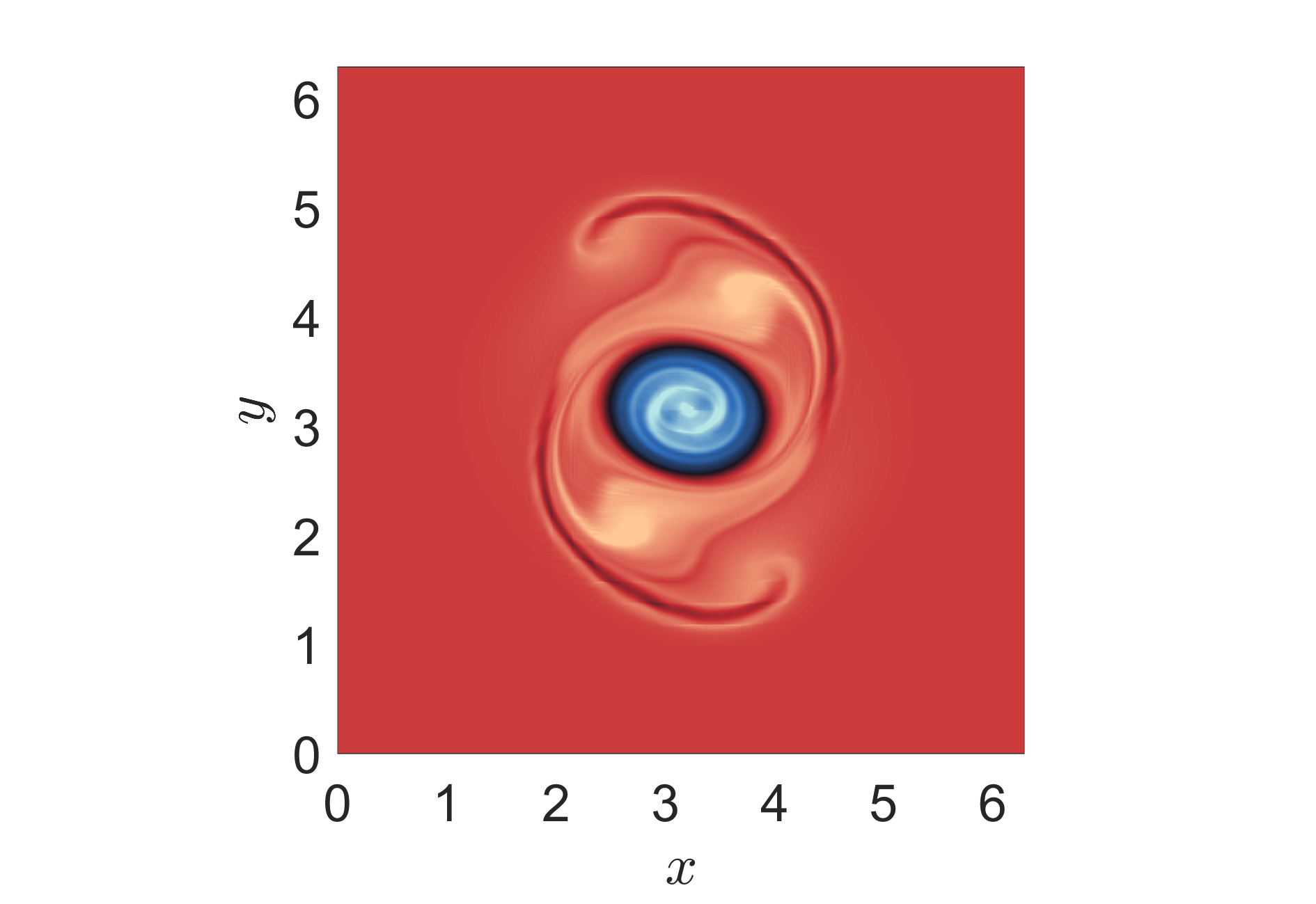}
}
\end{minipage}
\hspace{-8em}
\begin{minipage}[c]{0.08\textwidth}
\centering
\includegraphics[width=\linewidth]{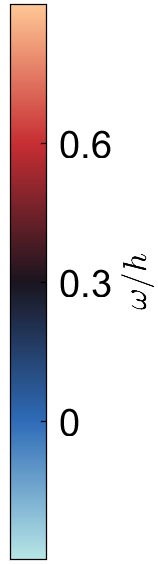}
\end{minipage}
\captionsetup{name={Fig.},labelsep=period}
\caption{
$6$th order for the merging vortices with
$2^5\times2^5$ subdomains in Example~\ref{ex8}.
}\label{sw_2D_merg}
\end{figure}
As shown in Fig.~\ref{sw_2D_merg}, both the SBP-SAT scheme and the oscillation-free SBP-SAT scheme capture the two vortices well and produce similar numerical solutions at the early stage of the simulation, such as at $t=2$. However, as the vortices evolve and interact, spurious oscillations gradually develop in the solution without the damping term, severely degrading the coherent vortex structures. Although the total entropy decreases, these oscillations are accompanied by the nonphysical growth in total enstrophy, as shown in Fig.~\ref{sw_2D_enes1}. In contrast, the oscillation-free SBP-SAT scheme effectively suppresses spurious oscillations and preserves the coherent vortex structures, while both total entropy and enstrophy decrease over time. These results demonstrate that the decay of total entropy alone does not prevent spurious oscillations and verify the effectiveness of the damping term in controlling enstrophy growth and preserving vortex structures during long-time simulations.
\end{example}

\begin{example}[Barotropic shear instability]\label{ex9}
We follow reference~\cite{Hew2025} and finally consider the barotropic shear instability for the two-dimensional rotating shallow water equations, 
This test describes the development of instability in a zonal shear flow triggered by small perturbations. The computational domain is $\Omega=[0,L]^2$ with $L=2\pi\times 6371.22~\mathrm{km}$.
The initial velocity field is prescribed as
\[
U_0(y)=\bar{U}_0[\operatorname{sech}(\frac{y-y_+}{\delta})-\operatorname{sech}(\frac{y-y_-}{\delta})], \quad  V_0=0,
\]
where $y_+=0.25L$, $y_-=0.75L$ and $\delta=10^6~\mathrm{m}$.

To maintain geostrophic balance in the initial state, the water height is defined by
\[
h_0(x,y)=H-\frac{f}{g}\int_0^y U_0(s)\,ds+\widetilde{h}(x,y),
\]
where
$\widetilde{h}(x,y)=
\bar{h}_0\sum_{i=1}^{2}\exp\left(-k d_i(x,y)\right)$, $
d_i(x,y)=\frac{(x-x_i)^2+(y-y_i)^2}{L^2}$, $ i=1,2$. 
The perturbation centers are $(x_1,y_1)=(0.15L,y_+)$
and $(x_2,y_2)=(0.85L,y_-)$. The physical parameters are
$\bar{U}_0=50~\mathrm{m\,s^{-1}}$, $f=7.292\times10^{-5}~\mathrm{s^{-1}}$,  $g=9.80616~\mathrm{m\,s^{-2}}$, $H=10~\mathrm{km}$, $\bar{h}_0=0.01H$ and$ k=10^3$.
\begin{figure}[!htb]  
  \centering
  \subfigure {
    \includegraphics[width=0.4\textwidth] {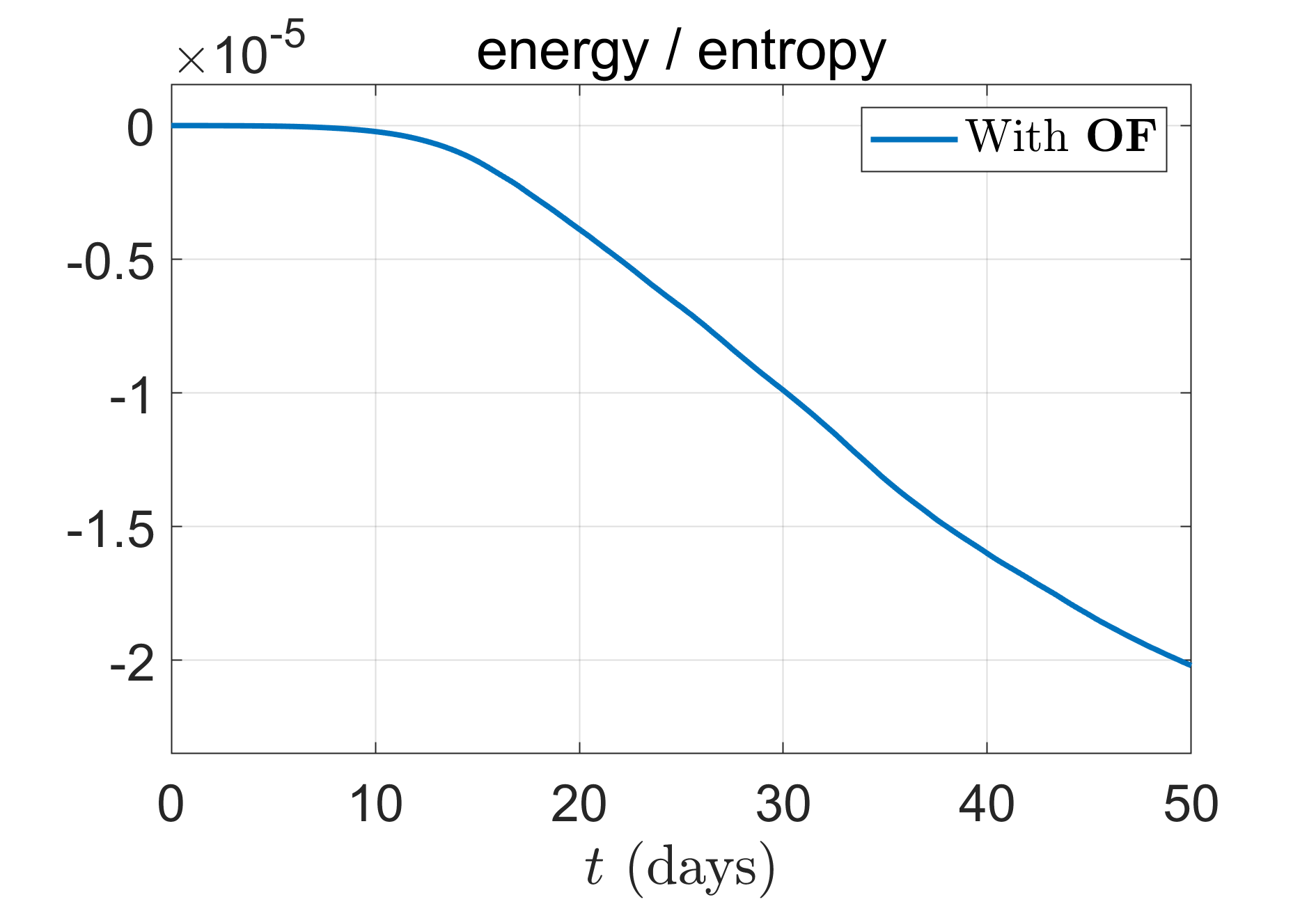}
  }
  \hspace{-2em}
  \subfigure{
    \includegraphics[width=0.4\textwidth]
    {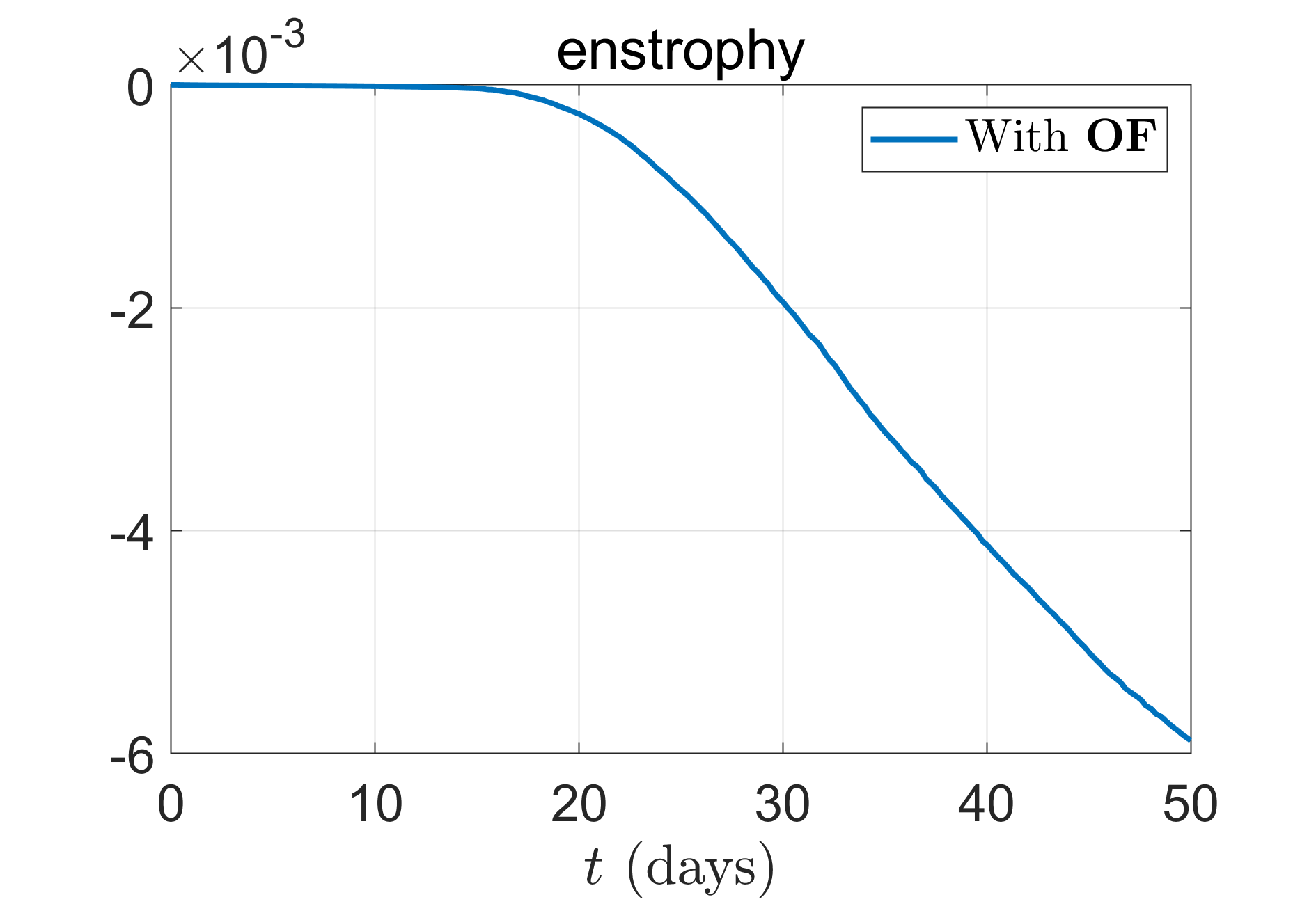}
  }
  \captionsetup{name={Fig.},labelsep=period}
  \caption{The relative change of Entropy and Enstrophy in Example~\ref{ex9}.}   
  \label{sw_2D_enes2}
\end{figure}

\begin{figure}[!htb]  
\centering
\begin{minipage}[c]{0.82\textwidth}
\centering
\subfigure[ $t=9$ with $\mathbf {OF}$]{
  \includegraphics[width=0.4\linewidth, trim=0cm 0.75cm 1cm 2cm, clip]
  {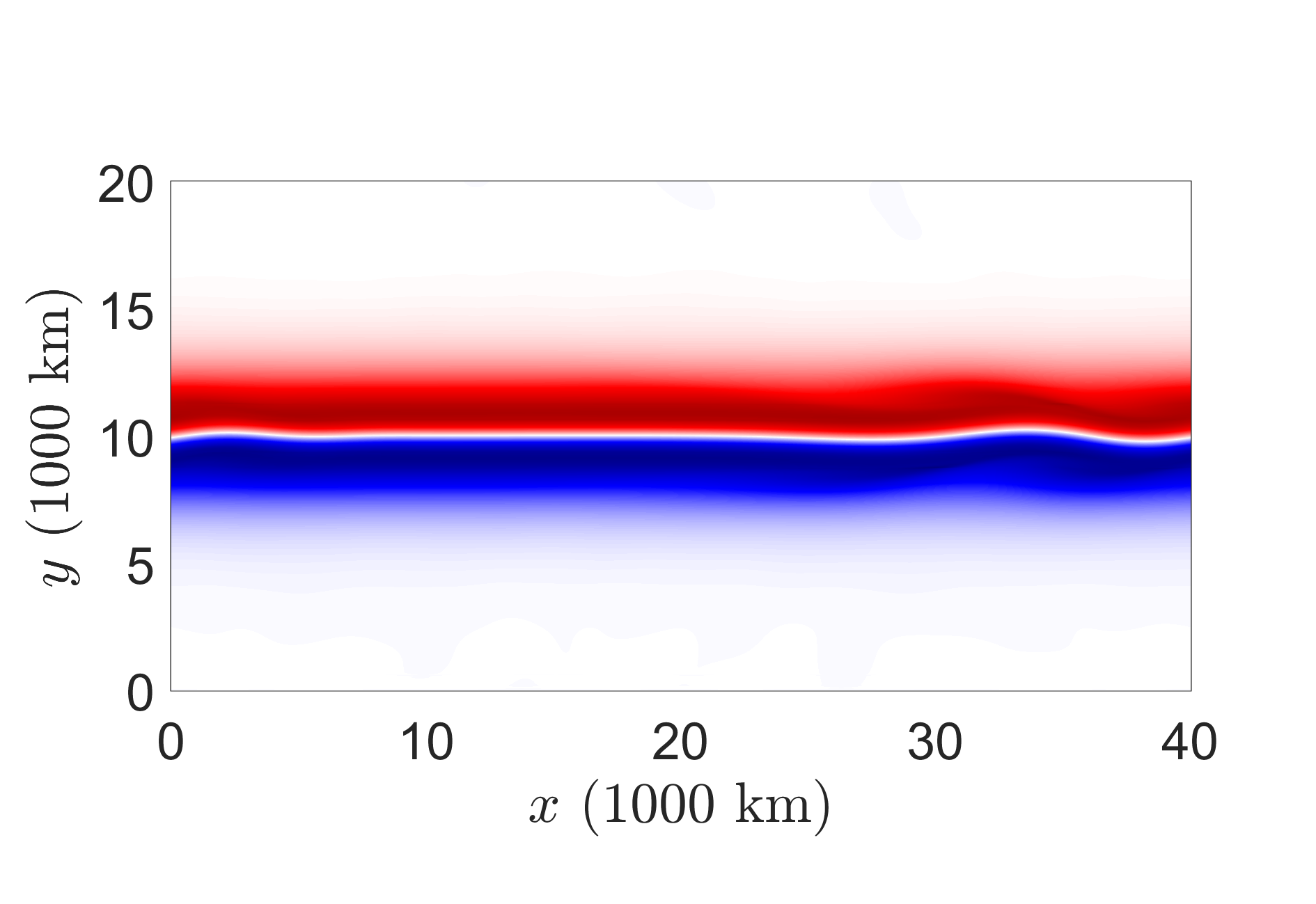}
}
\hspace{-0.2em}
\subfigure[ $t=20$ with $\mathbf {OF}$]{
  \includegraphics[width=0.4\linewidth, trim=0cm 2cm 1.5cm 2.3cm, clip]
  {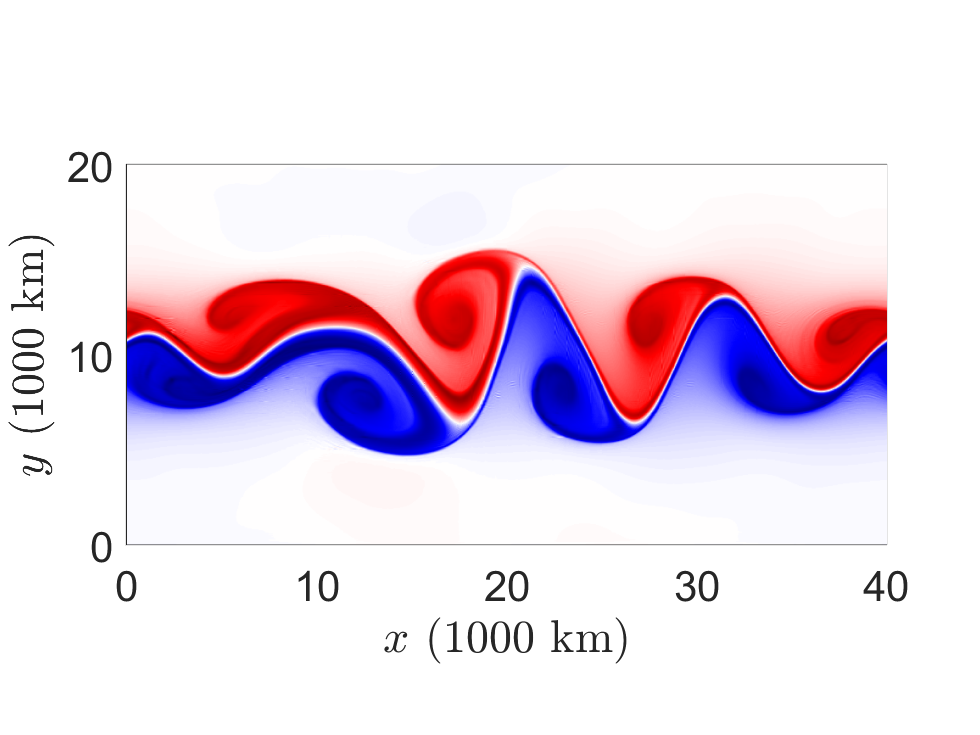}
}\\[-0.3em]
\subfigure[ $t=30$ with $\mathbf {OF}$]{
  \includegraphics[width=0.4\linewidth, trim=0cm 0.75cm 1cm 2cm, clip]
  {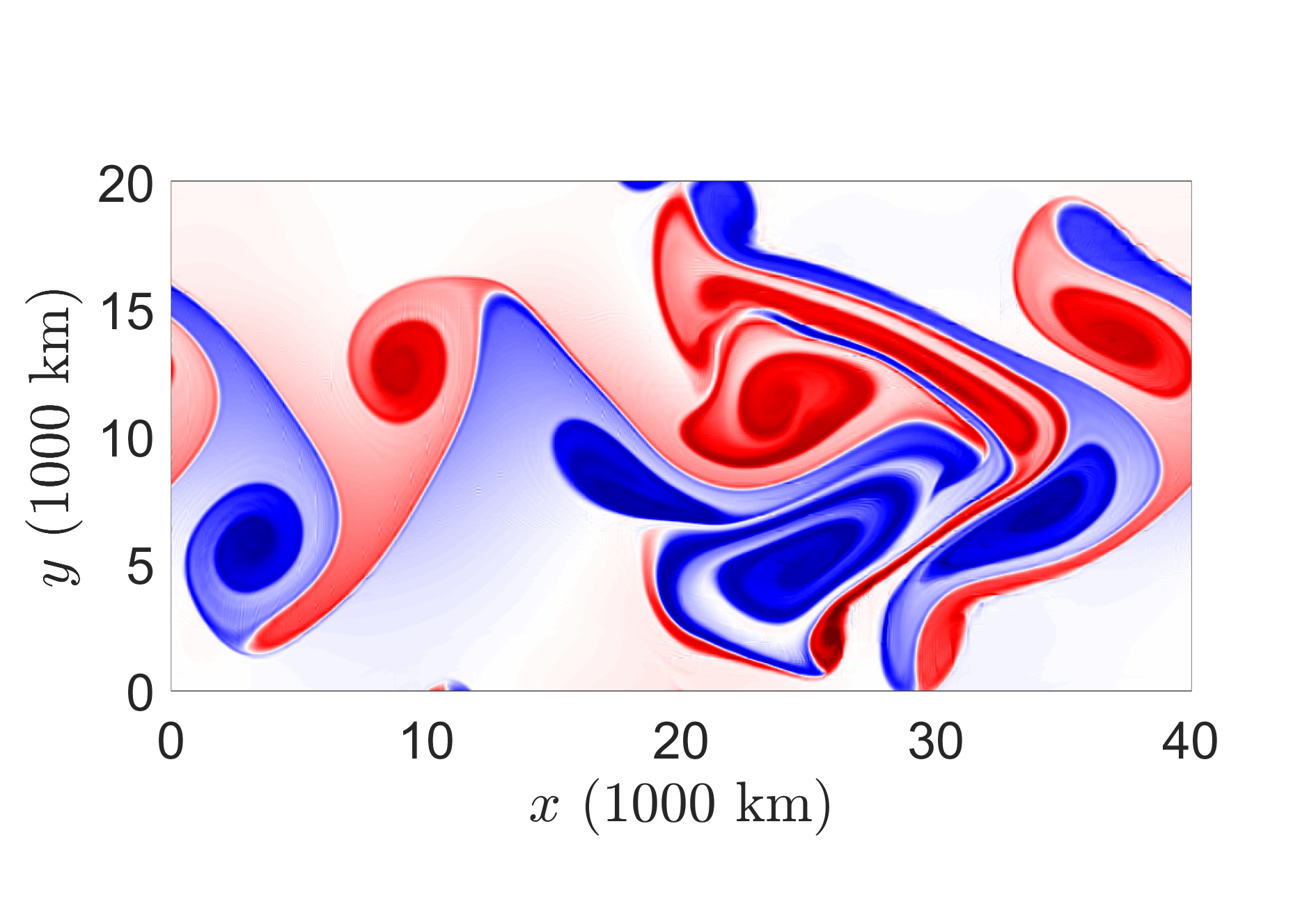}
}
\hspace{-0.2em}
\subfigure[ $t=50$ with $\mathbf {OF}$]{
  \includegraphics[width=0.4\linewidth, trim=0cm 0.8cm 1cm 2cm, clip]
  {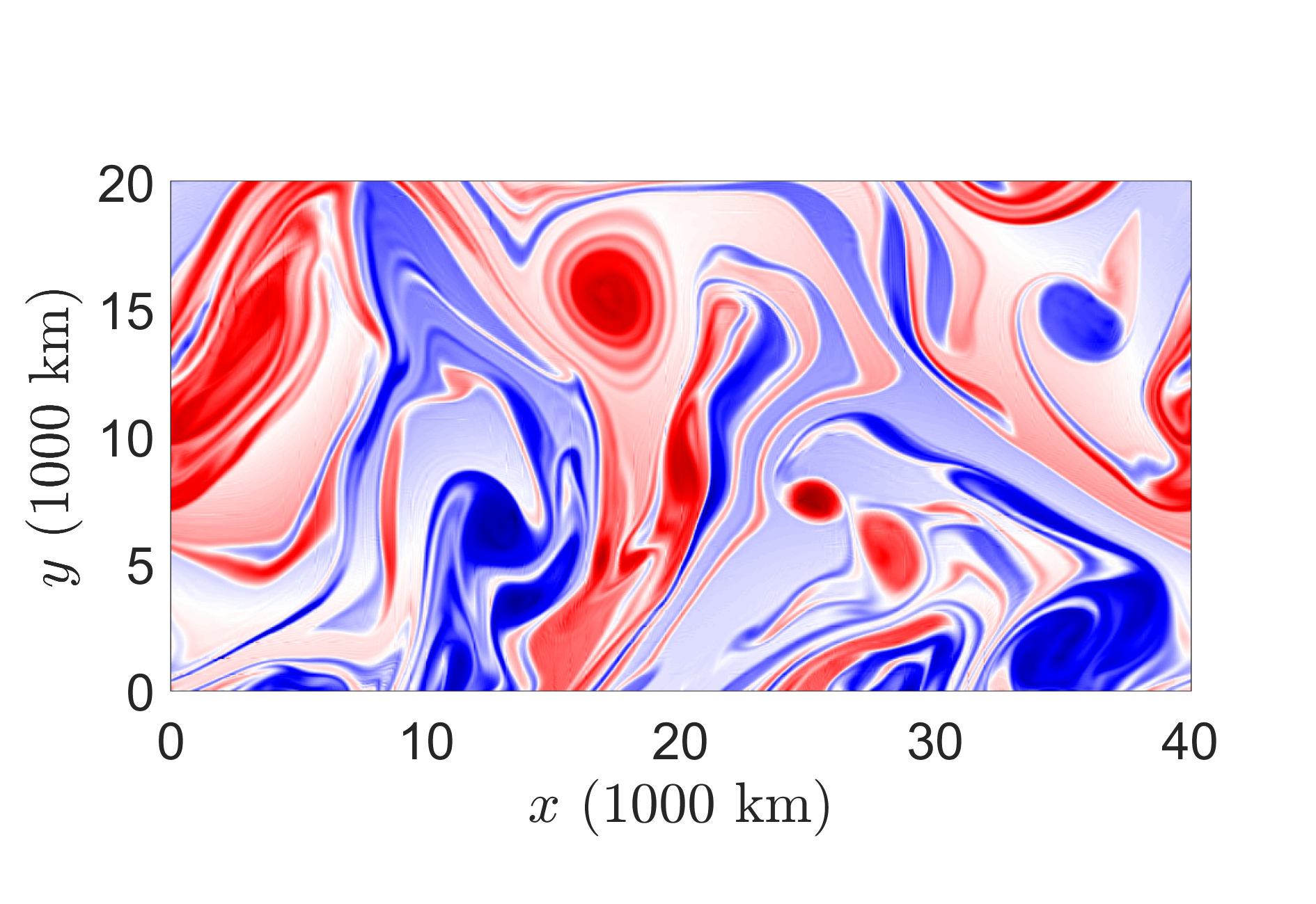}
}
\end{minipage}
\hspace{-3em}
\begin{minipage}[c]{0.1\textwidth}
\centering
\includegraphics[width=\linewidth]{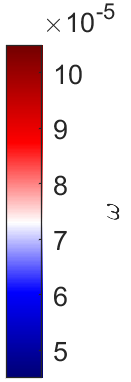}
\end{minipage}
\captionsetup{name={Fig.},labelsep=period}
\caption{
$6$th order for the barotropic shear instability  with $2^6\times2^6$ subdomains in Example~\ref{ex9}.
}\label{sw_2D_bar}
\end{figure}
To illustrate the subsequent formation and evolution of vortical structures, we plot the absolute vorticity $\omega=\partial_xV-\partial_yU+f$ up to $t=50$ days. We use the SBP operator $D_1$ with sixth-order interior accuracy, a time step of $\Delta t=0.002\Delta x$
and $K=2^6\times2^6$ subdomains. Due to the symmetry of the problem, only the lower half of the computational domain is displayed. In addition, we monitor the relative changes in the total entropy and total enstrophy during the simulation, where the total enstrophy is defined as in Example~\ref{ex8}.

We first perform the simulation without the damping term. In this case, the numerical solution becomes unstable and the computation breaks down at approximately $t=15$ days, before the vortical structures have fully developed. In contrast, the oscillation-free SBP-SAT scheme remains stable up to the final time $t=50$ days, capturing the formation and subsequent evolution of complex vortical structures, as shown in Fig.~\ref{sw_2D_bar}. Moreover, both the total entropy and enstrophy remain non-increasing throughout the simulation, as shown in Fig.~\ref{sw_2D_enes2}. These results demonstrate the effectiveness of the damping term in preventing numerical breakdown and enabling the long-time simulation of barotropic shear instability.
\end{example}

\section{Conclusion}
In this work, we have developed an oscillation-free SBP-SAT finite difference method for hyperbolic conservation laws within a multi-domain framework. We first present the underlying SBP-SAT semi-discrete formulation for a one-dimensional system with $m$ conserved variables on two subdomains and show that it satisfies discrete conservation and entropy stability. A projection-based damping term is then incorporated into the SBP-SAT framework for both one- and two-dimensional systems. We prove that the oscillation-free schemes preserve discrete conservation and entropy stability. Furthermore, the accuracy analysis shows that the damping term does not degrade the formal high-order accuracy of the scheme for smooth solutions.

Numerical experiments for the linear advection equation, Burgers' equation, and the shallow water equations demonstrate that the proposed method effectively suppresses spurious oscillations near discontinuities and steep gradients while retaining high-order accuracy in smooth regions. In addition, two physical problems involving the rotating shallow water equations further demonstrate the robustness of the proposed method in long-time simulations.

Future work will focus on extending the proposed method to more complex physical problems and investigating its construction and analysis under non-periodic boundary conditions.

\appendix
\section{SBP-SAT discretizations for the model problems}\label{app1}
In this appendix, we present the semi-discrete formulations of Eq.~\eqref{generalform_nondp} for the linear advection equation, Burgers' equation, and the shallow water equations. To determine the SAT coefficients through the semi-discrete entropy estimates, we first set the upwind coefficient $\boldsymbol{\alpha}=0$ and focus on the SAT terms. Since the upwind interface term contributes a non-positive quantity to the entropy estimates, its inclusion does not affect the entropy stability of the resulting schemes. The SAT coefficients in Appendices~\ref{app.adv}, \ref{app.bur}, and \ref{app.sw} are chosen to satisfy the corresponding semi-discrete entropy estimates. Since the choice of SAT terms is not unique, the formulations presented here represent one possible choice. In particular, the SAT terms for Burgers' equation and the shallow water equations are chosen following~\cite{Ranocha2017h} and~\cite{Ranocha2017}, respectively. The upwind coefficients used in the numerical experiments are specified in Section~\ref{Sec.5}.
 \subsection{Advection equation}\label{app.adv}
We begin with the advection equation, which serves as the simplest linear scalar conservation laws in our discussion.
For the one-dimensional advection equation \eqref{adv} with advection speed $a=1$ , the two-subdomain continuous problem under periodic boundary conditions takes the form
\begin{equation*}
\begin{cases}
U_t^{(1)}+U_x^{(1)}=0, 
& x\in \Omega_1,\\
U_t^{(2)} +U_x^{(2)}=0, 
& x\in \Omega_2,\\
U^{(1)}(b^{(1)},t)=U^{(2)}(a^{(2)},t),\\
U^{(1)}(a^{(1)},t)=U^{(2)}(b^{(2)},t).
\end{cases}
\end{equation*}
The semi-discrete SBP-SAT formulation for the two-subdomain advection equation is then given by
\begin{equation}\label{adv_if_semi}
\begin{cases}
\mathbf{u}^{(1)}_t + D_1\mathbf{u}^{(1)}= 
H^{-1}\left[
\dfrac12\boldsymbol\delta_0\left(u^{(2)}_N-u^{(1)}_0\right)
-\dfrac12\boldsymbol\delta_{N}\left(u^{(2)}_0-u^{(1)}_N\right)
\right],
\\
\mathbf{u}^{(2)}_t +D_1\mathbf{u}^{(2)}= 
H^{-1}\left[
-\dfrac12\boldsymbol\delta_{N}\left(u^{(1)}_0-u^{(2)}_N\right)
+\dfrac12\boldsymbol\delta_0\left(u^{(1)}_N-u^{(2)}_0\right)
\right].
\end{cases}
\end{equation}
The right-hand side of the above equations are SAT terms, where the first and second terms impose the periodic boundary condition and the interface condition, respectively.

 \subsection{Burgers' equation}\label{app.bur}
For the nonlinear scalar conservation laws, we consider the one-dimensional inviscid Burgers' equation introduced in Example~\ref{bur}. 
In the continuous energy estimate, the chain rule is used to cancel some positive terms and then get that the energy is non-increasing in time. However, this chain rule does not hold at the discrete level. Therefore, to mimic the continuous energy estimate in the discrete setting, the problem is written in the following split form
\begin{equation*}
\begin{cases}
U_t^{(1)}+ \dfrac13 \left(U^{(1)}U_x^{(1)}+
((U^{(1)})^2)_x\right)=0, & x\in\Omega_1,\\
U_t^{(2)}+ \dfrac13 \left(U^{(2)}U_x^{(2)}+
((U^{(2)})^2)_x\right)=0,& x\in\Omega_2,\\
U^{(1)}(b^{(1)},t)=U^{(2)}(a^{(2)},t),\\
U^{(1)}(a^{(1)},t)=U^{(2)}(b^{(2)},t).
\end{cases}
\end{equation*}
The associated semi-discrete SBP-SAT scheme is as follows
\begin{equation}\label{bur_if_semi}
\begin{cases}
\mathbf{u}^{(1)}_t+\dfrac{1}{3}\left(
\mathbf{u}^{(1)} \circ D_1 \mathbf{u}^{(1)}+ D_1((\mathbf{u}^{(1)})^2)\right)= H^{-1}\,\mathbf{SAT}^{(1)},\\
\mathbf{u}^{(2)}_t + \dfrac{1}{3}\left(\mathbf{u}^{(2)} \circ D_1 \mathbf{u}^{(2)} + D_1((\mathbf{u}^{(2)})^2)\right)= H^{-1}\,\mathbf{SAT}^{(2)},
\end{cases}
\end{equation}
  where $\mathbf{SAT}^{(1)}$ and $\mathbf{SAT}^{(2)}$ denote the SAT terms. Following the construction in \cite{Ranocha2017h}, they are chosen as,
 \begin{align*}
\mathbf{SAT}^{(1)}
=&-\dfrac16\boldsymbol\delta_{N}
\left(u^{(1)}_N(u^{(2)}_0-u^{(1)}_N)\right)
-\dfrac13\boldsymbol\delta_{N}
\left(\dfrac{(u^{(2)}_0)^2-(u^{(1)}_N)^2}{2}\right)\\
&+\dfrac16\boldsymbol\delta_0 \left(u^{(1)}_0(u^{(2)}_N-u^{(1)}_0)\right) +\dfrac13\boldsymbol\delta_0\left(\dfrac{(u^{(2)}_N)^2-(u^{(1)}_0)^2}{2}\right),\\
\mathbf{SAT}^{(2)}=&\dfrac16\boldsymbol\delta_0
\left(u^{(2)}_0(u^{(1)}_N-u^{(2)}_0)\right)
+\dfrac13\boldsymbol\delta_0\left(\dfrac{(u^{(1)}_N)^2-(u^{(2)}_0)^2}{2}\right)\\
&-\dfrac16\boldsymbol\delta_{N}\left(u^{(2)}_N(u^{(1)}_0-u^{(2)}_N)\right)-\dfrac13\boldsymbol\delta_{N}
\left(\dfrac{(u^{(1)}_0)^2-(u^{(2)}_N)^2}{2}\right).
\end{align*}
Here, the first two contributions in both $\mathbf{SAT}^{(1)}$ and $\mathbf{SAT}^{(2)}$ correspond to the interface coupling, while the last two correspond to the periodic boundary condition. Compared with the linear advection case, the form of the SAT terms is no longer linear in the boundary and interface values. Instead, it is written in the nonlinear split form to obtain that the discrete energy is non-increasing in time. 

\subsection{Shallow water equation}  \label{app.sw}
We further consider the shallow water equations in Example~\ref{sw} as a nonlinear system of conservation laws. This model involves a vector of conserved variables, let $ \mathbf{U}^{(k)}=[h^{(k)}, h^{(k)}U^{(k)}]^\top$ denote the conserved variable vector on $\Omega_k$, $k=1,2$. 
As in the Burgers case, the corresponding subdomain equations are written in the split form,
\[
\left\{
\begin{aligned}
\partial_t h^{(k)} &+\partial_x (h^{(k)}U^{(k)})=0,\\
\partial_t (h^{(k)}U^{(k)}) &+\dfrac12\partial_x\left(h^{(k)}(U^{(k)})^2\right)+\dfrac12U^{(k)}\partial_x (h^{(k)}U^{(k)})\\
&+\dfrac12(h^{(k)}U^{(k)})\partial_xU^{(k)}
+gh^{(k)}\partial_xh^{(k)}=0,
\end{aligned}
\right.
\]
with the interface condition 
$\mathbf U^{(1)}(b^{(1)},t)=\mathbf U^{(2)}(a^{(2)},t)$ and the periodic boundary condition $\mathbf U^{(1)}(a^{(1)},t)= \mathbf U^{(2)}(b^{(2)},t)$.

For $k=1,2$, the corresponding semi-discrete SBP-SAT scheme takes the form
\small{\begin{equation}\label{sw_if_semi}
\left\{\begin{aligned}
\partial_t \mathbf h^{(k)}+&D_1 (\mathbf h^{(k)}\circ\mathbf u^{(k)})
=H^{-1}\mathbf{SAT}_h^{(k)},\\
\partial_t (\mathbf h^{(k)}\circ\mathbf u^{(k)})+&\dfrac12D_1\left(
\mathbf h^{(k)}\circ(\mathbf{u}^{(k)})^2\right)
+\dfrac12\mathbf u^{(k)}\circ D_1(\mathbf h^{(k)}\circ\mathbf u^{(k)})\\
+&\dfrac12(\mathbf h^{(k)}\circ\mathbf u^{(k)}) \circ D_1\mathbf u^{(k)}
+g\mathbf h^{(k)}\circ D_1\mathbf h^{(k)}
=H^{-1}\mathbf{SAT}_{hu}^{(k)},
\end{aligned}\right. 
\end{equation}}
where $\mathbf{SAT}_h^{(k)}$ and $\mathbf{SAT}_{hu}^{(k)}$ denote the corresponding SAT terms for different conservative variables. Following the construction in \cite{Ranocha2017}, the $\mathbf{SAT}_h$ terms are given by
\begin{align*}
\mathbf{SAT}_h^{(1)}
=&-\dfrac12\boldsymbol\delta_{N}\left(h_0^{(2)}u_0^{(2)}-h_N^{(1)}u_N^{(1)}\right)
+\dfrac12\boldsymbol\delta_0\left(h_N^{(2)}u_N^{(2)}-h_0^{(1)}u_0^{(1)}\right),\\
\mathbf{SAT}_h^{(2)}
=&\dfrac12\boldsymbol\delta_0\left(h_N^{(1)}u_N^{(1)}-h_0^{(2)}u_0^{(2)}\right)
-\dfrac12\boldsymbol\delta_{N}\left(h_0^{(1)}u_0^{(1)}-h_N^{(2)}u_N^{(2)}\right).
\end{align*}
In both $\mathbf{SAT}_h^{(1)}$ and $\mathbf{SAT}_h^{(2)}$, the first term on the right-hand side is used to weakly impose the interface condition and the second term is used to weakly impose the periodic boundary condition.\\
For the $\mathbf{SAT}_{hu}$ terms, denote $G_i^{(k)}=h_i^{(k)}(u_i^{(k)})^2$, $i=0,1,\ldots,N$, $ k=1,2$. Then the SAT terms are given by
\begin{align*}
\mathbf{SAT}_{hu}^{(1)}
=&-\boldsymbol\delta_{N}\Big[
\dfrac14\left(G_0^{(2)}-G_N^{(1)}\right)
+\dfrac14 u_N^{(1)}\left(h_0^{(2)}u_0^{(2)}-h_N^{(1)}u_N^{(1)}\right)
+\dfrac14 h_N^{(1)}u_N^{(1)}\left(u_0^{(2)}-u_N^{(1)}\right) \\
& +\dfrac12 g h_N^{(1)}\left(h_0^{(2)}-h_N^{(1)}\right)
\Big] +\boldsymbol \delta_0\Big[
\dfrac14\left(G_N^{(2)}-G_0^{(1)}\right)
+\dfrac14 u_0^{(1)}\left(h_N^{(2)}u_N^{(2)}-h_0^{(1)}u_0^{(1)}\right)\\
&+\dfrac14 h_0^{(1)}u_0^{(1)}\left(u_N^{(2)}-u_0^{(1)}\right)
+\dfrac12g h_0^{(1)}\left(h_N^{(2)}-h_0^{(1)}\right)
\Big],
\end{align*}
and
\begin{align*}
\mathbf{SAT}_{hu}^{(2)}
=&\boldsymbol \delta_0\Big[
\dfrac14\left(G_N^{(1)}-G_0^{(2)}\right)
+\dfrac14 u_0^{(2)}\left(h_N^{(1)}u_N^{(1)}-h_0^{(2)}u_0^{(2)}\right)
+\dfrac14 h_0^{(2)}u_0^{(2)}\left(u_N^{(1)}-u_0^{(2)}\right)\\
&+\dfrac12 g h_0^{(2)}\left(h_N^{(1)}-h_0^{(2)}\right)\Big] -\boldsymbol \delta_{N}\Big[\dfrac14 \left(G_0^{(1)}-G_N^{(2)}\right)
+\dfrac14u_N^{(2)}\left(h_0^{(1)}u_0^{(1)}-h_N^{(2)}u_N^{(2)}\right)\\
&+\dfrac14h_N^{(2)}u_N^{(2)}\left(u_0^{(1)}-u_N^{(2)}\right)+\dfrac12 g h_N^{(2)}\left(h_0^{(1)}-h_N^{(2)}\right)\Big].
\end{align*}
Similarly, on the right-hand sides of both $\mathbf{SAT}_{hu}^{(1)}$ and $\mathbf{SAT}_{hu}^{(2)}$, the terms inside the first bracket weakly impose the interface condition and those inside the second bracket weakly impose the periodic boundary condition.

\section{Proofs of some lemmas and theorems}

\subsection{Proof of Lemma~\ref{lemma_OF_acc}}\label{app_OF_acc}

\begin{proof}

Let $\hat{x}^{(k)}\in\Omega_k$ be an arbitrary point. Based on the Taylor expansion of $g_{\mu}(\mathbf U)$ about point $\hat{x}^{(k)}$, the value of $g_{\mu}(\mathbf U)$ at each grid point $x_i^{(k)}$ can be written as
\begin{align}\label{Tylor}
    g_{\mu}(\mathbf U(x_i^{(k)})) =\sum_{\beta=0}^{l-1}\frac{\partial_x^{\beta} g_{\mu}(\mathbf U(\hat{x}^{(k)}))}{\beta!}(x_i^{(k)}-\hat{x}^{(k)})^{\beta} +
    \frac{\partial_x^l g_{\mu}(\mathbf U(\eta_i^{(k)}))}{l!}
    (x_i^{(k)}-\hat{x}^{(k)})^l,
\end{align}
where $\eta_i^{(k)}$ lies between $\hat{x}^{(k)}$ and $x_i^{(k)}$, $i=0,1,\cdots,N$.

Take $q^{(k)}_{\mu,i} =\sum_{\beta=0}^{l-1}\frac{\partial_x^{\beta} g_{\mu}(\mathbf U(\hat{x}^{(k)}))}{\beta!} (x_i^{(k)}-\hat{x}^{(k)})^{\beta}$ and $ \mathbf q_{\mu}^{(k)}=[q^{(k)}_{\mu,0},q^{(k)}_{\mu,1},\ldots,q^{(k)}_{\mu,N}]^\top$.
Since the corresponding polynomial has degree no greater than $l-1$, $\mathbf q_{\mu}^{(k)}\in V_{l-1}$. Then, by $\left\langle\mathbf g_{\mu}^{(k)} - P^{\,l-1}\mathbf g_{\mu}^{(k)}, P^{\,l-1}\mathbf g_{\mu}^{(k)} -\mathbf q_{\mu}^{(k)}\right\rangle_H =0$, we obtain
\begin{align*}
\left\| \mathbf g_{\mu}^{(k)}-\mathbf q_{\mu}^{(k)}\right\|_H^2=
\left\| \mathbf g_{\mu}^{(k)}-P^{\,l-1}\mathbf g_{\mu}^{(k)}
\right\|_H^2+\left\|P^{\,l-1}\mathbf g_{\mu}^{(k)}
-\mathbf q_{\mu}^{(k)}\right\|_H^2.
\end{align*}
Hence,
\begin{align}\label{EE_Pytha}
    \left\|\mathbf g_{\mu}^{(k)}- P^{\,l-1}\mathbf g_{\mu}^{(k)} \right\|_H\le
    \left\| \mathbf g_{\mu}^{(k)}-\mathbf q_{\mu}^{(k)} \right\|_H.
\end{align}
From equation~\eqref{Tylor}, we get
\begin{align*}
    g_{\mu}(\mathbf U(x_i^{(k)}))-q^{(k)}_{\mu,i} = \frac{\partial_x^l g_{\mu}(\mathbf U(\eta_i^{(k)}))}{l!}
    (x_i^{(k)}-\hat{x}^{(k)})^l, \quad  0\le i\le N
\end{align*}
By $|x_i^{(k)}-\hat{x}^{(k)}| \le |\Omega_k|$ and $ \left|\partial_x^l g_{\mu}(\mathbf U(\eta_i^{(k)}))\right| \le\max_{x\in \Omega_k}\left| \partial_x^l g_{\mu}(\mathbf U)\right|$, we get
\begin{align*}
   \left|g_{\mu}(\mathbf U(x_i^{(k)}))-q^{(k)}_{\mu,i}\right|\le C |\Omega_k|^l.
\end{align*}
Since the diagonal-norm SBP matrix $H$ satisfies
$\mathbf 1^\top H\mathbf 1=|\Omega_k|$, it follows that
\begin{align}\label{EE_dp_sub}
    \left\|\mathbf g_{\mu}^{(k)}-\mathbf q_{\mu}^{(k)}\right\|_H^2=\Delta x \sum_{i=0}^{N}
    w_{i}|g_{\mu}(\mathbf U(x_i^{(k)}))-q^{(k)}_{\mu,i}|^2\le C|\Omega_k|^{2l+1}.
\end{align}
By the equation~\eqref{EE_Pytha} and~\eqref{EE_dp_sub}, we yield
\begin{align}\label{pre}
    \left\| \mathbf r_{\mu}^{(k,l)}\right\|_H = \left\|\mathbf g_{\mu}^{(k)}- P^{\,l-1}\mathbf g_{\mu}^{(k)} \right\|_H\le C|\Omega_k|^{l+\frac12},
    \quad l\ge1. 
\end{align}
We next estimate the interface jumps appearing in
$\sigma_{\mu}^{(k,l)}$.  The vectors
$\mathbf d_{\mathcal I_L}^{(k,l)}$ and
$\mathbf d_{\mathcal I_R}^{(k,l)}$ are constructed as $p$-th order one-sided approximations of the $l$-th derivative.

For the interface jump, we have
\begin{align*}
    \left|[\mkern-2.8mu[
\mathcal D^l \mathbf z_{\mu}^{(k)}
]\mkern-2.8mu]_{R}
\right| =\left|(\mathbf d_{\mathcal I_L}^{(k+1,l)})^\top
    \mathbf z_{\mu}^{(k+1)} - (\mathbf d_{\mathcal I_R}^{(k,l)})^\top\mathbf z_{\mu}^{(k)}\right|
    \le C N^l|\Omega_k|^p.
\end{align*}
Here, for the neighbor subdomains $ \Omega_{k+1}$ and $\Omega_k$, they satisfy $x_L^{(k+1)}=x_R^{(k)}$. Since $\partial_x^lz_{\mu}(\mathbf U)$ is continuous across the interface, $\partial_x^lz_{\mu}(\mathbf U(x_L^{(k+1)})) = \partial_x^lz_{\mu}(\mathbf U(x_R^{(k)}))$.

Similarly, we have $\left|[\mkern-2.8mu[
\mathcal D^l \mathbf z_{\mu}^{(k)}
]\mkern-2.8mu]_{L}
\right|\le C N^l|\Omega_k|^p$.
For $l\ge1$, we have
\begin{align*}
 & \sqrt{[\mkern-2.8mu[
\mathcal D^l \mathbf z_{\mu}^{(k)}
]\mkern-2.8mu]_{L}^2 + [\mkern-2.8mu[
\mathcal D^l \mathbf z_{\mu}^{(k)}
]\mkern-2.8mu]_{R}^2 } \le C N^l|\Omega_k|^p.
\end{align*}
Therefore, from the definition
\begin{align}\label{jump_ee}
    \sigma_{\mu}^{(k,l)} =C_l(\Delta x)^l \sqrt{
    [\mkern-2.8mu[
\mathcal D^l \mathbf z_{\mu}^{(k)}
]\mkern-2.8mu]_{L}^2 + [\mkern-2.8mu[
\mathcal D^l \mathbf z_{\mu}^{(k)}
]\mkern-2.8mu]_{R}^2  }
   \le C |\Omega_k|^{l+p}.
\end{align}
We now combine \eqref{pre} and
\eqref{jump_ee}. Since $\sigma_{\mu}^{(k,0)}(\mathbf U)=0$, only the terms with $l\geq1$ contribute to the damping term evaluated at the exact solution.
For each $l\ge1$,
\begin{align*}
    \left\|\frac{\sigma_{\mu}^{(k,l)}}{|\Omega_k|}
    \mathbf r_{\mu}^{(k,l)}\right\|_H
    = \frac{\sigma_{\mu}^{(k,l)}}{|\Omega_k|}
    \left\| \mathbf r_{\mu}^{(k,l)}
    \right\|_H
    \le C|\Omega_k|^{2l+p-\frac12}.
\end{align*}
As $|\Omega_k|$ goes to zero, we have
\begin{align*}
    \left\|\mathbf{OF}_{\mu}^{(k)}(\mathbf U)\right\|_H
    \le C\sum_{l=1}^{s}|\Omega_k|^{2l+p-\frac12}\le C|\Omega_k|^{p+\frac32},\quad  \forall \,\ 1\le \mu\le m.
\end{align*}
Thus, for a sufficiently smooth exact solution, the local damping term is of order $p+\frac32$ in the discrete $L^2$ norm. Consequently, the
additional damping term does not reduce the order accuracy of the existing SBP-SAT discretization.
\end{proof}

\bibliographystyle{elsarticle-num} 
\bibliography{references}

\end{document}